\documentclass[graybox]{svmult}
\usepackage{mathptmx}
\usepackage{helvet}
\usepackage{courier}
\usepackage{type1cm}
\usepackage{makeidx}
\usepackage{graphicx}
\usepackage{multicol}
\usepackage[bottom]{footmisc}

\makeindex

\input{tcilatex}

\begin{document}

\title*{Parameterized Yang-Hilbert-Type Integral Inequalities and Their
Operator Expressions}
\titlerunning{Parameterized Yang-Hilbert-Type Integral Inequalities and Their Operator Expressions}
\author{Bicheng Yang\ and\ Michael Th. Rassias}
\institute{Bicheng Yang, Department of Mathematics, Guangdong University of
 Education, Guangzhou, Guangdong
 510303, P. R. China, \\
 e-mail: {bcyang@gdei.edu.cn\,\,\,\,bcyang818@163.com} \\
  Michael Th. Rassias, Department of Mathematics, Princeton University, Fine Hall, Washington Road, Princeton, NJ 08544-1000, USA\ \&  Department of Mathematics, ETH-Z\"urich, CH-8092, Switzerland,\\ %
e-mail:{ michailrassias@math.princeton.edu; michael.rassias@math.ethz.ch}}
%
%
\maketitle

\begin{abstract}
{ Applying methods of Real Analysis and Functional Analysis, we build two
weight functions with parameters and provide two kinds of parameterized
Yang-Hilbert-type integral inequalities with the best constant factors. Equivalent
forms, the reverses, and the operator expressions are also given. In
particular, the Hardy-type inequalities and Hardy-type operators are
studied. Additionally, a number of examples with two kinds of particular kernels are
considered.}
\end{abstract}

\section{Introduction}

\label{sec:1} If $f(x),\:g(y)\geq 0,$ satisfy $$0<\int_{0}^{\infty
}f^{2}(x)dx<\infty $$ and $$0<\int_{0}^{\infty }g^{2}(y)dy<\infty ,$$ then we
have the following well known Hilbert's integral inequality (cf. \cite{YM1})
\begin{equation}
\int_{0}^{\infty }\int_{0}^{\infty }\frac{f(x)g(y)}{x+y}dxdy<\pi \left(
\int_{0}^{\infty }f^{2}(x)dx\int_{0}^{\infty }g^{2}(y)dy\right) ^{\frac{1}{2}%
},  \label{1}
\end{equation}%
where the constant factor $\pi $ is the best possible. The operator
expression of (\ref{1}) was studied in \cite{YM1a} and \cite%
{YM1b}.

In 1925, by introducing one pair of conjugate exponents $(p,q)$, that is $\frac{1}{p}+%
\frac{1}{q}=1$, \ Hardy \cite{YM2} provided an extension of (\ref{1}) as
follows:

For $ p>1,$ $f(x),g(y)\geq 0,$ satisfying $$0<\int_{0}^{\infty
}f^{p}(x)dx<\infty $$ and $$0<\int_{0}^{\infty }g^{q}(y)dy<\infty ,$$ we have%
\begin{equation}
\int_{0}^{\infty }\int_{0}^{\infty }\frac{f(x)g(y)}{x+y}dxdy <\frac{\pi }{%
\sin (\pi /p)}\left( \int_{0}^{\infty }f^{p}(x)dx\right) ^{\frac{1}{p}%
}\left( \int_{0}^{\infty }g^{q}(y)dy\right) ^{\frac{1}{q}},  \label{2}
\end{equation}%
where the constant factor $\frac{\pi }{\sin (\pi /p)}$ is still the best
possible. Inequality (\ref{2}) is known as Hardy-Hilbert's integral
inequality, and is important in analysis and its applications (cf. \cite%
{YM3}, \cite{YM4}).

\begin{definition}
If $\lambda \in \mathbf{R=(-\infty ,\infty )}, k_{\lambda }(x,y)$ is a non-negative measurable function in \textbf{\ }$ \mathbf{R}_{+}^{2}=\mathbf{R}_{+}\times \mathbf{R}_{+}$, satisfying $$k_{\lambda }(tx,ty)=t^{-\lambda }k_{\lambda }(x,y),$$
for any $ t,x,y>0$, then we call $k_{\lambda }(x,y)$ the homogeneous function of degree $-\lambda $ in $\mathbf{R}_{+}^{2}$.
\end{definition}

In 1934, replacing $\frac{1}{x+y}$ in (\ref{2}) by a general homogeneous kernel of degree
-1, as $k_{1}(x,y)$, Hardy et al. presented an extension of (\ref{2})
with the best possible constant factor
\[
k_{p}=\int_{0}^{\infty }k_{1}(t,1)t^{\frac{-1}{p}}dt\in \mathbf{R}%
_{+}=(0,\infty )
\]
obtaining (cf. \cite{YM3}, Th. 319):%
\begin{equation}
\int_{0}^{\infty }\int_{0}^{\infty }k_{1}(x,y)f(x)g(y)dxdy <k_{p}\left(
\int_{0}^{\infty }f^{p}(x)dx\right) ^{\frac{1}{p}}\left( \int_{0}^{\infty
}g^{q}(y)dy\right) ^{\frac{1}{q}}.  \label{3}
\end{equation}
The following inequality with the non-homogeneous kernel $h(xy)$ similar to (\ref{3})
was studied (cf. \cite{YM3}, Th. 350):

For  $h(t)>0$, satisfying $\phi(s) :=\int_{0}^{\infty }h(t)t^{s-1}dt\in \mathbf{R}_{+}, f(x),g(y)\geq 0,$ $$ 0 < \int_{0}^{\infty }x^{p-2}f^{p}(x)dx<\infty $$ and  $$0<\int_{0}^{\infty }g^{q}(y)dy < \infty,$$ we have
\begin{equation}
\int_{0}^{\infty }\int_{0}^{\infty }h(xy)f(x)g(y)dxdy < \phi (\frac{1}{p})\left( \int_{0}^{\infty }x^{p-2}f^{p}(x)dx \right)^{\frac{1}{p}}\left(\int_{0}^{\infty }g^{q}(y)dy \right) ^{\frac{1}{q}}.  \label{4}
\end{equation}

\begin{remark}
Hardy could not prove that the constant factor in (\ref{4}) is the best possible and did not consider the operator expressions of (\ref{3}) and (\ref{4}) (cf. \cite{YM3}, Chapter 9). We shall call (\ref{3}) and (\ref{4}) Hardy-Hilbert-type integral inequalities, which only contain one pair of conjugate exponents $(p,q)$.
\end{remark}

In 1998, by introducing an independent parameter $\lambda >0,$ Yang \cite%
{YM5}, \cite{YM6} gave an extension of (\ref{1}) as follows:

For $f(x),\:g(y)\geq 0,$ such that $$0<\int_{0}^{\infty }x^{1-\lambda }f^{2}(x)dx<\infty $$ and
$$0<\int_{0}^{\infty }y^{1-\lambda }g^{2}(y)dy<\infty ,$$ we have
\begin{equation}
\int_{0}^{\infty }\int_{0}^{\infty }\frac{f(x)g(y)}{(x+y)^{\lambda }}dxdy <B\left(%
\frac{\lambda }{2},\frac{\lambda }{2}\right)\left( \int_{0}^{\infty }x^{1-\lambda
}f^{2}(x)dx\int_{0}^{\infty }y^{1-\lambda }g^{2}(y)dy\right) ^{\frac{1}{2}},
\label{5}
\end{equation}%
where, the constant factor $B\left(\frac{\lambda }{2},\frac{\lambda }{2}\right)$ is the
best possible, and
\begin{equation}
B(u,v):=\int_{0}^{\infty }\frac{t^{u-1}}{(t+1)^{u+v}}dt^{{}}(u,v>0)
\label{6}
\end{equation}%
is the beta function (cf. \cite{YM6a}).

In 2004, by introducing two pairs of conjugate exponents $(p,q)$ and $(r,s)$, that is
$\frac{1}{p}+\frac{1}{q}=\frac{1}{r}+\frac{1}{s}=1$, and an independent
parameter $\lambda >0,$ Yang \cite{YM7} gave the following extension of (\ref{3}):

For $p,r>1,$ $f(x),g(y)\geq 0,$ such that $$0<\int_{0}^{\infty }x^{p(1-\frac{%
\lambda }{r})-1}f^{p}(x)dx<\infty $$ and $$0<\int_{0}^{\infty }y^{q(1-\frac{%
\lambda }{s})-1}g^{q}(y)dy<\infty ,$$ it holds
\begin{eqnarray}
&&\int_{0}^{\infty }\int_{0}^{\infty }\frac{f(x)g(y)}{x^{\lambda
}+y^{\lambda }}dxdy  \nonumber \\
&<&\frac{\pi }{\lambda \sin (\pi /r)}\left[ \int_{0}^{\infty }x^{p(1-\frac{%
\lambda }{r})-1}f^{p}(x)dx\right] ^{\frac{1}{p}}\left[ \int_{0}^{\infty
}y^{q(1-\frac{\lambda }{s})-1}g^{q}(y)dy\right] ^{\frac{1}{q}},  \label{7}
\end{eqnarray}%
where, the constant factor $$\frac{\pi }{\lambda \sin (\pi /r)}$$ is the best
possible.

For $\lambda =1,r=q,s=p,$ inequality (\ref{7}) reduces to (\ref{2}); for $\lambda
=1,r=p,s=q,$ inequality (\ref{7}) reduces to the dual form of (\ref{2}) as follows:%
\begin{eqnarray}
&&\int_{0}^{\infty }\int_{0}^{\infty }\frac{f(x)g(y)}{x+y}dxdy  \nonumber \\
&<&\frac{\pi }{\sin (\pi /p)}\left( \int_{0}^{\infty
}x^{p-2}f^{p}(x)dx\right) ^{\frac{1}{p}}\left( \int_{0}^{\infty
}y^{q-2}g^{q}(y)dy\right) ^{\frac{1}{q}}.  \label{8}
\end{eqnarray}

In 2009, replacing $1/(x^{\lambda }+y^{\lambda })$ in (\ref{7}), by a general
homogeneous kernel of degree $-\lambda $, as $k_{\lambda }(x,y)$,
Yang \cite{YM8}, \cite{YM9} proved an extension of (\ref{7}) in the following form:

For $%
\lambda >0,$ $f(x),g(y)\geq 0,$ such that $$0<\int_{0}^{\infty }x^{p(1-\frac{\lambda }{r%
})-1}f^{p}(x)dx<\infty $$ and $$0<\int_{0}^{\infty }y^{q(1-\frac{\lambda }{s}%
)-1}g^{q}(y)dy<\infty ,$$ it follows
\begin{eqnarray}
&&\int_{0}^{\infty }\int_{0}^{\infty }k_{\lambda }(x,y)f(x)g(y)dxdy
\nonumber \\
&<&k_{\lambda }(r)\left[ \int_{0}^{\infty }x^{p(1-\frac{\lambda }{r}%
)-1}f^{p}(x)dx\right] ^{\frac{1}{p}}\left[ \int_{0}^{\infty }y^{q(1-\frac{%
\lambda }{s})-1}g^{q}(y)dy\right] ^{\frac{1}{q}},  \label{9}
\end{eqnarray}%
where, the constant factor $$k_{\lambda }(r):=\int_{0}^{\infty }k_{\lambda
}(t,1)t^{\frac{\lambda }{r}-1}dt\ \in \mathbf{R}_{+}$$ is the best possible.

In \cite{YM10}, Yang presented also the following new inequality with a non-homogeneous kernel
similar to (\ref{4}):

For $f(x),g(y)\geq 0,$ satisfying $$0<\int_{0}^{\infty
}x^{p(1-\sigma)-1}f^{p}(x)dx<\infty $$ and $$0<\int_{0}^{\infty
}y^{q(1-\sigma)-1}g^{q}(y)dy<\infty, $$ we have
\begin{eqnarray}
&&\int_{0}^{\infty }\int_{0}^{\infty }h(xy)f(x)g(y)dxdy  \nonumber \\
&<&\phi (\sigma )\left( \int_{0}^{\infty }x^{p(1-\sigma
)-1}f^{p}(x)dx\right) ^{\frac{1}{p}}\left( \int_{0}^{\infty }y^{q(1-\sigma
)-1}g^{q}(y)dy\right) ^{\frac{1}{q}},  \label{10}
\end{eqnarray}%
where, the constant factor $$\phi (\sigma )=\int_{0}^{\infty }h(t)t^{\sigma
-1}dt\in \mathbf{R}_{+}$$ is the best possible.

\begin{remark}
 For $\lambda =1,r=q,s=p,$ it follows that inequality (\ref{9}) reduces to (%
\ref{3}). Hence, (\ref{9}) is an extension of (\ref{3}) with two
pairs of conjugate exponents and an independent parameter. In 2014, Yang
\cite{YM11} proved that inequalities (\ref{9}) and (\ref{10}) are equivalent for $$%
h(u)=k_{\lambda }(u,1),$$ and considered the operator expressions of (\ref{9}%
) and (\ref{10}). We call (\ref{9}) together with (\ref{10}) as
Yang-Hilbert-type integral inequalities in the first quadrant. Also we can
call some similar inequalities as Yang-Hilbert-type integral inequalities in
the whole plane.
\end{remark}

In 2007, Yang \cite{YM12} introduced a Hilbert-type integral inequality in
the whole plane as follows:

For $f(x),g(y)\geq 0,$ such that $$0< \int_{-\infty
}^{\infty }e^{-\lambda x}f^{2}(x)dx<\infty $$ and $$0<\int_{0}^{\infty }
\int_{-\infty }^{\infty }e^{-\lambda y}g^{2}(y)dy<\infty ,$$ it holds
\begin{eqnarray}
&&\int_{-\infty }^{\infty }\int_{-\infty }^{\infty }\frac{f(x)g(y)}{%
(1+e^{x+y})^{\lambda }}dxdy  \nonumber \\
&<&B(\frac{\lambda }{2},\frac{\lambda }{2})\left( \int_{-\infty }^{\infty
}e^{-\lambda x}f^{2}(x)dx\int_{-\infty }^{\infty }e^{-\lambda
y}g^{2}(y)dy\right) ^{\frac{1}{2}},  \label{11}
\end{eqnarray}%
where, the constant factor $B(\frac{\lambda }{2},\frac{\lambda }{2})\ (\lambda
>0)$ is the best possible.

For the case when $0<\lambda <1$, $p>1,$ $\frac{1}{p}+\frac{1}{q}%
=1,$ Yang \cite{YM13} proved in 2008, the following Hilbert-type integral
inequality in the whole plane:

For $p>1,$ $f(x),g(y)\geq 0,$ satisfying $$%
0<\int_{-\infty }^{\infty }|x|^{p(1-\frac{\lambda }{2})-1}f^{p}(x)dx<\infty $$
and $$0< \int_{-\infty }^{\infty }|y|^{q(1-\frac{\lambda }{2}%
)-1}g^{q}(y)dy<\infty ,$$ we have
\begin{eqnarray}
&&\int_{-\infty }^{\infty }\int_{-\infty }^{\infty }\frac{1}{|1+xy|^{\lambda
}}f(x)g(y)dxdy  \nonumber \\
&<&k_{\lambda }\left[ \int_{-\infty }^{\infty }|x|^{p(1-\frac{\lambda }{2}%
)-1}f^{p}(x)dx\right] ^{\frac{1}{p}}\left[ \int_{-\infty }^{\infty }|y|^{q(1-%
\frac{\lambda }{2})-1}g^{q}(y)dy\right] ^{\frac{1}{q}},  \label{12}
\end{eqnarray}%
where, the constant $$k_{\lambda }=B(\frac{\lambda }{2},\frac{\lambda }{2}%
)+2B(1-\lambda ,\frac{\lambda }{2})$$ is still the best possible.

Additionally, Yang et
al. \cite{YM14}-\cite{YM22} provided also some other Hilbert-type integral
inequalities in the whole plane. Rassias, M. Th. et al. \cite{MY1}-\cite{MY6}
presented as well some different new Hilbert-type inequalities.

In this paper, applying methods of Real Analysis and Functional Analysis, we build two
weight functions with parameters, and provide two kinds of parameterized
Yang-Hilbert-type integral inequalities with the best constant factors. Equivalent
forms, the reverses, and the operator expressions are also given. In
particular, the Hardy-type inequalities and Hardy-type operators are
studied. Furthermore, a number of examples with two kinds of particular kernels are
considered.

\section{Yang-Hilbert-Type Integral Inequalities in the First Quadrant}

\label{sec:2} \bigskip In this section, we present a weight function and study
some Yang-Hilbert-type integral inequalities in the first quadrant with
parameters and the best constant factors. Equivalent forms, the reverses,
the Hardy-type inequalities, the operator expressions and some particular
examples are also discussed.

\subsection{Definition of Weight Function and a Lemma}

\begin{definition}
 If $\sigma \in \mathbf{R},h(t)$ is a
non-negative measurable function in $\mathbf{R}_{+},$ define the following
weight function:%
\begin{equation}
\omega (\sigma ,y):=y^{\sigma }\int_{0}^{\infty }h(xy)x^{\sigma
-1}dx^{{}}(y\in \mathbf{R}_{+}).  \label{13}
\end{equation}
\end{definition}

Setting $t=xy$ in (\ref{13}), we obtain%
\begin{equation}
\omega (\sigma ,y)=k(\sigma ):=\int_{0}^{\infty }h(t)t^{\sigma -1}dt.
\label{14}
\end{equation}

\begin{lemma}
 If $p>0\ (p\neq 1),\frac{1}{p}+\frac{1}{q}=1,\sigma \in
\mathbf{R},$  both $h(t)$ and $f(t)$
are non-negative measurable functions in $\mathbf{R}_{+},$ and $k(\sigma )$ is defined by (\ref{14}),
then, (i) for $%
p>1,$ we have the following inequality:%
\begin{equation}
J :=\int_{0}^{\infty }y^{p\sigma -1}\left( \int_{0}^{\infty
}h(xy)f(x)dx\right) ^{p}dy
\leq k^{p}(\sigma )\int_{0}^{\infty }x^{p(1-\sigma )-1}f^{p}(x)dx;
\label{15}
\end{equation}%
(ii) for $0<p<1,$ we have the reverse of (\ref{15}).
\end{lemma}

\begin{proof}
 (i) By the weighted H\"{o}lder's inequality (cf. \cite{YM23})
and (\ref{13}), it follows that%
\begin{eqnarray}
&&\int_{0}^{\infty }h(xy)f(x)dx
=\int_{0}^{\infty }h(xy)\left[\frac{x^{(1-\sigma )/q}}{y^{(1-\sigma )/p}}f(x)\right]\left[%
\frac{y^{(1-\sigma )/p}}{x^{(1-\sigma )/q}}\right]dx \nonumber\\
&\leq &\left[ \int_{0}^{\infty }h(xy)\frac{x^{(1-\sigma )p/q}}{y^{1-\sigma }}%
f^{p}(x)dx\right] ^{\frac{1}{p}}
 \left[ \int_{0}^{\infty }h(xy)\frac{y^{(1-\sigma )q/p}}{x^{1-\sigma
}}dx\right] ^{\frac{1}{q}}   \nonumber\\
&=&(\omega (\sigma ,y))^{\frac{1}{q}}y^{\frac{1}{p}-\sigma }\left[
\int_{0}^{\infty }h(xy)\frac{x^{(1-\sigma )(p-1)}}{y^{1-\sigma }}f^{p}(x)dx%
\right] ^{\frac{1}{p}}.  \label{16}
\end{eqnarray}%
Then by (\ref{14}) and Fubini's theorem (cf. \cite{YM24}), we have%
\begin{eqnarray}
J &\leq &k^{p-1}(\sigma )\int_{0}^{\infty }\int_{0}^{\infty }h(xy)\frac{%
x^{(1-\sigma )(p-1)}}{y^{1-\sigma }}f^{p}(x)dxdy  \nonumber \\
&=&k^{p-1}(\sigma )\int_{0}^{\infty }\left[ \int_{0}^{\infty }h(xy)\frac{%
x^{(1-\sigma )(p-1)}}{y^{1-\sigma }}dy\right] f^{p}(x)dx   \nonumber \\
&=&k^{p-1}(\sigma )\int_{0}^{\infty }\omega (\sigma ,x)x^{p(1-\sigma
)-1}f^{p}(x)dx.  \label{17}
\end{eqnarray}%
 By (\ref{14}), we obtain (\ref{15}).

(ii) For $0<p<1,$ by the reverse of the weighted H\"{o}lder's inequality (cf.
\cite{YM23}), combined with (\ref{13}) and (\ref{14}), we obtain the reverse of (\ref%
{16}) and (\ref{17}). Then we get the reverse of (\ref{15}) by using (\ref%
{14}). This completes the proof of the lemma.
\end{proof}

\subsection{Two Equivalent Inequalities as well as the Reverses with the
Best Possible Constant Factors}

\begin{theorem}
 Suppose that $p>1,\frac{1}{p}+\frac{1}{q}%
=1,\sigma\in \textbf{R}, h(t)\geq 0,$ and $$k(\sigma )=\int_{0}^{\infty }h(t)t^{\sigma -1}dt\in \mathbf{R}%
_{+}.$$ If $f(x),g(y)\geq0, $ such that
$$0<\int_{0}^{\infty }x^{p(1-\sigma )-1}f^{p}(x)dx<\infty$$ and $$
0<\int_{0}^{\infty }y^{q(1-\sigma )-1}g^{q}(y)dy<\infty,$$
then we have the following equivalent inequalities:%
\begin{eqnarray}
I &:&=\int_{0}^{\infty }\int_{0}^{\infty }h(xy)f(x)g(y)dxdy   \nonumber \\
&<&k(\sigma )\left[ \int_{0}^{\infty }x^{p(1-\sigma )-1}f^{p}(x)dx\right] ^{%
\frac{1}{p}}\left[ \int_{0}^{\infty }y^{q(1-\sigma )-1}g^{q}(y)dy\right] ^{%
\frac{1}{q}},  \label{18}
\end{eqnarray}%
\begin{equation}
J=\int_{0}^{\infty }y^{p\sigma -1}\left( \int_{0}^{\infty
}h(xy)f(x)dx\right) ^{p}dy
<k^{p}(\sigma )\int_{0}^{\infty }x^{p(1-\sigma )-1}f^{p}(x)dx,  \label{19}
\end{equation}%
where, the constant factors $k(\sigma )$ and $k^{p}(\sigma )$ are the best
possible.
\end{theorem}

\begin{proof}
 We first proved that (\ref{16}) preserves the form of a
strict inequality for any $y\in \mathbf{R}_{+}$. Otherwise, there exists a $%
y>0,$ such that (\ref{16}) becomes an equality. Then, there exist two
constants $A$ and $B,$ such that they are not all zero, and (cf. \cite{YM23}%
)
$$
A\frac{x^{(1-\sigma )p/q}}{y^{1-\sigma }}f^{p}(x)=B\frac{y^{(1-\sigma )q/p}}{%
x^{1-\sigma }}\text{ \,\,a. e.\,\, in\,\, }\mathbf{R}_{+}.
$$
If $A=0,$ then $B=0,$ which is impossible. Suppose that $A\neq 0.$ Then it
follows that%
$$
x^{p(1-\sigma )-1}f^{p}(x)=y^{(1-\sigma )q}\frac{B}{Ax}\text{ \,\,a. e.\,\, in\,\, }%
\mathbf{R}_{+},
$$
which contradicts the fact that $$0<\int_{0}^{\infty }x^{p(1-\sigma
)-1}f^{p}(x)dx<\infty ,$$ in virtue of $$\int_{0}^{\infty }\frac{1}{x}dx=\infty .$$
Hence, both (\ref{16}) and (\ref{17}) preserve the forms of strict inequalities,
and thus we have (\ref{19}).

By H\"{o}lder's inequality (cf. \cite{YM23}), we obtain
\begin{eqnarray}
I &=&\int_{0}^{\infty }\left( y^{\sigma -\frac{1}{p}}\int_{0}^{\infty
}h(xy)f(x)dx\right) (y^{\frac{1}{p}-\sigma }g(y))dy   \nonumber \\
&\leq &J^{\frac{1}{p}}\left[ \int_{0}^{\infty }y^{q(1-\sigma )-1}g^{q}(y)dy%
\right] ^{\frac{1}{q}}.  \label{20}
\end{eqnarray}%
Then by (\ref{19}), we get (\ref{18}). On the other hand, assuming that (%
\ref{18}) is valid, we set%
$$
g(y):=y^{p\sigma -1}\left( \int_{0}^{\infty }h(xy)f(x)dx\right) ^{p-1},y\in
\mathbf{R}_{+}.
$$
 Then we obtain $$J=\int_{0}^{\infty }y^{q(1-\sigma )-1}g^{q}(y)dy.$$ By (\ref{15}%
), in view of $$0<\int_{0}^{\infty }x^{p(1-\sigma )-1}f^{p}(x)dx<\infty,$$
it follows that $J<\infty.$
If $J=0,$ then, (\ref{19}) is trivially valid;
if $J>0,$ then by (\ref{18}), we have%
\begin{eqnarray*}
0 &<&\int_{0}^{\infty }y^{q(1-\sigma )-1}g^{q}(y)dy=J=I \\
&<&k(\sigma )\left[ \int_{0}^{\infty }x^{p(1-\sigma )-1}f^{p}(x)dx\right] ^{%
\frac{1}{p}}\left[ \int_{0}^{\infty }y^{q(1-\sigma )-1}g^{q}(y)dy\right] ^{%
\frac{1}{q}},
\end{eqnarray*}
\begin{eqnarray*}
J^{\frac{1}{p}} =\left[ \int_{0}^{\infty }y^{q(1-\sigma )-1}g^{q}(y)dy%
\right] ^{\frac{1}{p}}
&<&k(\sigma )\left[ \int_{0}^{\infty }x^{p(1-\sigma )-1}f^{p}(x)dx\right] ^{%
\frac{1}{p}},
\end{eqnarray*}
and then (\ref{19}) follows, which is equivalent to (\ref{18}).

For any $n\in \mathbf{N}$ (where $\mathbf{N}$ is the set of positive integers), we
define the functions $f_{n}(x)$ and $g_{n}(y)$ as follows:%
$$
f_{n}(x):=\left\{
\begin{array}{c}
0,x\in (0,1) \\
x^{\sigma -\frac{1}{np}-1},x\in \lbrack 1,\infty )%
\end{array}%
,\right. g_{n}(y):=\left\{
\begin{array}{c}
y^{\sigma +\frac{1}{nq}-1},y\in (0,1] \\
0,y\in (1,\infty )%
\end{array}%
.\right.
$$
Then we find%
\begin{eqnarray*}
L_{n} &:&=\left[ \int_{0}^{\infty }x^{p(1-\sigma )-1}f_{n}^{p}(x)dx\right] ^{%
\frac{1}{p}}\left[ \int_{0}^{\infty }y^{q(1-\sigma )-1}g_{n}^{q}(y)dy\right]
^{\frac{1}{q}} \\
&=&\left( \int_{1}^{\infty }x^{-\frac{1}{n}-1}dx\right) ^{\frac{1}{p}}\left(
\int_{0}^{1}y^{\frac{1}{n}-1}dy\right) ^{\frac{1}{q}}=n.
\end{eqnarray*}
By Fubini's theorem, it follows that%
\begin{eqnarray*}
I_{n} &:&=\int_{0}^{\infty }\int_{0}^{\infty }h(xy)f_{n}(x)g_{n}(y)dxdy \\
&=&\int_{1}^{\infty }x^{\sigma -\frac{1}{np}-1}\left(
\int_{0}^{1}h(xy)y^{\sigma +\frac{1}{nq}-1}dy\right) dx \\
&=&\int_{1}^{\infty }x^{-\frac{1}{n}-1}\left(
\int_{0}^{x}h(t)t^{\sigma +\frac{1}{nq}-1}dt\right) dx\\
&=&\int_{1}^{\infty }x^{-\frac{1}{n}-1}\left( \int_{0}^{1}h(t)t^{\sigma +%
\frac{1}{nq}-1}dt+\int_{1}^{x}h(t)t^{\sigma +\frac{1}{nq}-1}dt\right) dx \\
&=&n\int_{0}^{1}h(t)t^{\sigma +\frac{1}{nq}-1}dt+\int_{1}^{\infty }x^{-\frac{%
1}{n}-1}\left( \int_{1}^{x}h(t)t^{\sigma +\frac{1}{nq}-1}dt\right) dx\\
&=&n\int_{0}^{1}h(t)t^{\sigma +\frac{1}{nq}-1}dt+\int_{1}^{\infty }\left(
\int_{t}^{\infty }x^{-\frac{1}{n}-1}dx\right) h(t)t^{\sigma +\frac{1}{nq}%
-1}dt \\
&=&n\left( \int_{0}^{1}h(t)t^{\sigma +\frac{1}{nq}-1}dt+\int_{1}^{\infty
}h(t)t^{\sigma -\frac{1}{np}-1}dt\right) .
\end{eqnarray*}

If there exists a positive number $k\leq k(\sigma ),$ such that (\ref{18})
is still valid when replacing $k(\sigma )$ by $k,$ then in particular, it
follows that $$\frac{1}{n}I_{n}<k\frac{1}{n}L_{n},$$ and%
$$
\int_{0}^{1}h(t)t^{\sigma +\frac{1}{nq}-1}dt+\int_{1}^{\infty }h(t)t^{\sigma
-\frac{1}{np}-1}dt<k.
$$
Since both $$\{h(t)t^{\sigma +\frac{1}{nq}-1}\}_{n=1}^{\infty }\ (t\in (0,1])$$
and $$\{h(t)t^{\sigma -\frac{1}{np}-1}\}_{n=1}^{\infty }\ (t\in (1,\infty ))$$
are non-negative and increasing, then by Levi's theorem (cf. \cite{YM24}), we get
\begin{eqnarray*}
k(\sigma ) &=&\int_{0}^{1}h(t)t^{\sigma -1}dt+\int_{1}^{\infty
}h(t)t^{\sigma -1}dt \\
&=&\lim_{n\rightarrow \infty }\left( \int_{0}^{1}h(t)t^{\sigma +\frac{1}{nq}%
-1}dt+\int_{1}^{\infty }h(t)t^{\sigma -\frac{1}{np}-1}dt\right) \leq k.
\end{eqnarray*}%
Thus $k=k(\sigma )$ is the best possible constant factor of (\ref{18}).

The constant factor in (\ref{19}) is still the best possible. Otherwise, by (\ref{20}) we
would reach the contradiction that the constant factor in (\ref%
{18}) is not the best possible.\\
This completes the proof of the theorem.
\end{proof}

\begin{theorem}
 Replacing $p>1$ by $0<p<1$ in Theorem 1, we obtain the
equivalent reverses of (\ref{18}) and (\ref{19}). If there exists a constant
$\delta _{\:0}>0,$ such that for any $\widetilde{\sigma }\in (\sigma -\delta
_{\:0},\sigma ],$ $$k(\widetilde{\sigma })=\int_{0}^{\infty }h(t)t^{\widetilde{%
\sigma }-1}dt\in \mathbf{R}_{+},$$ then the constant factors in the reverses
of (\ref{18}) and (\ref{19}) are the best possible.
\end{theorem}

\begin{proof}
By Lemma 1 and the reverse of H\"{o}lder's inequality, we get
the reverses of (\ref{18}), (\ref{19}) and (\ref{20}). Similarly, we
can set $g(y)$ as in Theorem 1, and prove that the reverses of (\ref{18}) and (%
\ref{19}) are equivalent.

For $n>\frac{2}{\delta _{\:0}|q|}(n\in \mathbf{N}),$ we set $f_{n}(x)$ and $%
g_{n}(y)$ as in Theorem 1. If there exists a positive number $k\geq k(\sigma ),$
such that the reverse of (\ref{18}) is valid when replacing $k(\sigma )$ by $%
k,$ then it follows that $$\frac{1}{n}I_{n}>k\frac{1}{n}L_{n},$$
and%
\begin{equation}
\int_{0}^{1}h(t)t^{\sigma +\frac{1}{nq}-1}dt+\int_{1}^{\infty }h(t)t^{\sigma
-\frac{1}{np}-1}dt>k.  \label{21}
\end{equation}

Since $\{h(t)t^{\sigma -\frac{1}{np}-1}\}_{n=1}^{\infty }(t\in (1,\infty ))$
is still non-negative and increasing, by Levi's theorem it follows that%
$$
\lim_{n\rightarrow \infty }\int_{1}^{\infty }h(t)t^{\sigma -\frac{1}{np}%
-1}dt=\int_{1}^{\infty }h(t)t^{\sigma -1}dt.
$$
Since $$0\leq h(t)t^{\sigma +\frac{1}{nq}-1}\leq h(t)t^{(\sigma -\frac{\delta
_{\:0}}{2})-1}\ \left(t\in (0,1],n>\frac{2}{\delta _{\:0}|q|}\right),$$ and
$$
0\leq \int_{0}^{1}h(t)t^{(\sigma -\frac{\delta _{\:0}}{2})-1}dt\leq k(\sigma -%
\frac{\delta _{\:0}}{2})<\infty ,
$$
then by Lebesgue's dominated convergence theorem (cf. \cite{YM24}), it follows that%
$$
\lim_{n\rightarrow \infty }\int_{0}^{1}h(t)t^{\sigma +\frac{1}{nq}%
-1}dt=\int_{0}^{1}h(t)t^{\sigma -1}dt.
$$

In view of the above results and (\ref{21}), we have
$$
k(\sigma )=\lim_{n\rightarrow \infty }\left( \int_{0}^{1}h(t)t^{\sigma +%
\frac{1}{nq}-1}dt+\int_{1}^{\infty }h(t)t^{\sigma -\frac{1}{np}-1}dt\right)
\geq k.
$$
Then $k=k(\sigma )$ is the best possible constant factor for the reverse
of (\ref{18}).

Following the same method, we can prove that the constant factor in the reverse of (%
\ref{19}) is the best possible, by the use of the reverse of (\ref{20}).\\
This completes the proof of the theorem.
\end{proof}

\subsection{\protect Yang-Hilbert-Type Integral Inequalities in the
First Quadrant with Multi-Variables}
\begin{theorem}
 Suppose that $p>1,\:\frac{1}{p}+\frac{1}{q}%
=1,\:h(t)\geq 0,\sigma\in \textbf{R},$ $$k(\sigma )=\int_{0}^{\infty }h(t)t^{\sigma -1}dt\in \mathbf{R}%
_{+},$$ $\delta _{i}\in \{-1,1\},$ $0\leq a_{i}<b_{i}\leq \infty ,\:v_{i}^{\prime
}(s)>0\ (s\in (a_{i},b_{i}))$ and $$v_{i}(a_{i}^{+})=\lim_{s\rightarrow a_{i}+}v_{i}(s)=0,$$ $$v_{i}(b_{i}^{-})=\lim_{s\rightarrow b_{i}-}v_{i}(s)=\infty\
(i=1,2).$$
If $f(x),g(y)\geq 0,$ such that%
$$
0 <\int_{a_{1}}^{b_{1}}\frac{(v_{1}(x))^{p(1-\delta _{1}\sigma )-1}}{%
(v_{1}^{\prime }(x))^{p-1}}f^{p}(x)dx<\infty $$
and $$
0 <\int_{a_{2}}^{b_{2}}\frac{(v_{2}(y))^{q(1-\delta _{2}\sigma )-1}}{%
(v_{2}^{\prime }(y))^{q-1}}g^{q}(y)dy<\infty,
$$
then we have the following equivalent inequalities:%
\begin{eqnarray}
&&\int_{a_{2}}^{b_{2}}\int_{a_{1}}^{b_{1}}h(v_{1}^{\delta
_{1}}(x)v_{2}^{\delta _{2}}(y))f(x)g(y)dxdy   \nonumber \\
&<&k(\sigma )\left[ \int_{a_{1}}^{b_{1}}\frac{(v_{1}(x))^{p(1-\delta
_{1}\sigma )-1}}{(v_{1}^{\prime }(x))^{p-1}}f^{p}(x)dx\right] ^{\frac{1}{p}}
 \nonumber \\
&&\times \left[ \int_{a_{2}}^{b_{2}}\frac{(v_{2}(y))^{q(1-\delta _{2}\sigma
)-1}}{(v_{2}^{\prime }(y))^{q-1}}g^{q}(y)dy\right] ^{\frac{1}{q}},
\label{22}
\end{eqnarray}%
\begin{eqnarray}
&&\int_{a_{2}}^{b_{2}}\frac{v_{2}^{\prime }(y)}{(v_{2}(y))^{1-p\delta
_{2}\sigma }}\left( \int_{a_{1}}^{b_{1}}h(v_{1}^{\delta
_{1}}(x)v_{2}^{\delta _{2}}(y))f(x)dx\right) ^{p}dy   \nonumber \\
&<&k^{p}(\sigma )\int_{a_{1}}^{b_{1}}\frac{(v_{1}(x))^{p(1-\delta _{1}\sigma
)-1}}{(v_{1}^{\prime }(x))^{p-1}}f^{p}(x)dx,  \label{23}
\end{eqnarray}%
where the constant factors $k(\sigma )$ and $k^{p}(\sigma )$ are the best
possible.
\end{theorem}

\begin{proof}
Setting $$x=v_{1}^{\delta _{1}}(s),\:y=v_{2}^{\delta _{2}}(t)$$ in
(\ref{18}), since $\delta _{i}\in \{-1,1\},$ we get $$dx=\delta
_{1}v_{1}^{\delta _{1}-1}(s)v_{1}^{\prime }(s)ds,\ \ dy=\delta
_{2}v_{2}^{\delta _{2}-1}(t)v_{2}^{\prime }(t)dt,$$ and%
$$
I=|\delta_{1}\delta_{2}|\int_{a_{2}}^{b_{2}}\int_{a_{1}}^{b_{1}}h(v_{1}^{\delta
_{1}}(s)v_{2}^{\delta _{2}}(t))f(v_{1}^{\delta _{1}}(s))g(v_{2}^{\delta
_{2}}(t))v_{1}^{\delta _{1}-1}(s)v_{1}^{\prime }(s)v_{2}^{\delta
_{2}-1}(t)v_{2}^{\prime }(t)dsdt
$$
$$
=\int_{a_{2}}^{b_{2}}\int_{a_{1}}^{b_{1}}h(v_{1}^{\delta
_{1}}(s)v_{2}^{\delta _{2}}(t))(f(v_{1}^{\delta _{1}}(s))v_{1}^{\delta
_{1}-1}(s)v_{1}^{\prime }(s))(g(v_{2}^{\delta _{2}}(t))v_{2}^{\delta
_{2}-1}(t)v_{2}^{\prime }(t))dsdt,
$$
\begin{eqnarray*}
I_{1} &:&=\int_{0}^{\infty }x^{p(1-\sigma
)-1}f^{p}(x)dx=\int_{a_{1}}^{b_{1}}(v_{1}^{\delta _{1}}(s))^{p(1-\sigma
)-1}f^{p}(v_{1}^{\delta _{1}}(s))v_{1}^{\delta _{1}-1}(s)v_{1}^{\prime
}(s)ds, \\
I_{2} &:&=\int_{0}^{\infty }y^{q(1-\sigma
)-1}g^{q}(y)dy=\int_{a_{2}}^{b_{2}}(v_{2}^{\delta _{2}}(t))^{q(1-\sigma
)-1}g^{q}(v_{2}^{\delta _{2}}(t))v_{2}^{\delta _{2}-1}(t)v_{2}^{\prime
}(t)dt.
\end{eqnarray*}%
Setting $$F(s)=f(v_{1}^{\delta _{1}}(s))v_{1}^{\delta _{1}-1}(s)v_{1}^{\prime
}(s),\ \
G(t)=g(v_{2}^{\delta _{2}}(t))v_{2}^{\delta _{2}-1}(t)v_{2}^{\prime }(t),
$$
we obtain%
\begin{eqnarray*}
f^{p}(v_{1}^{\delta _{1}}(s)) &=&v_{1}^{p(1-\delta _{1})}(s)(v_{1}^{\prime
}(s))^{-p}F^{p}(s), \\
g^{q}(v_{2}^{\delta _{2}}(t)) &=&v_{2}^{q(1-\delta _{2})}(t)(v_{2}^{\prime
}(t))^{-q}G^{q}(t),
\end{eqnarray*}
and then it follows that
\begin{eqnarray*}
I &=&\int_{a_{2}}^{b_{2}}\int_{a_{1}}^{b_{1}}h(v_{1}^{\delta
_{1}}(s)v_{2}^{\delta _{2}}(t))F(s)G(t)dsdt, \\
I_{1} &=&\int_{a_{1}}^{b_{1}}\frac{(v_{1}(s))^{p(1-\delta _{1}\sigma )-1}}{%
(v_{1}^{\prime }(s))^{p-1}}F^{p}(s)ds,\ \
I_{2}=\int_{a_{2}}^{b_{2}}\frac{(v_{2}(t))^{q(1-\delta _{2}\sigma )-1}}{%
(v_{2}^{\prime }(t))^{q-1}}G^{q}(t)dt.
\end{eqnarray*}%
\

Substituting the above results in (\ref{18}), resetting $$%
s=x,\:t=y,\:F(s)=f(x),\:G(t)=g(y),$$ we obtain (\ref{22}). Similarly, we have
(\ref{23}).

On the other hand, if we set $$v_{1}^{\delta
_{1}}(x)=x,\:v_{2}^{\delta _{2}}(y)=y,\:a_{i}=0,\:b_{i}=\infty\ (i=1,2)$$ in (%
\ref{22}), we obtain (\ref{18}). Hence, the inequalities (\ref{22}) and (\ref{18})
are equivalent. It is evident that the inequalities (\ref{23}) and (\ref{19})
are equivalent. Hence, the inequalities (\ref{22}) and (\ref{23}) are
equivalent. Since the constant factors in (\ref{18}) and (\ref{19}) are the
best possible, it follows that the constant factors in (\ref{22}) and (\ref{23}) are
also the best possible by using the equivalency.

This completes the proof of the theorem.
\end{proof}
\begin{theorem}
 Replacing $p>1$ by $0<p<1$ in Theorem 3, we
obtain the equivalent reverses of (\ref{22}) and (\ref{23}). If there exists a
constant $\delta _{\:0}>0,$ such that for any value of $\widetilde{\sigma }\in (\sigma
-\delta _{\:0},\sigma ],$ $$k(\widetilde{\sigma })=\int_{0}^{\infty }h(t)t^{%
\widetilde{\sigma }-1}dt\in \mathbf{R}_{+},$$ then the constant factors in
the reverses of (\ref{22}) and (\ref{23}) are the best possible.
\end{theorem}

In particular, for $\delta _{1}=\delta _{2}=1$ in Theorem 3 and Theorem 4,
we get the following integral inequalities with the non-homogeneous kernel:

\begin{corollary}
 Suppose that $p>1,\:\frac{1}{p}+\frac{1}{q}=1,\: h(t)\geq
0, \sigma\in \textbf{R},$
\begin{eqnarray*}
k(\sigma )=\int_{0}^{\infty }h(t)t^{\sigma -1}dt\in \mathbf{R}_{+},
\end{eqnarray*}
$0\leq a_{i}<b_{i}\leq \infty ,v_{i}^{\prime }(s)>0(s\in
(a_{i},b_{i})),v_{i}(a_{i}^{+})=0,v_{i}(b_{i}^{-})=\infty\ (i=1,2).$
If $%
f(x),g(y)\geq 0,$ such that%
$$
0 <\int_{a_{1}}^{b_{1}}\frac{(v_{1}(x))^{p(1-\sigma )-1}}{(v_{1}^{\prime
}(x))^{p-1}}f^{p}(x)dx<\infty $$
and
$$0<\int_{a_{2}}^{b_{2}}\frac{(v_{2}(y))^{q(1-\sigma )-1}}{(v_{2}^{\prime
}(y))^{q-1}}g^{q}(y)dy<\infty,$$
then we have the following equivalent inequalities:%
\begin{eqnarray}
&&\int_{a_{2}}^{b_{2}}\int_{a_{1}}^{b_{1}}h(v_{1}(x)v_{2}(y))f(x)g(y)dxdy
 \nonumber \\
&<&k(\sigma )\left[ \int_{a_{1}}^{b_{1}}\frac{(v_{1}(x))^{p(1-\sigma )-1}}{%
(v_{1}^{\prime }(x))^{p-1}}f^{p}(x)dx\right] ^{\frac{1}{p}}   \nonumber \\
&&\times \left[ \int_{a_{2}}^{b_{2}}\frac{(v_{2}(y))^{q(1-\sigma )-1}}{%
(v_{2}^{\prime }(y))^{q-1}}g^{q}(y)dy\right] ^{\frac{1}{q}},  \label{24}
\end{eqnarray}%
\begin{eqnarray}
&&\int_{a_{2}}^{b_{2}}\frac{v_{2}^{\prime }(y)}{(v_{2}(y))^{1-p\sigma }}%
\left( \int_{a_{1}}^{b_{1}}h(v_{1}(x)v_{2}(y))f(x)dx\right) ^{p}dy   \nonumber \\
&<&k^{p}(\sigma )\int_{a_{1}}^{b_{1}}\frac{(v_{1}(x))^{p(1-\sigma )-1}}{%
(v_{1}^{\prime }(x))^{p-1}}f^{p}(x)dx,  \label{25}
\end{eqnarray}%
where the constant factors $k(\sigma )$ and $k^{p}(\sigma )$ are the best
possible.
\end{corollary}

Replacing $p>1$ by $0<p<1$ in the above inequalities, we obtain the equivalent
reverses of (\ref{24}) and (\ref{25}).

 If there exists a constant $\delta
_{\:0}>0,$ such that for any $\widetilde{\sigma }\in (\sigma -\delta
_{\:0},\sigma ],$ $$k(\widetilde{\sigma })=\int_{0}^{\infty }h(t)t^{\widetilde{%
\sigma }-1}dt\in \mathbf{R}_{+},$$ then the constant factors in the reverses
of (\ref{24}) and (\ref{25}) are the best possible.

In particular, for $\delta _{1}=-1,\delta _{2}=1$ in Theorem 3 and Theorem
4, setting $$h(t)=k_{\lambda }(1,t)$$ (cf. Definition 1), we find%
\[
h(\frac{v_{2}(y)}{v_{1}(x)})=k_{\lambda }(1,\frac{v_{2}(y)}{v_{1}(x)}%
)=v_{1}^{\lambda }(x)k_{\lambda }(v_{1}(x),v_{2}(y)).
\]
Replacing $f(x)$ by $v_{1}^{-\lambda }(x)f(x),$ it follows that $%
[v_{1}(x)]^{p(1+\sigma )-1}f^{p}(x)$ is replaced by
\[
\lbrack v_{1}(x)]^{p(1+\sigma )-1}[v_{1}^{-\lambda
}(x)f(x)]^{p}=[v_{1}(x)]^{p(1-\mu )-1}f^{p}(x),
\]%
and we have the following integral inequalities with the homogeneous kernel:
\begin{corollary}
 Suppose that $p>1,\frac{1}{p}+\frac{1}{q}=1,\mu
,\sigma \in \mathbf{R},\mu +\sigma =\lambda ,k_{\lambda }(x,y)$ is a
homogeneous function of degree $-\lambda $ in $\mathbf{R}_{+}^{2},$%
\begin{eqnarray*}
k_{\lambda }(\sigma )=\int_{0}^{\infty }k_{\lambda }(1,t)t^{\sigma -1}dt\in
\mathbf{R}_{+},
\end{eqnarray*}%
$0\leq a_{i}<b_{i}\leq \infty ,v_{i}^{\prime }(s)>0$ $(s\in
(a_{i},b_{i})),v_{i}(a_{i}^{+})=0,v_{i}(b_{i}^{-})=\infty\ (i=1,2).$ If $%
f(x),g(y)\geq 0,$ such that%
$$
0 <\int_{a_{1}}^{b_{1}}\frac{(v_{1}(x))^{p(1-\mu )-1}}{(v_{1}^{\prime
}(x))^{p-1}}f^{p}(x)dx<\infty$$
and
$$0 <\int_{a_{2}}^{b_{2}}\frac{(v_{2}(y))^{q(1-\sigma )-1}}{(v_{2}^{\prime
}(y))^{q-1}}g^{q}(y)dy<\infty,$$
then we have the following equivalent inequalities:%
\begin{eqnarray}
&&\int_{a_{2}}^{b_{2}}\int_{a_{1}}^{b_{1}}k_{\lambda
}(v_{1}(x),v_{2}(y))f(x)g(y)dxdy   \nonumber \\
&<&k_{\lambda }(\sigma )\left[ \int_{a_{1}}^{b_{1}}\frac{(v_{1}(x))^{p(1-\mu
)-1}}{(v_{1}^{\prime }(x))^{p-1}}f^{p}(x)dx\right] ^{\frac{1}{p}}   \nonumber \\
&&\times \left[ \int_{a_{2}}^{b_{2}}\frac{(v_{2}(y))^{q(1-\sigma )-1}}{%
(v_{2}^{\prime }(y))^{q-1}}g^{q}(y)dy\right] ^{\frac{1}{q}},  \label{26}
\end{eqnarray}%
\begin{eqnarray}
&&\int_{a_{2}}^{b_{2}}\frac{v_{2}^{\prime }(y)}{(v_{2}(y))^{1-p\sigma }}%
\left( \int_{a_{1}}^{b_{1}}k_{\lambda }(v_{1}(x),v_{2}(y))f(x)dx\right)
^{p}dy   \nonumber \\
&<&k_{\lambda }^{p}(\sigma )\int_{a_{1}}^{b_{1}}\frac{(v_{1}(x))^{p(1-\mu
)-1}}{(v_{1}^{\prime }(x))^{p-1}}f^{p}(x)dx,  \label{27}
\end{eqnarray}%
where the constant factors $k_{\lambda }(\sigma )$ and $k_{\lambda
}^{p}(\sigma )$ are the best possible.
\end{corollary}

Replacing $p>1$ by $0<p<1$ in the above inequalities, we obtain the equivalent
reverses of (\ref{26}) and (\ref{27}). If there exists a constant $\delta
_{\:0}>0,$ such that for any $\widetilde{\sigma }\in (\sigma -\delta
_{\:0},\sigma ],$ $$k_{\lambda }(\widetilde{\sigma })=\int_{0}^{\infty
}k_{\lambda }(1,t)t^{\widetilde{\sigma }-1}dt\in \mathbf{R}_{+},$$ then the
constant factors in the reverses of (\ref{26}) and (\ref{27}) are the best
possible.

Setting $$a_{i}=0,\:b_{i}=\infty\ (i=1,2),\:v_{1}(x)=x,\:v_{2}(y)=y$$ in Corollary 2,
we have

\begin{corollary}
 Suppose that $p>1,\frac{1}{p}+\frac{1}{q}%
=1,\mu ,\sigma \in \mathbf{R},\mu +\sigma =\lambda ,$ $k_{\lambda }(x,y)$ is
a homogeneous function of degree $-\lambda $ in $\mathbf{R}_{+}^{2},$ and
\begin{eqnarray*}
k_{\lambda }(\sigma )=\int_{0}^{\infty }k_{\lambda }(1,t)t^{\sigma -1}dt\in
\mathbf{R}_{+}.
\end{eqnarray*}
If $f(x),g(y)\geq 0,$
$$
0<\int_{0}^{\infty }x^{p(1-\mu )-1}f^{p}(x)dx<\infty$$
and $$  0<\int_{0}^{\infty
}g^{q(1-\sigma )-1}(y)dy<\infty,
$$
 then we have the following equivalent inequalities:%
\begin{eqnarray}
&&\int_{0}^{\infty }\int_{0}^{\infty }k_{\lambda }(x,y)f(x)g(y)dxdy   \nonumber
\\
&<&k_{\lambda }(\sigma )\left[ \int_{0}^{\infty }x^{p(1-\mu )-1}f^{p}(x)dx%
\right] ^{\frac{1}{p}}\left[ \int_{0}^{\infty }y^{q(1-\sigma )-1}g^{q}(y)dy%
\right] ^{\frac{1}{q}},  \label{28}
\end{eqnarray}%
\begin{equation}
\int_{0}^{\infty }y^{p\sigma -1}\left( \int_{0}^{\infty }k_{\lambda
}(x,y)f(x)dx\right) ^{p}dy
<k_{\lambda }^{p}(\sigma )\int_{0}^{\infty }x^{p(1-\mu )-1}f^{p}(x)dx,
\label{29}
\end{equation}%
where, the constant factors $k_{\lambda }(\sigma )$ and $k_{\lambda
}^{p}(\sigma )$ are the best possible.
\end{corollary}

Replacing $p>1$ by $0<p<1$ in the above inequalities, we obtain the equivalent
reverses of (\ref{28}) and (\ref{29}). If there exists a constant $\delta
_{\:0}>0,$ such that for any $\widetilde{\sigma }\in (\sigma -\delta
_{\:0},\sigma ],$
\[
k_{\lambda }(\widetilde{\sigma })=\int_{0}^{\infty }k_{\lambda }(1,t)t^{%
\widetilde{\sigma }-1}dt\in \mathbf{R}_{+},
\]
then the constant factors in the reverses of (\ref{28}) and (\ref{29}) are
the best possible.

\begin{remark}
 (a) It is evident that (\ref{18}) and (\ref{28}) are
equivalent for $h(t)=k_{\lambda }(1,t).$ The same holds for (\ref{19}) and (\ref{29}).
\end{remark}

(b) In the following, we list the functions $v_{i}(s)\ (i=1,2)$ which satisfy the
conditions of Theorem 3 and Theorem 4:

(i) $v_{i}(s)=s^{a},\:s\in (0,\infty )\ (a\in \mathbf{R}_{+}\mathbf{)},$
with $v_{i}^{\prime }(s)=as^{a-1}>0;$

(ii) $v_{i}(s)=\tan ^{a}s,\:s\in (0,\frac{\pi }{2})\ (a\in \mathbf{R}_{+}\mathbf{%
)},$ with $v_{i}^{\prime }(s)=a\tan ^{a-1}s\:\sec ^{2}s>0;$

(iii) $v_{i}(s)=\ln ^{a}s,\:s\in (1,\infty )\ (a\in \mathbf{R}_{+}\mathbf{)},$
with $v_{i}^{\prime }(s)=\frac{a}{s}\ln ^{a-1}s>0;$

(iv) $v_{i}(s)=e^{as}-1,\:s\in (0,\infty )\ (a\in \mathbf{R}_{+}\mathbf{)},$
with $v_{i}^{\prime }(s)=ae^{as}>0.$

\subsection{Hardy-Type Integral Inequalities with Multi-Variables}

In the following two sections, if the constant factors in the inequalities
(operator inequalities) are related to $k^{(1)}(\sigma)$ ( or $%
k_{\lambda}^{(1)}(\sigma)$ ), then we call them Hardy-type
inequalities (operators) of the \textit{first} kind; if the constant factors in the inequalities
(operator inequalities) are related to $k^{(2)}(\sigma)$ ( or $%
k_{\lambda}^{(2)}(\sigma)$ ), then we call them Hardy-type
inequalities (operators) of the \textit{second} kind.

\bigskip If $h(t)=0\ (t>1)$, then we have $h(xy)=0\ (x>\frac{1}{y}>0),$ and%
\[
k(\sigma )=\int_{0}^{\infty }h(t)t^{\sigma -1}dt=\int_{0}^{1}h(t)t^{\sigma
-1}dt.
\]
Setting
\begin{equation}
k^{(1)}(\sigma ):=\int_{0}^{1}h(t)t^{\sigma -1}dt,  \label{30}
\end{equation}%
by Theorem 1 and Theorem 2, we obtain the following first kind Hardy-type
integral inequalities with non-homogeneous kernel:

\begin{corollary}
 Suppose that $p>1,\frac{1}{p}+\frac{1}{q}%
=1,\:h(t)\geq 0,\sigma\in \textbf{R},$%
\begin{eqnarray*}
k^{(1)}(\sigma )=\int_{0}^{1}h(t)t^{\sigma -1}dt\in \mathbf{R}_{+}.
\end{eqnarray*}
If $f(x),g(y)\geq 0,$%
$$
0<\int_{0}^{\infty }x^{p(1-\sigma )-1}f^{p}(x)dx<\infty$$ and $$0<\int_{0}^{\infty
}y^{q(1-\sigma )-1}g^{q}(y)dy<\infty,
$$
then we have the following equivalent inequalities:%
\begin{eqnarray}
&&\int_{0}^{\infty }\left( \int_{0}^{\frac{1}{y}}h(xy)f(x)dx\right) g(y)dy
=\int_{0}^{\infty }\left( \int_{0}^{\frac{1}{x}}h(xy)g(y)dy\right) f(x)dx\nonumber\\
&<&k^{(1)}(\sigma )\left[ \int_{0}^{\infty }x^{p(1-\sigma )-1}f^{p}(x)dx\right]
^{\frac{1}{p}}\left[ \int_{0}^{\infty }y^{q(1-\sigma )-1}g^{q}(y)dy\right] ^{%
\frac{1}{q}},  \label{31}
\end{eqnarray}%
\begin{equation}
\int_{0}^{\infty }y^{p\sigma -1}\left( \int_{0}^{\frac{1}{y}%
}h(xy)f(x)dx\right) ^{p}dy
<(k^{(1)}(\sigma ))^{p}\int_{0}^{\infty }x^{p(1-\sigma )-1}f^{p}(x)dx;
\label{32}
\end{equation}%
where, the constant factors $k^{(1)}(\sigma )$ and $(k^{(1)}(\sigma ))^{p}$
are the best possible.
\end{corollary}

Replacing $p>1$ by $0<p<1$ in the above inequalities, we obtain the equivalent
reverses of (\ref{31}) and (\ref{32}). If there exists a constant $\delta
_{\:0}>0,$ such that for any $\widetilde{\sigma }\in (\sigma -\delta
_{\:0},\sigma ],$%
\[
k^{(1)}(\widetilde{\sigma })=\int_{0}^{1}h(t)t^{\widetilde{\sigma }-1}dt\in
\mathbf{R}_{+},
\]
then the constant factors in the reverses of (\ref{31}) and (\ref{32}) are
the best possible.

If $h(t)=0\ (t>1)$ in Corollary 1, then $$h(v_{1}(x)v_{2}(y))=0\ \ (v_{1}(x)>\frac{1%
}{v_{2}(y)}>0),$$ and therefore we obtain the following general results:

\begin{corollary}
 Suppose that $p>1,\:\frac{1}{p}+\frac{1}{q}=1,\:h(t)\geq
0, \sigma\in \textbf{R},$%
\begin{eqnarray*}
k^{(1)}(\sigma )=\int_{0}^{1}h(t)t^{\sigma -1}dt\in \mathbf{R}_{+},
\end{eqnarray*}
$0\leq a_{i}<b_{i}\leq \infty ,\:v_{i}^{\prime }(s)>0\ (s\in
(a_{i},b_{i})),v_{i}(a_{i}^{+})=0,v_{i}(b_{i}^{-})=\infty\ (i=1,2).$
If $%
f(x),\:g(y)\geq 0,$ such that%
$$
0 <\int_{a_{1}}^{b_{1}}\frac{(v_{1}(x))^{p(1-\sigma )-1}}{(v_{1}^{\prime
}(x))^{p-1}}f^{p}(x)dx<\infty$$ and $$
0 <\int_{a_{2}}^{b_{2}}\frac{(v_{2}(y))^{q(1-\sigma )-1}}{(v_{2}^{\prime
}(y))^{q-1}}g^{q}(y)dy<\infty ,
$$
then we have the following equivalent inequalities:%
\begin{eqnarray}
&&\int_{a_{2}}^{b_{2}}\left( \int_{a_{1}}^{v_{1}^{-1}(\frac{1}{v_{2}(y)}%
)}h(v_{1}(x)v_{2}(y))f(x)dx\right) g(y)dy   \nonumber \\
&=&\int_{a1}^{b_{1}}\left( \int_{a_{2}}^{v_{2}^{-1}(\frac{1}{v_{1}(x)}%
)}h(v_{1}(x)v_{2}(y))g(y)dy\right) f(x)dx   \nonumber \\
&<&k^{(1)}(\sigma )\left[ \int_{a_{1}}^{b_{1}}\frac{(v_{1}(x))^{p(1-\sigma
)-1}}{(v_{1}^{\prime }(x))^{p-1}}f^{p}(x)dx\right] ^{\frac{1}{p}}   \nonumber \\
&&\times \left[ \int_{a_{2}}^{b_{2}}\frac{(v_{2}(y))^{q(1-\sigma )-1}}{%
(v_{2}^{\prime }(y))^{q-1}}g^{q}(y)dy\right] ^{\frac{1}{q}},  \label{33}
\end{eqnarray}%
\begin{eqnarray}
&&\int_{a_{2}}^{b_{2}}\frac{v_{2}^{\prime }(y)}{(v_{2}(y))^{1-p\sigma }}%
\left( \int_{a_{1}}^{v_{1}^{-1}(\frac{1}{v_{2}(y)}%
)}h(v_{1}(x)v_{2}(y))f(x)dx\right) ^{p}dy   \nonumber \\
&<&(k^{(1)}(\sigma ))^{p}\int_{a_{1}}^{b_{1}}\frac{(v_{1}(x))^{p(1-\sigma
)-1}}{(v_{1}^{\prime }(x))^{p-1}}f^{p}(x)dx,  \label{34}
\end{eqnarray}%
where, the constant factors $k^{(1)}(\sigma )$ and $(k^{(1)}(\sigma ))^{p}$
are the best possible.
\end{corollary}

Replacing $p>1$ by $0<p<1$ in the above inequalities, we obtain the equivalent
reverses of (\ref{33}) and (\ref{34}). If there exists a constant $\delta
_{\:0}>0,$ such that for any $\widetilde{\sigma }\in (\sigma -\delta
_{\:0},\sigma ],$ $$k^{(1)}(\widetilde{\sigma })=\int_{0}^{1}h(t)t^{\widetilde{%
\sigma }-1}dt\in \mathbf{R}_{+},$$ then the constant factors in the reverses
of (\ref{33}) and (\ref{34}) are the best possible.

If $h(t)=0\ (0<t<1)$, then $h(xy)=0\ (0<x<\frac{1}{y}),$ and%
\[
k(\sigma )=\int_{0}^{\infty }h(t)t^{\sigma -1}dt=\int_{1}^{\infty
}h(t)t^{\sigma -1}dt.
\]
Setting%
\begin{equation}
k^{(2)}(\sigma ):=\int_{1}^{\infty }h(t)t^{\sigma -1}dt,  \label{35}
\end{equation}%
by Theorem 1 and Theorem 2, we have the following second kind Hardy-type
integral inequalities with the non-homogeneous kernel:

\begin{corollary}
 Suppose that $p>1,\frac{1}{p}+\frac{1}{q}=1,h(t)\geq
0, \sigma\in \textbf{R},$%
\begin{eqnarray*}
k^{(2)}(\sigma )=\int_{1}^{\infty }h(t)t^{\sigma -1}dt\in \mathbf{R}_{+}.
\end{eqnarray*}
If $f(x),g(y)\geq 0,$%
$$
0<\int_{0}^{\infty }x^{p(1-\sigma )-1}f^{p}(x)dx<\infty $$ and $$ 0<\int_{0}^{\infty
}y^{q(1-\sigma )-1}g^{q}(y)dy<\infty,
$$
 then we have the following equivalent inequalities:%
\begin{eqnarray*}
\int_{0}^{\infty }\left( \int_{\frac{1}{y}}^{\infty }h(xy)f(x)dx\right)
g(y)dy
=\int_{0}^{\infty }\left( \int_{\frac{1}{x}}^{\infty }h(xy)g(y)dy\right)
f(x)dx
\end{eqnarray*}
\begin{equation}
<k^{(2)}(\sigma )\left[ \int_{0}^{\infty }x^{p(1-\sigma )-1}f^{p}(x)dx\right]
^{\frac{1}{p}}\left[ \int_{0}^{\infty }y^{q(1-\sigma )-1}g^{q}(y)dy\right] ^{%
\frac{1}{q}},  \label{36}
\end{equation}%
\begin{equation}
\int_{0}^{\infty }y^{p\sigma -1}\left( \int_{\frac{1}{y}}^{\infty
}h(xy)f(x)dx\right) ^{p}dy
<(k^{(2)}(\sigma ))^{p}\int_{0}^{\infty }x^{p(1-\sigma )-1}f^{p}(x)dx,
\label{37}
\end{equation}%
where, the constant factors $k^{(2)}(\sigma )$ and $(k^{(2)}(\sigma ))^{p}$
are the best possible.
\end{corollary}

Replacing $p>1$ by $0<p<1$ in the above inequalities, we derive the equivalent
reverses of (\ref{36}) and (\ref{37}).

If there exists a constant $\delta
_{\:0}>0,$ such that for any $\widetilde{\sigma }\in (\sigma -\delta
_{\:0},\sigma ],$
\[
k^{(2)}(\widetilde{\sigma })=\int_{1}^{\infty }h(t)t^{\widetilde{\sigma }%
-1}dt\in \mathbf{R}_{+},
\]
then the constant factors in the reverses of (\ref{36}) and (\ref{37}) are
the best possible.

If $h(t)=0\ (0<t<1)$ in Corollary 1, then $$h(v_{1}(x)v_{2}(y))=0\ (0<v_{1}(x)<%
\frac{1}{v_{2}(y)}).$$
We have the following general results:

\begin{corollary}
 Suppose that $p>1,\frac{1}{p}+\frac{1}{q}=1,h(t)\geq
0, \sigma\in \textbf{R},$%
\begin{eqnarray*}
k^{(2)}(\sigma )=\int_{1}^{\infty }h(t)t^{\sigma -1}dt\in \mathbf{R}_{+},
\end{eqnarray*}
$0\leq a_{i}<b_{i}\leq \infty ,\:v_{i}^{\prime }(s)>0\ (s\in
(a_{i},b_{i})),\:v_{i}(a_{i}^{+})=0,\:v_{i}(b_{i}^{-})=\infty\ (i=1,2).$ If $%
f(x),g(y)\geq 0,$ such that%
$$
0 <\int_{a_{1}}^{b_{1}}\frac{(v_{1}(x))^{p(1-\sigma )-1}}{(v_{1}^{\prime
}(x))^{p-1}}f^{p}(x)dx<\infty$$ and $$
0 <\int_{a_{2}}^{b_{2}}\frac{(v_{2}(y))^{q(1-\sigma )-1}}{(v_{2}^{\prime
}(y))^{q-1}}g^{q}(y)dy<\infty,
$$
then we have the following equivalent inequalities:%
\begin{eqnarray}
&&\int_{a_{2}}^{b_{2}}\left( \int_{v_{1}^{-1}(\frac{1}{v_{2}(y)}%
)}^{b_{1}}h(v_{1}(x)v_{2}(y))f(x)dx\right) g(y)dy \nonumber\\
&=&\int_{a1}^{b_{1}}\left( \int_{v_{2}^{-1}(\frac{1}{v_{1}(x)}%
)}^{b_{2}}h(v_{1}(x)v_{2}(y))g(y)dy\right) f(x)dx\nonumber\\
&<&k^{(2)}(\sigma )\left[ \int_{a_{1}}^{b_{1}}\frac{(v_{1}(x))^{p(1-\sigma
)-1}}{(v_{1}^{\prime }(x))^{p-1}}f^{p}(x)dx\right] ^{\frac{1}{p}}   \nonumber \\
&&\times \left[ \int_{a_{2}}^{b_{2}}\frac{(v_{2}(y))^{q(1-\sigma )-1}}{%
(v_{2}^{\prime }(y))^{q-1}}g^{q}(y)dy\right] ^{\frac{1}{q}},  \label{38}
\end{eqnarray}%
\begin{eqnarray}
&&\int_{a_{2}}^{b_{2}}\frac{v_{2}^{\prime }(y)}{(v_{2}(y))^{1-p\sigma }}%
\left( \int_{v_{1}^{-1}(\frac{1}{v_{2}(y)}%
)}^{b_{1}}h(v_{1}(x)v_{2}(y))f(x)dx\right) ^{p}dy   \nonumber \\
&<&(k^{(2)}(\sigma ))^{p}\int_{a_{1}}^{b_{1}}\frac{(v_{1}(x))^{p(1-\sigma
)-1}}{(v_{1}^{\prime }(x))^{p-1}}f^{p}(x)dx,  \label{39}
\end{eqnarray}%
where the constant factors $k^{(2)}(\sigma )$ and $(k^{(2)}(\sigma ))^{p}$
are the best possible.
\end{corollary}

Replacing $p>1$ by $0<p<1$ in the above inequalities, we have the equivalent
reverses of (\ref{38}) and (\ref{39}). If there exists a constant $\delta
_{\:0}>0,$ such that for any $\widetilde{\sigma }\in (\sigma -\delta
_{\:0},\sigma ],$
\[
k^{(2)}(\widetilde{\sigma })=\int_{1}^{\infty }h(t)t^{\widetilde{\sigma }%
-1}dt\in \mathbf{R}_{+},
\]
then the constant factors in the reverses of (\ref{38}) and (\ref{39}) are
the best possible.

Similarly, if $k_{\lambda }(1,t)=0\ (t>1),$ then $$k_{\lambda
}(x,y)=x^{-\lambda }k_{\lambda }\left(1,\frac{y}{x}\right)=0\ (y>x>0),$$ by Corollary 3,
we have the following First kind of Hardy-type integral inequalities with
the homogeneous kernel:

\begin{corollary}
 Suppose that $p>1,\:\frac{1}{p}+\frac{1}{q}%
=1,\:\mu ,\:\sigma \in \mathbf{R},\:\mu +\sigma =\lambda ,\:k_{\lambda }(x,y)$ is
a homogeneous function of degree $-\lambda $ in $\mathbf{R}_{+}^{2}$ $,$%
\begin{eqnarray*}
k_{\lambda }^{(1)}(\sigma )=\int_{0}^{1}k_{\lambda }(1,t)t^{\sigma -1}dt\in
\mathbf{R}_{+}.
\end{eqnarray*}
If $f(x),g(y)\geq 0,$
$$
0<\int_{0}^{\infty }x^{p(1-\mu )-1}f^{p}(x)dx<\infty$$ and $$0<\int_{0}^{\infty
}y^{q(1-\sigma )-1}g^{q}(y)dy<\infty,
$$
 then we have the following equivalent inequalities:%
\begin{eqnarray}
&&\int_{0}^{\infty }\left( \int_{y}^{\infty }k_{\lambda }(x,y)f(x)dx\right)
g(y)dy
=\int_{0}^{\infty }\left( \int_{x}^{\infty }k_{\lambda }(x,y)g(y)dy\right)
f(x)dx   \nonumber \\
&<&k_{\lambda }^{(1)}(\sigma )\left[ \int_{0}^{\infty }x^{p(1-\mu
)-1}f^{p}(x)dx\right] ^{\frac{1}{p}}\left[ \int_{0}^{\infty }y^{q(1-\sigma
)-1}g^{q}(y)dy\right] ^{\frac{1}{q}},  \label{40}
\end{eqnarray}%
\begin{equation}
\int_{0}^{\infty }y^{p\sigma -1}\left( \int_{y}^{\infty }k_{\lambda
}(x,y)f(x)dx\right) ^{p}dy
<(k_{\lambda }^{(1)}(\sigma ))^{p}\int_{0}^{\infty }x^{p(1-\mu
)-1}f^{p}(x)dx,  \label{41}
\end{equation}%
where, the constant factors $k_{\lambda }^{(1)}(\sigma )$ and $(k_{\lambda
}^{(1)}(\sigma ))^{p}$ are the best possible.
\end{corollary}

Replacing $p>1$ by $0<p<1$ in the above cases, we obtain the equivalent
reverses of (\ref{28}) and (\ref{29}). If there exists a constant $\delta
_{\:0}>0,$ such that for any $\widetilde{\sigma }\in (\sigma -\delta
_{\:0},\sigma ],$%
\[
k_{\lambda }^{(1)}(\widetilde{\sigma })=\int_{0}^{1}k_{\lambda }(1,t)t^{%
\widetilde{\sigma }-1}dt\in \mathbf{R}_{+},
\]
then the constant factors in the reverses of (\ref{40}) and (\ref{41}) are
the best possible.

If $k_{\lambda }(1,t)=0\ (t>1)$ in Corollary 2, then $$k_{\lambda
}(v_{1}(x),v_{2}(y))=0\ (0<v_{1}(x)<v_{2}(y)),$$ and we have the following general
results:

\begin{corollary}
 Suppose that $p>1,\:\frac{1}{p}+\frac{1}{q}=1,\:\mu
,\:\sigma \in \mathbf{R},\:\mu +\sigma =\lambda ,\:k_{\lambda }(x,y)$ is a
homogeneous function of degree $-\lambda $ in $\mathbf{R}_{+}^{2},$%
\begin{eqnarray*}
k_{\lambda }^{(1)}(\sigma )=\int_{0}^{1}k_{\lambda }(1,t)t^{\sigma -1}dt\in
\mathbf{R}_{+},
\end{eqnarray*}
$0\leq a_{i}<b_{i}\leq \infty ,\:v_{i}^{\prime }(s)>0$ $(s\in
(a_{i},b_{i})),v_{i}(a_{i}^{+})=0,\:v_{i}(b_{i}^{-})=\infty\ (i=1,2).$ If $%
f(x),g(y)\geq 0,$ such that%
$$
0 <\int_{a_{1}}^{b_{1}}\frac{(v_{1}(x))^{p(1-\mu )-1}}{(v_{1}^{\prime
}(x))^{p-1}}f^{p}(x)dx<\infty$$ and $$
0 <\int_{a_{2}}^{b_{2}}\frac{(v_{2}(y))^{q(1-\sigma )-1}}{(v_{2}^{\prime
}(y))^{q-1}}g^{q}(y)dy<\infty,
$$
then we have the following equivalent inequalities:%
\begin{eqnarray}
&&\int_{a_{2}}^{b_{2}}\left( \int_{v_{1}^{-1}(v_{2}(y))}^{b_{1}}k_{\lambda
}(v_{1}(x),v_{2}(y))f(x)dx\right) g(y)dy\nonumber \\
&=&\int_{a_{1}}^{b_{1}}\left( \int_{v_{2}^{-1}(v_{1}(x))}^{b_{2}}k_{\lambda
}(v_{1}(x),v_{2}(y))g(y)dy\right) f(x)dx\nonumber\\
&<&k_{\lambda }^{(1)}(\sigma )\left[ \int_{a_{1}}^{b_{1}}\frac{%
(v_{1}(x))^{p(1-\mu )-1}}{(v_{1}^{\prime }(x))^{p-1}}f^{p}(x)dx\right] ^{%
\frac{1}{p}}   \nonumber \\
&&\times \left[ \int_{a_{2}}^{b_{2}}\frac{(v_{2}(y))^{q(1-\sigma )-1}}{%
(v_{2}^{\prime }(y))^{q-1}}g^{q}(y)dy\right] ^{\frac{1}{q}},  \label{42}
\end{eqnarray}%
\begin{eqnarray}
&&\int_{a_{2}}^{b_{2}}\frac{v_{2}^{\prime }(y)}{(v_{2}(y))^{1-p\sigma }}%
\left( \int_{v_{1}^{-1}(v_{2}(y))}^{b_{1}}k_{\lambda
}(v_{1}(x),v_{2}(y))f(x)dx\right) ^{p}dy   \nonumber \\
&<&(k_{\lambda }^{(1)}(\sigma ))^{p}\int_{a_{1}}^{b_{1}}\frac{%
(v_{1}(x))^{p(1-\mu )-1}}{(v_{1}^{\prime }(x))^{p-1}}f^{p}(x)dx,  \label{43}
\end{eqnarray}%
where, the constant factors $k_{\lambda }^{(1)}(\sigma )$ and $(k_{\lambda
}^{(1)}(\sigma ))^{p}$ are the best possible.
\end{corollary}

Replacing $p>1$ by $0<p<1$ in the above cases, we have the equivalent
reverses of (\ref{42}) and (\ref{43}). If there exists a constant $\delta
_{\:0}>0,$ such that for any $\widetilde{\sigma }\in (\sigma -\delta
_{\:0},\sigma ],$%
\[
k_{\lambda }^{(1)}(\widetilde{\sigma })=\int_{0}^{1}k_{\lambda }(1,t)t^{%
\widetilde{\sigma }-1}dt\in \mathbf{R}_{+},
\]
then the constant factors in the reverses of (\ref{42}) and (\ref{43}) are
the best possible.

Similarly, if $k_{\lambda }(1,t)=0\ (0<t<1)$ in Corollary 3, then $$%
k_{\lambda }(x,y)=0\ (x>y>0),$$ and we have the following second kind Hardy-type
integral inequalities with the homogeneous kernel:

\begin{corollary}
 Suppose that $p>1,\:\frac{1}{p}+\frac{1}{q}=1,\:\mu
,\:\sigma \in \mathbf{R},\:\mu +\sigma =\lambda ,$ $k_{\lambda }(x,y)$ is a
homogeneous function of degree $-\lambda $ in $\mathbf{R}_{+}^{2},$%
\begin{eqnarray*}
k_{\lambda }^{(2)}(\sigma )=\int_{1}^{\infty }k_{\lambda }(1,t)t^{\sigma
-1}dt\in \mathbf{R}_{+}.
\end{eqnarray*}
If $f(x),g(y)\geq 0,$
$$
0<\int_{0}^{\infty }x^{p(1-\mu )-1}f^{p}(x)dx<\infty$$ and $$0<\int_{0}^{\infty
}y^{q(1-\sigma )-1}g^{q}(y)dy<\infty,
$$
then we have the following equivalent inequalities:%
$$
\int_{0}^{\infty }\left( \int_{0}^{y}k_{\lambda }(x,y)f(x)dx\right) g(y)dy
=\int_{0}^{\infty }\left( \int_{0}^{x}k_{\lambda }(x,y)g(y)dy\right) f(x)dx
$$
\begin{equation}
<k_{\lambda }^{(2)}(\sigma )\left[ \int_{0}^{\infty }x^{p(1-\mu
)-1}f^{p}(x)dx\right] ^{\frac{1}{p}}\left[ \int_{0}^{\infty }y^{q(1-\sigma
)-1}g^{q}(y)dy\right] ^{\frac{1}{q}},  \label{44}
\end{equation}%
\begin{equation}
\int_{0}^{\infty }y^{p\sigma -1}\left( \int_{0}^{y}k_{\lambda
}(x,y)f(x)dx\right) ^{p}dy
<(k_{\lambda }^{(2)}(\sigma ))^{p}\int_{0}^{\infty }x^{p(1-\mu
)-1}f^{p}(x)dx;  \label{45}
\end{equation}%
where, the constant factors $k_{\lambda }^{(2)}(\sigma )$ and $(k_{\lambda
}^{(2)}(\sigma ))^{p}$ are the best possible.
\end{corollary}

Replacing $p>1$ by $0<p<1$ in the above cases, we have the equivalent
reverses of (\ref{44}) and (\ref{45}). If there exists a constant $\delta
_{\:0}>0,$ such that for any $\widetilde{\sigma }\in (\sigma -\delta
_{\:0},\sigma ],$
\[
k_{\lambda }^{(2)}(\widetilde{\sigma })=\int_{1}^{\infty }k_{\lambda
}(1,t)t^{\widetilde{\sigma }-1}dt\in \mathbf{R}_{+},
\]
then the constant factors in the reverses of (\ref{44}) and (\ref{45}) are
the best possible.

If $k_{\lambda }(1,t)=0\ (0<t<1)$ in Corollary 2, then $$k_{\lambda
}(v_{1}(x),v_{2}(y))=0\ (v_{1}(x)>v_{2}(y)>0),$$ we have the following general
results:

\begin{corollary}
 Suppose that $p>1,\:\frac{1}{p}+\frac{1}{q}=1,\:\mu
,\sigma \in \mathbf{R},\:\mu +\sigma =\lambda ,\:k_{\lambda }(x,y)$ is a
homogeneous function of degree $-\lambda $ in $\mathbf{R}_{+}^{2},$%
\begin{eqnarray*}
k_{\lambda }^{(2)}(\sigma )=\int_{1}^{\infty }k_{\lambda }(1,t)t^{\sigma
-1}dt\in \mathbf{R}_{+},
\end{eqnarray*}
$0\leq a_{i}<b_{i}\leq \infty ,\:v_{i}^{\prime }(s)>0$ $(s\in
(a_{i},b_{i})),\:v_{i}(a_{i}^{+})=0,\:v_{i}(b_{i}^{-})=\infty\ (i=1,2).$ If $%
f(x),\:g(y)\geq 0,$ such that%
$$
0 <\int_{a_{1}}^{b_{1}}\frac{(v_{1}(x))^{p(1-\mu )-1}}{(v_{1}^{\prime
}(x))^{p-1}}f^{p}(x)dx<\infty$$ and $$
0 <\int_{a_{2}}^{b_{2}}\frac{(v_{2}(y))^{q(1-\sigma )-1}}{(v_{2}^{\prime
}(y))^{q-1}}g^{q}(y)dy<\infty,
$$
then we have the following equivalent inequalities:%
\begin{eqnarray}
&&\int_{a_{2}}^{b_{2}}\left( \int_{a_{1}}^{v_{1}^{-1}(v_{2}(y))}k_{\lambda
}(v_{1}(x),v_{2}(y))f(x)dx\right) g(y)dy\nonumber \\
&=&\int_{a_{1}}^{b_{1}}\left( \int_{a_{2}}^{v_{2}^{-1}(v_{1}(x))}k_{\lambda
}(v_{1}(x),v_{2}(y))g(y)dy\right) f(x)dx\nonumber\\
&<&k_{\lambda }^{(2)}(\sigma )\left[ \int_{a_{1}}^{b_{1}}\frac{%
(v_{1}(x))^{p(1-\mu )-1}}{(v_{1}^{\prime }(x))^{p-1}}f^{p}(x)dx\right] ^{%
\frac{1}{p}}   \nonumber \\
&&\times \left[ \int_{a_{2}}^{b_{2}}\frac{(v_{2}(y))^{q(1-\sigma )-1}}{%
(v_{2}^{\prime }(y))^{q-1}}g^{q}(y)dy\right] ^{\frac{1}{q}},  \label{46}
\end{eqnarray}%
\begin{eqnarray}
&&\int_{a_{2}}^{b_{2}}\frac{v_{2}^{\prime }(y)}{(v_{2}(y))^{1-p\sigma }}%
\left( \int_{a_{1}}^{v_{1}^{-1}(v_{2}(y))}k_{\lambda
}(v_{1}(x),v_{2}(y))f(x)dx\right) ^{p}dy   \nonumber \\
&<&(k_{\lambda }^{(2)}(\sigma ))^{p}\int_{a_{1}}^{b_{1}}\frac{%
(v_{1}(x))^{p(1-\mu )-1}}{(v_{1}^{\prime }(x))^{p-1}}f^{p}(x)dx,  \label{47}
\end{eqnarray}%
where, the constant factors $k_{\lambda }^{(2)}(\sigma )$ and $(k_{\lambda
}^{(2)}(\sigma ))^{p}$ are the best possible.
\end{corollary}
Replacing $p>1$ by $0<p<1$ in the above inequalities, we have the equivalent
reverses of (\ref{46}) and (\ref{47}). If there exists a constant $\delta
_{\:0}>0,$ such that for any $\widetilde{\sigma }\in (\sigma -\delta
_{\:0},\sigma ],$
\[
k_{\lambda }^{(2)}(\widetilde{\sigma })=\int_{1}^{\infty }k_{\lambda
}(1,t)t^{\widetilde{\sigma }-1}dt\in \mathbf{R}_{+},
\]
then the constant factors in the reverses of (\ref{46}) and (\ref{47}) are
the best possible.

\subsection{Yang-Hilbert-Type Operators and Hardy-Type Operators}

Suppose that $p>1,\:\frac{1}{p}+\frac{1}{q}=1,\:\mu ,\:\sigma \in \mathbf{R},\:\mu
+\sigma =\lambda .$ We set the following functions:%
\[
\varphi (x):=x^{p(1-\sigma)-1},\:\psi (y):=y^{q(1-\sigma )-1},\:\phi
(x):=x^{p(1-\mu )-1}(x,y\in \mathbf{R}_{+}),
\]
from which we obtain that $\psi ^{1-p}(y)=y^{p\sigma -1}.$

 Define the following real normed
linear space:
\[
L_{p,\varphi }(\mathbf{R}_{+}):=\left\{f:||f||_{p,\varphi }:=\left\{\int_{0}^{\infty
}\varphi (x)|f(x)|^{p}dx\right\}^{\frac{1}{p}}<\infty\right\}.
\]
Therefore,%
\begin{eqnarray*}
L_{p,\psi ^{1-p}}(\mathbf{R}_{+}) &=&\left\{h:||h||_{p,\psi
^{1-p}}:=\left\{\int_{0}^{\infty }\psi ^{1-p}(y)|h(y)|^{p}dy\right\}^{\frac{1}{p}%
}<\infty\right\}, \\
L_{p,\phi}(\mathbf{R}_{+}) &=&\left\{g:||g||_{p,\phi }:=\left\{\int_{0}^{\infty
}\phi(x)|g(x)|^{p}dx\right\}^{\frac{1}{p}}<\infty\right\}.
\end{eqnarray*}

(a) In view of Theorem 1, for $f\in L_{p,\varphi }(\mathbf{R}_{+}),$ setting
$$H_{1}(y):=\int_{0}^{\infty }h(xy)|f(x)|dx^{{}}\ \ (y\in \mathbf{R}_{+}), $$
by (\ref{19}), we have
\begin{equation}
||H_{1}||_{p,\psi ^{1-p}} :=\left( \int_{0}^{\infty }\psi
^{1-p}(y)H_{1}^{p}(y)dy\right) ^{\frac{1}{p}} <k(\sigma )||f||_{p,\varphi
}<\infty .  \label{48}
\end{equation}

\begin{definition}
 Let us define the Yang-Hilbert-type integral operator with the
non-homogeneous kernel $T_{1}:$ $L_{p,\varphi }(\mathbf{R}_{+})\rightarrow
L_{p,\psi ^{1-p}}(\mathbf{R}_{+})$ as follows:\\ For any $f\in L_{p,\varphi }(%
\mathbf{R}_{+}),$ there exists a unique representation $T_{1}f=H_{1}\in
L_{p,\psi ^{1-p}}(\mathbf{R}_{+}),$ satisfying
$$T_{1}f(y)=H_{1}(y),$$
for any $y\in \mathbf{R}%
_{+}$.
\end{definition}

In view of (\ref{48}), it follows that $$||T_{1}f||_{p,\psi
^{1-p}}=||H_{1}||_{p,\psi ^{1-p}}\leq k(\sigma )||f||_{p,\varphi } $$ and
then the operator $T_{1}$ is bounded satisfying%
\[
||T_{1}||=\sup_{f(\neq \theta )\in L_{p,\varphi }(\mathbf{R}_{+})}\frac{%
||T_{1}f||_{p,\psi ^{1-p}}}{||f||_{p,\varphi }}\leq k(\sigma ).
\]
Since the constant factor $k(\sigma )$ in (\ref{48}) is the best possible,
we have%
\begin{equation}
||T_{1}||=k(\sigma )=\int_{0}^{\infty }h(t)t^{\sigma -1}dt.  \label{49}
\end{equation}

If we define the formal inner product of $T_{1}f$ and $g$ as
\[
(T_{1}f,g) :=\int_{0}^{\infty }\left( \int_{0}^{\infty }h(xy)f(x)dx\right)
g(y)dy =\int_{0}^{\infty }\int_{0}^{\infty }h(xy)f(x)g(y)dxdy,
\]
then we can rewrite (\ref{18}) and (\ref{19}) as follows:%
\[
(T_{1}f,g)<||T_{1}||\cdot ||f||_{p,\varphi }||g||_{q,\psi
},\ \ ||T_{1}f||_{p,\psi ^{1-p}}<||T_{1}||\cdot ||f||_{p,\varphi }.
\]

(b) In view of Corollary 4, for $f\in L_{p,\varphi }(\mathbf{R}_{+}),$
setting $$H_{1}^{(1)}(y):=\int_{0}^{\frac{1}{y}}h(xy)|f(x)|dx^{{}}\ \ (y\in
\mathbf{R}_{+}),$$ by (\ref{32}), we obtain
\begin{equation}
||H_{1}^{(1)}||_{p,\psi ^{1-p}} :=\left( \int_{0}^{\infty }\psi
^{1-p}(y)(H_{1}^{(1)}(y))^{p}dy\right) ^{\frac{1}{p}} <k^{(1)}(\sigma
)||f||_{p,\varphi }<\infty .  \label{50}
\end{equation}

\begin{definition}
 Define the Hardy-type integral operator of the first kind
with the non-\\homogeneous kernel $$T_{1}^{(1)}\::\:L_{p,\varphi }(\mathbf{R}%
_{+})\rightarrow L_{p,\psi ^{1-p}}(\mathbf{R}_{+})$$ as follows:

For any $%
f\in L_{p,\varphi }(\mathbf{R}_{+}),$ there exists a unique representation $%
T_{1}^{(1)}f=H_{1}^{(1)}\in L_{p,\psi ^{1-p}}(\mathbf{R}_{+}),$ satisfying
$$T_{1}^{(1)}f(y)=H_{1}^{(1)}(y),$$
for any $y\in \mathbf{R}_{+}.$
\end{definition}

In view of (\ref{50}), it follows that $$||T_{1}^{(1)}f||_{p,\psi
^{1-p}}=||H_{1}^{(1)}||_{p,\psi ^{1-p}}\leq k^{(1)}(\sigma )||f||_{p,\varphi
} $$ and thus the operator $T_{1}^{(1)}$ is bounded satisfying%
\[
||T_{1}^{(1)}||=\sup_{f(\neq \theta )\in L_{p,\varphi }(\mathbf{R}_{+})}%
\frac{||T_{1}^{(1)}f||_{p,\psi ^{1-p}}}{||f||_{p,\varphi }}\leq
k^{(1)}(\sigma ).
\]
Since the constant factor $k^{(1)}(\sigma )$ in (\ref{50}) is the best
possible, we have%
\begin{equation}
||T_{1}^{(1)}||=k^{(1)}(\sigma )=\int_{0}^{1}h(t)t^{\sigma -1}dt.  \label{51}
\end{equation}

Setting the formal inner product of $T_{1}^{(1)}f$ and $g$ as
\[
(T_{1}^{(1)}f,g)=\int_{0}^{\infty }\left( \int_{0}^{\frac{1}{y}%
}h(xy)f(x)dx\right) g(y)dy,
\]
we can rewrite (\ref{31}) and (\ref{32}) as follows:%
\begin{equation}
(T_{1}^{(1)}f,g)<||T_{1}^{(1)}||\cdot ||f||_{p,\varphi }||g||_{q,\psi
},\ \ ||T_{1}^{(1)}f||_{p,\psi ^{1-p}}<||T_{1}^{(1)}||\cdot ||f||_{p,\varphi }.
\label{52}
\end{equation}

(c) In view of Corollary 6, for $f\in L_{p,\varphi }(\mathbf{R}_{+}),$
setting $$H_{1}^{(2)}(y):=\int_{\frac{1}{y}}^{\infty }h(xy)|f(x)|dx\ \ (y\in
\mathbf{R}_{+}),$$ by (\ref{37}), we have
\begin{equation}
||H_{1}^{(2)}||_{p,\psi ^{1-p}} :=\left( \int_{0}^{\infty }\psi
^{1-p}(y)(H_{1}^{(2)}(y))^{p}dy\right) ^{\frac{1}{p}} <k^{(2)}(\sigma
)||f||_{p,\varphi }<\infty .  \label{53}
\end{equation}

\begin{definition}
 Define the second kind Hardy-type integral operator
with the non-homogeneous kernel $$T_{1}^{(2)}\::\:L_{p,\varphi }(\mathbf{R}%
_{+})\rightarrow L_{p,\psi ^{1-p}}(\mathbf{R}_{+})$$ as follows:

For any $%
f\in L_{p,\varphi }(\mathbf{R}_{+}),$ there exists a unique representation $%
T_{1}^{(2)}f=H_{1}^{(2)}\in L_{p,\psi ^{1-p}}(\mathbf{R}_{+}),$ satisfying
$$T_{1}^{(2)}f(y)=H_{1}^{(2)}(y),$$
for any $y\in \mathbf{R}_{+}.$
\end{definition}

In view of (\ref{37}), it follows that
\[
||T_{1}^{(2)}f||_{p,\psi ^{1-p}}=||H_{1}^{(2)}||_{p,\psi ^{1-p}}\leq
k^{(2)}(\sigma )||f||_{p,\varphi }
\]
and hence the operator $T_{1}^{(2)}$ is bounded satisfying%
\[
||T_{1}^{(2)}||=\sup_{f(\neq \theta )\in L_{p,\varphi }(\mathbf{R}_{+})}%
\frac{||T_{1}^{(2)}f||_{p,\psi ^{1-p}}}{||f||_{p,\varphi }}\leq
k^{(2)}(\sigma ).
\]
Since the constant factor $k^{(2)}(\sigma )$ in (\ref{53}) is the best
possible, we have%
\begin{equation}
||T_{1}^{(2)}||=k^{(2)}(\sigma )=\int_{1}^{\infty }h(t)t^{\sigma -1}dt.
\label{54}
\end{equation}

Setting the formal inner product of $T_{1}^{(2)}f$ and $g$ as
\[
(T_{1}^{(2)}f,g)=\int_{0}^{\infty }\left( \int_{\frac{1}{y}}^{\infty
}h(xy)f(x)dx\right) g(y)dy,
\]
we can rewrite (\ref{36}) and (\ref{37}) as follows:%
\begin{equation}
(T_{1}^{(2)}f,g)<||T_{1}^{(2)}||\cdot ||f||_{p,\varphi }||g||_{q,\psi
},\ \ ||T_{1}^{(2)}f||_{p,\psi ^{1-p}}<||T_{1}^{(2)}||\cdot ||f||_{p,\varphi }.
\label{55}
\end{equation}

(d) In view of Corollary 3, for $f\in L_{p,\phi }(\mathbf{R}_{+}),$ setting $$%
H_{2}(y):=\int_{0}^{\infty }k_{\lambda }(x,y)|f(x)|dx\ \ (y\in \mathbf{R}%
_{+}),$$ by (\ref{29}), we have
\begin{equation}
||H_{2}||_{p,\psi ^{1-p}} :=\left( \int_{0}^{\infty }\psi
^{1-p}(y)H_{2}^{p}(y)dy\right) ^{\frac{1}{p}} <k_{\lambda }(\sigma
)||f||_{p,\phi }<\infty .  \label{56}
\end{equation}
\begin{definition}
 Define the Yang-Hilbert-type integral operator with the
homogeneous kernel $T_{2}:$ $L_{p,\phi }(\mathbf{R}_{+})\rightarrow
L_{p,\psi ^{1-p}}(\mathbf{R}_{+})$ as follows:

 For any $f\in L_{p,\phi }(%
\mathbf{R}_{+}),$ there exists a unique representation $T_{2}f=H_{2}\in
L_{p,\psi^{1-p}}(\mathbf{R}_{+}),$ satisfying $$T_{2}f(y)=H_{2}(y),$$
for any $y\in \mathbf{R}%
_{+}.$
\end{definition}

In view of (\ref{56}), it follows that
\[
||T_{2}f||_{p,\psi ^{1-p}}=||H_{2}||_{p,\psi ^{1-p}}\leq k_{\lambda }(\sigma
)||f||_{p,\phi }
\]
and thus the operator $T_{2}$ is bounded satisfying%
\[
||T_{2}||=\sup_{f(\neq \theta )\in L_{p,\phi }(\mathbf{R}_{+})}\frac{%
||T_{2}f||_{p,\psi ^{1-p}}}{||f||_{p,\phi }}\leq k_{\lambda }(\sigma ).
\]
Since the constant factor $k_{\lambda }(\sigma )$ in (\ref{56}) is the best
possible, we have%
\begin{equation}
||T_{2}||=k_{\lambda }(\sigma )=\int_{0}^{\infty }k_{\lambda }(1,t)t^{\sigma
-1}dt.  \label{57}
\end{equation}

Setting the formal inner product of $T_{2}f$ and $g$ as
\[
(T_{2}f,g)=\int_{0}^{\infty }\left( \int_{0}^{\infty }k_{\lambda
}(x,y)f(x)dx\right) g(y)dy =\int_{0}^{\infty }\int_{0}^{\infty }k_{\lambda
}(x,y)f(x)g(y)dxdy,
\]
we can rewrite (\ref{28}) and (\ref{29}) as follows:%
\[
(T_{2}f,g)<||T_{2}||\cdot ||f||_{p,\phi }||g||_{q,\psi },\ \ ||T_{2}f||_{p,\psi
^{1-p}}<||T_{2}||\cdot ||f||_{p,\phi }.
\]

(e) Due to Corollary 8, for $f\in L_{p,\phi}(\mathbf{R}_{+}),$  $$%
H_{2}^{(1)}(y):=\int_{y}^{\infty }k_{\lambda }(x,y)|f(x)|dx\ \ (y\in \mathbf{R%
}_{+}), $$ by (\ref{41}), we have
\begin{equation}
||H_{2}^{(1)}||_{p,\psi ^{1-p}} :=\left( \int_{0}^{\infty }\psi
^{1-p}(y)(H_{2}^{(1)}(y))^{p}dy\right) ^{\frac{1}{p}} <k_{\lambda
}^{(1)}(\sigma )||f||_{p,\phi }<\infty .  \label{58}
\end{equation}

\begin{definition}
 Define the Hardy-type integral operator of the first kind,
with the homogeneous kernel $T_{2}^{(1)}:$ $L_{p,\phi}(\mathbf{R}%
_{+})\rightarrow L_{p,\psi ^{1-p}}(\mathbf{R}_{+})$ as follows:

For any $%
f\in L_{p,\phi }(\mathbf{R}_{+}),$ there exists a unique representation $%
T_{2}^{(1)}f=H_{2}^{(1)}\in L_{p,\psi^{1-p}}(\mathbf{R}_{+}),$ satisfying
$$T_{2}^{(1)}f(y)=H_{2}^{(1)}(y),$$
for any $y\in \mathbf{R}_{+}.$
\end{definition}

By (\ref{41}), it follows that%
\[
||T_{2}^{(1)}f||_{p,\psi^{1-p}}=||H_{2}^{(1)}||_{p,\psi ^{1-p}}\leq
k_{\lambda }^{(1)}(\sigma )||f||_{p,\phi }
\]
and then the operator $T_{2}^{(1)}$ is bounded satisfying%
\[
||T_{2}^{(1)}||=\sup_{f(\neq \theta )\in L_{p,\phi }(\mathbf{R}_{+})}\frac{%
||T_{2}^{(1)}f||_{p,\psi^{1-p}}}{||f||_{p,\phi}}\leq k_{\lambda
}^{(1)}(\sigma ).
\]
Since the constant factor $k_{\lambda }^{(1)}(\sigma )$ in (\ref{58}) is the
best possible, we have%
\begin{equation}
||T_{2}^{(1)}||=k_{\lambda }^{(1)}(\sigma )=\int_{0}^{1}k_{\lambda
}(1,t)t^{\sigma -1}dt.  \label{59}
\end{equation}

Setting the formal inner product of $T_{2}^{(1)}f$ and $g$ as
\[
(T_{2}^{(1)}f,g)=\int_{0}^{\infty }\left( \int_{y}^{\infty }k_{\lambda
}(x,y)f(x)dx\right) g(y)dy,
\]
we can rewrite (\ref{40}) and (\ref{41}) as follows:%
\[
(T_{2}^{(1)}f,g)<||T_{2}^{(1)}||\cdot ||f||_{p,\phi}||g||_{q,\psi
},\ \ ||T_{2}^{(1)}f||_{p,\psi ^{1-p}}<||T_{2}^{(1)}||\cdot ||f||_{p,\phi }.
\]

\bigskip (f) By Corollary 10, for $f\in L_{p,\phi }(\mathbf{R}_{+}),$
$$H_{2}^{(2)}(y):=\int_{0}^{y}k_{\lambda }(x,y)|f(x)|dx\ \ (y\in \mathbf{R}%
_{+}), $$ and by (\ref{45}), we have
\begin{equation}
||H_{2}^{(2)}||_{p,\psi ^{1-p}} :=\left( \int_{0}^{\infty }\psi
^{1-p}(y)(H_{2}^{(2)}(y))^{p}dy\right) ^{\frac{1}{p}} <k_{\lambda
}^{(2)}(\sigma )||f||_{p,\phi }<\infty .  \label{60}
\end{equation}

\begin{definition}
 Define the Hardy-type integral operator of the second kind
with the homogeneous kernel $T_{2}^{(2)}\::\:L_{p,\phi }(\mathbf{R}%
_{+})\rightarrow L_{p,\psi^{1-p}}(\mathbf{R}_{+})$ as follows:

 For any $%
f\in L_{p,\phi}(\mathbf{R}_{+}),$ there exists a unique representation $%
T_{2}^{(2)}f=H_{2}^{(2)}\in L_{p,\psi^{1-p}}(\mathbf{R}_{+}),$ satisfying
$$T_{2}^{(2)}f(y)=H_{2}^{(2)}(y),$$
for any $y\in \mathbf{R}_{+}$.
\end{definition}

In view of (\ref{45}), it follows that%
\[
||T_{2}^{(2)}f||_{p,\psi ^{1-p}}=||H_{2}^{(2)}||_{p,\psi ^{1-p}}\leq
k_{\lambda }^{(2)}(\sigma )||f||_{p,\phi }
\]
and then the operator $T_{2}^{(2)}$ is bounded satisfying%
\[
||T_{2}^{(2)}||=\sup_{f(\neq \theta )\in L_{p,\phi }(\mathbf{R}_{+})}\frac{%
||T_{2}^{(2)}f||_{p,\psi ^{1-p}}}{||f||_{p,\phi}}\leq k_{\lambda
}^{(2)}(\sigma ).
\]
Since the constant factor $k_{\lambda }^{(2)}(\sigma )$ in (\ref{60}) is the
best possible, we have%
\begin{equation}
||T_{2}^{(2)}||=k_{\lambda }^{(2)}(\sigma )=\int_{1}^{\infty }k_{\lambda
}(1,t)t^{\sigma -1}dt.  \label{61}
\end{equation}
Setting the formal inner product of $T_{2}^{(2)}f$ and $g$ as
\[
(T_{2}^{(2)}f,g)=\int_{0}^{\infty }\left( \int_{0}^{y}k_{\lambda
}(x,y)f(x)dx\right) g(y)dy,
\]
we can rewrite (\ref{44}) and (\ref{45}) as follows:%
\[
(T_{2}^{(2)}f,g)<||T_{2}^{(2)}||\cdot ||f||_{p,\phi}||g||_{q,\psi
},||T_{2}^{(2)}f||_{p,\psi ^{1-p}}<||T_{2}^{(2)}||\cdot ||f||_{p,\phi }.
\]

\subsection{ Some Examples}

\begin{example}
 (a) Set $$h(t)=k_{\lambda }(1,t)=\frac{1}{%
(1+t)^{\lambda }}\ \ (\mu ,\sigma >0,\mu +\sigma =\lambda ).$$ Then we have the
kernels $$
h(xy)=\frac{1}{(1+xy)^{\lambda }},k_{\lambda }(x,y)=\frac{1}{(x+y)^{\lambda }%
}
$$
and obtain the constant factors
$$
k(\sigma )=k_{\lambda }(\sigma )=\int_{0}^{\infty }\frac{t^{\sigma -1}}{%
(1+t)^{\lambda }}dt=B(\mu ,\sigma )\in \mathbf{R}_{+}.
$$
 By (\ref{49}) and (\ref{57}), we have $||T_{1}||=||T_{2}||=B(\mu
,\sigma ).$
\end{example}

(b) Set $$h(t)=k_{\lambda }(1,t)=\frac{-\ln t}{1-t^{\lambda }}\ (\mu ,\sigma
>0,\mu +\sigma =\lambda ).$$ Then we have the kernels $$h(xy)=\frac{\ln (xy)}{%
(xy)^{\lambda }-1},\ \ k_{\lambda }(x,y)=\frac{\ln (x/y)}{x^{\lambda
}-y^{\lambda }} $$ and obtain the constant factors%
\begin{eqnarray*}
k(\sigma )&=&k_{\lambda }(\sigma )=\int_{0}^{\infty }\frac{(\ln t)t^{\sigma
-1}}{t^{\lambda }-1}dt \\
&=&\frac{1}{\lambda ^{2}}\int_{0}^{\infty }\frac{(\ln u)u^{(\sigma /\lambda
)-1}}{u-1}du=\left[\frac{\pi }{\lambda \sin \pi (\sigma /\lambda )}\right]^{2}\in
\mathbf{R}_{+}.
\end{eqnarray*}
In view of (\ref{49}) and (\ref{57}), we have $$||T_{1}||=||T_{2}||=[\frac{%
\pi }{\lambda \sin \pi (\sigma /\lambda )}]^{2}.$$

(c) Set $$h(t)=k_{\lambda }(1,t)=\frac{|\ln t|^{\beta }}{(\max
\{1,t\})^{\lambda }}\ \ (\beta >-1,\mu ,\sigma >0,\mu +\sigma =\lambda ).$$ Then
we have the kernels $$h(xy)=\frac{|\ln (xy)|^{\beta }}{(\max
\{1,xy\})^{\lambda }},\ \ k_{\lambda }(x,y)=\frac{|\ln (x/y)|^{\beta }}{(\max
\{x,y\})^{\lambda }} $$ and by using
the formula (cf. \cite{YM6a}) $$\Gamma(\alpha):=\int_0^\infty e^{-t}
t^{\alpha-1}dt\ \ (\alpha>0)$$
we obtain the following constant factors
\begin{eqnarray*}
k(\sigma ) &=&k_{\lambda }(\sigma )=\int_{0}^{\infty }\frac{|\ln t|^{\beta
}t^{\sigma -1}}{(\max \{1,t\})^{\lambda }}dt \\
&=&\int_{0}^{1}(-\ln t)^{\beta }t^{\sigma -1}dt+\int_{1}^{\infty }\frac{(\ln
t)^{\beta }t^{\sigma -1}}{t^{\lambda }}dt \\
&=&\int_{0}^{1}(-\ln t)^{\beta }(t^{\sigma -1}+t^{\mu -1})dt=\left(\frac{1}{%
\sigma ^{\beta +1}}\frac{1}{\mu ^{\beta +1}}\right)\int_{0}^{\infty }v^{\beta
}e^{-v}dv \\
&=&\left(\frac{1}{\sigma ^{\beta +1}}+\frac{1}{\mu ^{\beta +1}}\right)\Gamma (\beta
+1)\in \mathbf{R}_{+}.
\end{eqnarray*}%
By (\ref{49}) and (\ref{57}), we have $$||T_{1}||=||T_{2}||=\left(\frac{1}{%
\sigma ^{\beta +1}}+\frac{1}{\mu ^{\beta +1}}\right)\Gamma (\beta +1). $$ Due to (\ref%
{51}) and (\ref{59}), we have $$||T_{1}^{(1)}||=||T_{2}^{(1)}||=\frac{1}{%
\sigma ^{\beta +1}}\Gamma (\beta +1), $$ and by (\ref{54}) and (\ref{61}), it
follows that $$||T_{1}^{(2)}||=||T_{2}^{(2)}||=\frac{1}{\mu ^{\beta +1}}%
\Gamma (\beta +1). $$

(d) Set $$h(t)=k_{\lambda }(1,t)=\frac{|\ln t|^{\beta }}{(\min
\{1,t\})^{\lambda }}\ \ (\beta >-1,\mu ,\sigma <0,\mu +\sigma =\lambda ).$$ Then
we have the kernels $$h(xy)=\frac{|\ln (xy)|^{\beta }}{(\min
\{1,xy\})^{\lambda }},\ \ k_{\lambda }(x,y)=\frac{|\ln (x/y)|^{\beta }}{(\min
\{x,y\})^{\lambda }}$$ and obtain the constant factors%
\begin{eqnarray*}
k(\sigma ) &=&k_{\lambda }(\sigma )=\int_{0}^{\infty }\frac{|\ln t|^{\beta
}t^{\sigma -1}}{(\min \{1,t\})^{\lambda }}dt \\
&=&\int_{0}^{1}\frac{(-\ln t)^{\beta }t^{\sigma -1}}{t^{\lambda }}%
dt+\int_{1}^{\infty }(\ln t)^{\beta }t^{\sigma -1}dt \\
&=&\int_{0}^{1}(-\ln t)^{\beta }(t^{-\mu -1}+t^{-\sigma -1})dt=\left[\frac{1}{%
(-\mu )^{\beta +1}}\frac{1}{(-\sigma )^{\beta +1}}\right]\int_{0}^{\infty
}v^{\beta }e^{-v}dv \\
&=&\left[\frac{1}{(-\mu )^{\beta +1}}+\frac{1}{(-\sigma )^{\beta +1}}\right]\Gamma
(\beta +1)\in \mathbf{R}_{+}.
\end{eqnarray*}%
In view of (\ref{49}) and (\ref{57}), we have $$||T_{1}||=||T_{2}||=\left[\frac{1}{%
(-\mu )^{\beta +1}}+\frac{1}{(-\sigma )^{\beta +1}}\right]\Gamma (\beta +1).$$ By (%
\ref{51}) and (\ref{59}), we have $$||T_{1}^{(1)}||=||T_{2}^{(1)}||=\frac{1}{%
(-\mu )^{\beta +1}}\Gamma (\beta +1), $$ and by (\ref{54}) and (\ref{61}), it
follows that $$||T_{1}^{(2)}||=||T_{2}^{(2)}||=\frac{1}{(-\sigma )^{\beta +1}}%
\Gamma (\beta +1).$$

(e) Set $$h(t)=k_{\lambda }(1,t)=\frac{|\ln t|}{1+t^{\lambda }}\ \ (\mu
,\sigma >0,\mu +\sigma =\lambda ).$$ Then we have the kernels $$h(xy)=\frac{%
|\ln (xy)|}{1+(xy)^{\lambda }},\ \ k_{\lambda }(x,y)=\frac{|\ln (x/y)|}{%
x^{\lambda }+y^{\lambda }} $$ and obtain the constant factors%
\begin{eqnarray*}
&&k(\sigma )=k_{\lambda }(\sigma )=\int_{0}^{\infty }\frac{|\ln t|t^{\sigma
-1}}{1+t^{\lambda }}dt \\
&=&\int_{0}^{1}\frac{(-\ln t)t^{\sigma -1}}{t^{\lambda }+1}%
dt+\int_{1}^{\infty }\frac{(\ln t)t^{\sigma -1}}{t^{\lambda }+1}%
dt=\int_{0}^{1}\frac{(-\ln t)(t^{\sigma -1}+t^{\mu -1})}{t^{\lambda }+1}dt \\
&=&\int_{0}^{1}(-\ln t)\sum_{k=0}^{\infty }(-1)^{k}(t^{k\lambda +\sigma
-1}+t^{k\lambda +\mu -1})dt.
\end{eqnarray*}
By the fact that
\begin{eqnarray*}
&&\sum_{k=0}^{\infty }\int_{0}^{1}|(-1)^{k}(-\ln t)(t^{k\lambda +\sigma
-1}+t^{k\lambda +\mu -1})|dt \\
&=&\sum_{k=0}^{\infty }\int_{0}^{1}(-\ln t)(\frac{1}{k\lambda +\sigma }%
dt^{k\lambda +\sigma }+\frac{1}{k\lambda +\mu }dt^{k\lambda +\mu }) \\
&=&\sum_{k=0}^{\infty }\left[\frac{1}{(k\lambda +\sigma )^{2}}+\frac{1}{(k\lambda
+\mu )^{2}}\right]\in \mathbf{R}_{+},
\end{eqnarray*}%
in combination with Theorem 7 (cf. \cite{YM25}, Chapter 5), we obtain%
\begin{eqnarray*}
k(\sigma ) &=&k_{\lambda }(\sigma )=\int_{0}^{1}(-\ln t)\sum_{k=0}^{\infty
}(-1)^{k}(t^{k\lambda +\sigma -1}+t^{k\lambda +\mu -1})dt \\
&=&\sum_{k=0}^{\infty }\int_{0}^{1}(-1)^{k}(-\ln t)(t^{k\lambda +\sigma
-1}+t^{k\lambda +\mu -1})dt \\
&=&\sum_{k=0}^{\infty }(-1)^{k}\left[\frac{1}{(k\lambda +\sigma )^{2}}+\frac{1}{%
(k\lambda +\mu )^{2}}\right]\in \mathbf{R}_{+}.
\end{eqnarray*}%
By (\ref{49}) and (\ref{57}), we have $$||T_{1}||=||T_{2}||=%
\sum_{k=0}^{\infty }(-1)^{k}\left[\frac{1}{(k\lambda +\sigma )^{2}}+\frac{1}{%
(k\lambda +\mu )^{2}}\right]. $$ By (\ref{51}) and (\ref{59}), we have $$%
||T_{1}^{(1)}||=||T_{2}^{(1)}||=\sum_{k=0}^{\infty }(-1)^{k}\frac{1}{%
(k\lambda +\sigma )^{2}}, $$ and by (\ref{54}) and (\ref{61}), it follows
that $$||T_{1}^{(2)}||=||T_{2}^{(2)}||=\sum_{k=0}^{\infty }(-1)^{k}\frac{1}{%
(k\lambda +\mu )^{2}}. $$

(f) Set $$h(t)=k_{\lambda }(t,1)=\frac{1}{|1-t|^{\lambda }}\ \ (\ \mu ,\sigma
>0,\mu +\sigma =\lambda <1).$$ Then, we have kernels $$h(xy)=\frac{1}{%
|1-xy|^{\lambda }},\ \ k_{\lambda }(x,y)=\frac{1}{|x-y|^{\lambda }} $$ and obtain
the constant factors%
\begin{eqnarray*}
k(\sigma )&=&k_{\lambda }(\sigma )=\int_{0}^{\infty }\frac{t^{\sigma -1}}{%
|1-t|^{\lambda }}dt \\
&=&\int_{0}^{1}\frac{t^{\sigma -1}+t^{\mu -1}}{(1-t)^{\lambda }}dt
=B(1-\lambda ,\sigma )+B(1-\lambda ,\mu )\in \mathbf{R}_{+}.
\end{eqnarray*}
In view of (\ref{49}) and (\ref{57}), we obtain $$||T_{1}||=||T_{2}||=B(1-%
\lambda ,\sigma )+B(1-\lambda ,\mu ). $$ By (\ref{51}) and (\ref{59}), we
have $$||T_{1}^{(1)}||=||T_{2}^{(1)}||=B(1-\lambda ,\sigma ), $$ and by (\ref%
{54}) and (\ref{61}), it follows that $$||T_{1}^{(2)}||=||T_{2}^{(2)}||=B(1-%
\lambda ,\mu ). $$

For (a)-(f), we can obtain the equivalent inequalities with the kernels
and the best possible constant factors in Theorems 1-4. Setting $\delta _{\:0}=%
\frac{|\sigma |}{2}>0,$ we can still obtain the equivalent reverse
inequalities with the kernels and the best possible constant factors in
Theorems 1-4.

(g) Set $$h(t)=k_{\lambda }(1,t)=\frac{(\min \{t,1\})^{\eta }}{(\max
\{t,1\})^{\lambda +\eta }}\ \ (\eta >-\min \{\mu ,\sigma \},\mu +\sigma =\lambda
).$$ Then we have the kernels $$h(xy)=\frac{(\min \{1,xy\})^{\eta }}{(\max
\{1,xy\})^{\lambda +\eta }},\ \ k_{\lambda }(x,y)=\frac{(\min \{x,y\})^{\eta }}{%
(\max \{x,y\})^{\lambda +\eta }} $$ and obtain the constant factors%
\begin{eqnarray*}
k(\sigma ) &=&k_{\lambda }(\sigma )=\int_{0}^{\infty }\frac{(\min
\{1,t\})^{\eta }t^{\sigma -1}}{(\max \{1,t\})^{\lambda +\eta }}dt \\
&=&\int_{0}^{1}t^{\eta }t^{\sigma -1}dt+\int_{1}^{\infty }\frac{t^{\sigma -1}%
}{t^{\lambda +\eta }}dt=\frac{1}{\sigma +\eta }+\frac{1}{\mu +\eta } \\
&=&\frac{\lambda +2\eta }{(\sigma +\eta )(\mu +\eta )}\in \mathbf{R}_{+}.
\end{eqnarray*}%
In view of (\ref{49}) and (\ref{57}), we get $$||T_{1}||=||T_{2}||=\frac{%
\lambda +2\eta }{(\sigma +\eta )(\mu +\eta )}. $$ By (\ref{51}) and (\ref{59}%
), we have $$||T_{1}^{(1)}||=||T_{2}^{(1)}||=\frac{1}{\sigma +\eta }, $$ and
by (\ref{54}) and (\ref{61}), it follows that $$%
||T_{1}^{(2)}||=||T_{2}^{(2)}||=\frac{1}{\mu +\eta }. $$

Then we can derive the equivalent inequalities with the kernels and the best
possible constant factors in Theorems 1-4. Setting $\delta _{\:0}=\frac{\sigma
+\eta }{2}>0,$ we can still obtain the equivalent reverse inequalities with
the kernels and the best possible constant factors in Theorems 1-4.

In particular, (i) for $\eta =0$,$$h(t)=k_{\lambda }(1,t)=\frac{1}{(\max
\{1,t\})^{\lambda }}\ \ (\ \mu ,\sigma >0,\mu +\sigma =\lambda ),$$ we have $$%
h(xy)=\frac{1}{(\max \{1,xy\})^{\lambda }},\ \ k_{\lambda }(x,y)=\frac{1}{(\max
\{x,y\})^{\lambda }} $$ and $$||T_{1}||=||T_{2}||=\frac{\lambda }{\sigma \mu };
$$

(ii) for $\eta =-\lambda ,$ $$h(t)=k_{\lambda }(t,1)=\frac{1}{(\min
\{1,t\})^{\lambda }}\ \ (\ \mu ,\sigma <0,\mu +\sigma =\lambda ),$$ we have  $$%
h(xy)=\frac{1}{(\min \{1,xy\})^{\lambda }},\ \ k_{\lambda }(x,y)=\frac{1}{(\min
\{x,y\})^{\lambda }}$$ and $$||T_{1}||=||T_{2}||=\frac{-\lambda }{\sigma \mu };
$$

(iii) for $\lambda =0,$ $$h(t)=k_{0}(1,t)=(\frac{\min \{1,t\}}{\max \{1,t\}}%
)^{\eta }(\eta >|\sigma |),$$ we have%
\[
h(xy)=\left( \frac{\min \{1,xy\}}{\max \{1,xy\}}\right) ^{\eta },k_{\lambda
}(x,y)=\left( \frac{\min \{x,y\}}{\max \{x,y\}}\right) ^{\eta }
\]
and $$||T_{1}||=||T_{2}||=\frac{2\eta }{\eta ^{2}-\sigma ^{2}}.$$

\begin{example}
 (a) Set $$h(t)=k_{0}(1,t)=\ln (1+\frac{\rho }{%
t^{\eta }})\ \ (\rho >0,0<\sigma <\eta ).$$ Then we have the kernels $$
h(xy)=\ln [1+\frac{\rho }{(xy)^{\eta }}],\ \ k_{0}(x,y)=\ln \left[1+\rho \left(\frac{x}{y}%
\right)^{\eta }\right]
$$
and obtain the constant factors%
\begin{eqnarray*}
k(\sigma ) &=&k_{0}(\sigma )=\int_{0}^{\infty }t^{\sigma -1}\ln \left(1+%
\frac{\rho }{t^{\eta }}\right)dt \\
&=&\frac{1}{\sigma }\int_{0}^{\infty }\ln \left(1+\frac{\rho }{t^{\eta }}%
\right)dt^{\sigma }
=\frac{1}{\sigma }\left[ t^{\sigma }\ln \left(1+\frac{\rho }{t^{\eta }}%
\right)|_{0}^{\infty }-\int_{0}^{\infty }t^{\sigma }d\ln \left(1+\frac{\rho }{t^{\eta }}%
\right)\right]\\
&=&\frac{\eta }{\sigma }\int_{0}^{\infty }\frac{t^{\sigma -1}}{(t^{\eta
}/\rho )+1}dt=\frac{\rho ^{\sigma /\eta }}{\sigma }\int_{0}^{\infty }\frac{%
u^{(\sigma /\eta )-1}}{u+1}du\\
&=&\frac{\rho ^{\sigma /\eta }\pi }{\sigma \sin \pi (\sigma /\eta )}\in
\mathbf{R}_{+}.
\end{eqnarray*}%
In view of (\ref{49}) and (\ref{57}), we have%
$$
||T_{1}||=||T_{2}||=\frac{\rho ^{\sigma /\eta }\pi }{\sigma \sin \pi (\sigma
/\eta )}.
$$
\end{example}

(b) Set $$h(t)=k_{0}(1,t)=\arctan \left(\frac{\rho }{t^{\eta }}\right)\ \ (\rho
>0,0<\sigma <\eta ).$$ Then we have the kernels $$h(xy)=\arctan \left(\frac{\rho }{%
(xy)^{\eta }}\right),\ \ k_{0}(x,y)=\arctan \rho \left(\frac{x}{y}\right)^{\eta } $$ and obtain
the constant factors%
\begin{eqnarray*}
k(\sigma ) &=&k_{0 }(\sigma )=\int_{0}^{\infty }t^{\sigma -1}\arctan \left(\frac{%
\rho }{t^{\eta }}\right)dt \\
&=&\frac{1}{\sigma }\int_{0}^{\infty }\arctan \left(\frac{\rho }{t^{\eta }}%
\right)dt^{\sigma } =\frac{1}{\sigma }\left[ t^{\sigma }\arctan \left(\frac{\rho }{%
t^{\eta }}\right)|_{0}^{\infty }-\int_{0}^{\infty }t^{\sigma }d\arctan \left(\frac{\rho
}{t^{\eta }}\right)\right] \\
&=&\frac{\eta }{\sigma \rho }\int_{0}^{\infty }\frac{t^{\eta +\sigma -1}}{%
(t^{2\eta }/\rho ^{2})+1}dt=\frac{\rho ^{\sigma /\eta }}{2\sigma }%
\int_{0}^{\infty }\frac{u^{[(\eta +\sigma )/(2\eta )]-1}}{u+1}du \\
&=&\frac{\rho ^{\sigma /\eta }\pi }{2\sigma \sin \pi \lbrack (\eta +\sigma
)/(2\eta )]}=\frac{\rho ^{\sigma /\eta }\pi }{2\sigma \cos \pi \lbrack
\sigma /(2\eta )]}\in \mathbf{R}_{+}.
\end{eqnarray*}%
By (\ref{49}) and (\ref{57}), we have $$||T_{1}||=||T_{2}||=\frac{%
\rho ^{\sigma /\eta }\pi }{2\sigma \cos \pi \lbrack \sigma /(2\eta )]}.$$

(c) Set $$h(t)=k_{0}(1,t)=e^{-\rho t^{\eta }}\ \ (\rho ,\sigma ,\eta >0).$$
Then we have the kernels $$h(xy)=e^{-\rho (xy)^{\eta }},\ \ k_{0}(x,y)=e^{-\rho (%
\frac{y}{x})^{\eta }} $$ and obtain the constant factors
\begin{eqnarray*}
k(\sigma ) &=&k_{0}(\sigma )=\int_{0}^{\infty }t^{\sigma -1}e^{-\rho t^{\eta
}}dt \\
&=&\frac{1}{\eta \rho ^{\sigma /\eta }}\int_{0}^{\infty }e^{-u}u^{\frac{%
\sigma }{\eta }-1}du=\frac{1}{\eta \rho ^{\sigma /\eta }}\Gamma \left(\frac{%
\sigma }{\eta }\right)\in \mathbf{R}_{+}.
\end{eqnarray*}%
In view of (\ref{49}) and (\ref{57}), we have $$||T_{1}||=||T_{2}||=\frac{1}{%
\eta \rho ^{\sigma /\eta }}\Gamma \left(\frac{\sigma }{\eta }\right). $$

Then for (a)-(c), we can obtain the equivalent inequalities with the kernels
and the best possible constant factors in Theorems 1-4. Setting $\delta _{\:0}=%
\frac{\sigma }{2}>0,$ we can still obtain the equivalent reverse
inequalities with the kernels and the best possible constant factors in
Theorems 1-4.

\begin{example}
(a) Set $$h(t)=k_{0}(1,t)=\csc h(\rho t^{\eta })=\frac{%
2}{e^{\rho t^{\eta }}-e^{-\rho t^{\eta }}}\ \ (\rho >0,0<\eta <\sigma ).$$ Where
$\csc h(\cdot )$ stands for the hyperbolic cosecant function (cf. \cite{YM26}).
Then we have the kernels
 $$h(xy)=\frac{2}{e^{\rho (xy)^{\eta }}-e^{-\rho (xy)^{\eta }}},\ \ k_{0}(x,y)=%
\frac{2}{e^{\rho (\frac{y}{x})^{\eta }}-e^{-\rho (\frac{y}{x})^{\eta }}}.
$$
By the Lebesgue term by term integration theorem, we obtain the constant factors%
\begin{eqnarray*}
k(\sigma ) &=&k_{0}(\sigma )=\int_{0}^{\infty }\frac{2t^{\sigma -1}dt%
}{e^{\rho t^{\eta }}-e^{-\rho t^{\eta }}}=\int_{0}^{\infty }\frac{2t^{\sigma
-1}dt}{e^{\rho t^{\eta }}(1-e^{-2\rho t^{\eta }})} \\
&=&2\int_{0}^{\infty }t^{\sigma -1}\sum_{k=0}^{\infty }e^{-(2k+1)\rho
t^{\eta }}dt=2\sum_{k=0}^{\infty }\int_{0}^{\infty }t^{\sigma
-1}e^{-(2k+1)\rho t^{\eta }}dt\\
&=&\frac{2}{\eta \rho ^{\sigma /\eta }}\sum_{k=0}^{\infty }\frac{1}{%
(2k+1)^{\sigma /\eta }}\int_{0}^{\infty }e^{-u}u^{\frac{\sigma }{\eta }-1}du\\
&=&\frac{2}{\eta \rho ^{\sigma /\eta }}\Gamma \left(\frac{\sigma }{\eta }%
\right)\sum_{k=0}^{\infty }\frac{1}{(2k+1)^{\sigma /\eta }}\\
&=&\frac{2}{\eta \rho ^{\sigma /\eta }}\Gamma \left(\frac{\sigma }{\eta }\right)\left[
\sum_{k=1}^{\infty }\frac{1}{k^{\sigma /\eta }}-\sum_{k=1}^{\infty }\frac{1}{%
(2k)^{\sigma /\eta }}\right] \\
&=&\frac{2}{\eta \rho ^{\sigma /\eta }}\left(1-\frac{1}{2^{\sigma /\eta }}\right)\Gamma
\left(\frac{\sigma }{\eta }\right)\zeta \left(\frac{\sigma }{\eta }\right)\in \mathbf{R}_{+},
\end{eqnarray*}%
where, $\zeta \left(\frac{\sigma }{\eta }\right)=\sum_{k=1}^{\infty }\frac{1}{k^{\sigma
/\eta }}$ ($\zeta \left(\cdot \right)$ is Riemann zeta function). In view of (\ref{49})
and (\ref{57}), we have%
$$
||T_{1}||=||T_{2}||=\frac{2}{\eta \rho ^{\sigma /\eta }}\left(1-\frac{1}{%
2^{\sigma /\eta }}\right)\Gamma \left(\frac{\sigma }{\eta }\right)\zeta \left(\frac{\sigma }{\eta }%
\right).
$$
\end{example}

(b) Set%
\begin{eqnarray*}
h(t) &=&k_{0}(1,t)=e^{-\rho t^{\eta }}\cot h(\rho t^{\eta })=e^{-\rho
t^{\eta }}\frac{e^{\rho t^{\eta }}+e^{-\rho t^{\eta }}}{e^{\rho t^{\eta
}}-e^{-\rho t^{\eta }}} \\
&=&\frac{1+e^{-2\rho t^{\eta }}}{e^{\rho t^{\eta }}-e^{-\rho t^{\eta }}}=%
\frac{e^{-\rho t^{\eta }}+e^{-3\rho t^{\eta }}}{1-e^{-2\rho t^{\eta }}}\ \ (\rho
>0,0<\eta <\sigma ).
\end{eqnarray*}
Let $\cot h(\cdot )$ stand for the hyperbolic cotangent function (cf. \cite{YM26}).
Then we have the kernels $$h(xy)=\frac{1+e^{-2\rho (xy)^{\eta }}}{e^{\rho
(xy)^{\eta }}-e^{-\rho (xy)^{\eta }}},\ \ k_{0}(x,y)=\frac{1+e^{-2\rho (\frac{x}{%
y})^{\eta }}}{e^{\rho (\frac{y}{x})^{\eta }}-e^{-\rho (\frac{y}{x})^{\eta }}}%
. $$ By the Lebesgue term by term integration theorem, we obtain the constant
factors%
\begin{eqnarray*}
k(\sigma ) &=&k_{0 }(\sigma )=\int_{0}^{\infty }\frac{(e^{-\rho t^{\eta
}}+e^{-3\rho t^{\eta }})t^{\sigma -1}}{1-e^{-2\rho t^{\eta }}}dt \\
&=&\int_{0}^{\infty }t^{\sigma -1}\sum_{k=0}^{\infty }(e^{-(2k+1)\rho
t^{\eta }}+e^{-(2k+3)\rho t^{\eta }})dt \\
&=&\frac{1}{\eta \rho ^{\sigma /\eta }}\sum_{k=0}^{\infty }\left[\frac{1}{%
(2k+1)^{\sigma /\eta }}+\frac{1}{(2k+3)^{\sigma /\eta }}\right]\int_{0}^{\infty
}e^{-u}u^{\frac{\sigma }{\eta }-1}du \\
&=&\frac{1}{\eta \rho ^{\sigma /\eta }}\Gamma \left(\frac{\sigma }{\eta }\right)\left[
2\sum_{k=0}^{\infty }\frac{1}{(2k+1)^{\sigma /\eta }}-1\right] \\
&=&\frac{1}{\eta \rho ^{\sigma /\eta }}\Gamma \left(\frac{\sigma }{\eta }\right)\left[(2-%
\frac{1}{2^{(\sigma /\eta )-1}})\zeta \left(\frac{\sigma }{\eta }\right)-1\right]\in \mathbf{R%
}_{+}.
\end{eqnarray*}%
In view of (\ref{49}) and (\ref{57}), we have
\[
||T_{1}||=||T_{2}||=\frac{1}{\eta \rho ^{\sigma /\eta }}\Gamma \left(\frac{\sigma
}{\eta }\right)\left[\left(2-\frac{1}{2^{(\sigma /\eta )-1}}\right)\zeta \left(\frac{\sigma }{\eta }%
\right)-1\right].
\]
Then for (a)-(b), we can obtain the equivalent inequalities with the kernels
and the best possible constant factors in Theorems 1-4. Setting $\delta _{\:0}=%
\frac{\sigma -\eta }{2}>0,$ we can still obtain the equivalent reverse
inequalities with the kernels and the best possible constant factors in
Theorems 1-4.

(c) Set $$h(t)=k_{0}(1,t)=\sec h(\rho t^{\eta })=\frac{2}{e^{\rho t^{\eta
}}+e^{-\rho t^{\eta }}}\ \ (\rho ,\eta ,\sigma >0).$$ Let $\sec h(\cdot )$ stand for the
hyperbolic secant function (cf. \cite{YM26}). Then we have the kernels
\[
h(xy)=\frac{2}{e^{\rho (xy)^{\eta }}+e^{-\rho (xy)^{\eta }}},\ \ k_{0}(x,y)=%
\frac{2}{e^{\rho (\frac{y}{x})^{\eta }}+e^{-\rho (\frac{y}{x})^{\eta }}}.
\]
By the Lebesgue term by term integration theorem, we obtain the constant factors%
\begin{eqnarray*}
k(\sigma ) &=&k_{0}(\sigma )=\int_{0}^{\infty }\frac{2t^{\sigma -1}dt}{%
e^{\rho t^{\eta }}+e^{-\rho t^{\eta }}}=\int_{0}^{\infty }\frac{2t^{\sigma
-1}dt}{e^{\rho t^{\eta }}(1+e^{-2\rho t^{\eta }})} \\
&=&2\int_{0}^{\infty }t^{\sigma -1}\sum_{k=0}^{\infty
}(-1)^{k}e^{-(2k+1)\rho t^{\eta }}dt \\
&=&2\int_{0}^{\infty }t^{\sigma -1}\sum_{k=0}^{\infty }[e^{-(4k+1)\rho
t^{\eta }}-e^{-(4k+3)\rho t^{\eta }}]dt \\
&=&2\sum_{k=0}^{\infty }\int_{0}^{\infty }t^{\sigma -1}[e^{-(4k+1)\rho
t^{\eta }}-e^{-(4k+3)\rho t^{\eta }}]dt \\
&=&2\sum_{k=0}^{\infty }(-1)^{k}\int_{0}^{\infty }t^{\sigma
-1}e^{-(2k+1)\rho t^{\eta }}dt \\
&=&\frac{2}{\eta \rho ^{\sigma /\eta }}\sum_{k=0}^{\infty }\frac{(-1)^{k}}{%
(2k+1)^{\sigma /\eta }}\int_{0}^{\infty }e^{-u}u^{\frac{\sigma }{\eta }-1}du
\\
&=&\frac{2}{\eta \rho ^{\sigma /\eta }}\Gamma \left(\frac{\sigma }{\eta }\right)\xi \left(%
\frac{\sigma }{\eta }\right)\in \mathbf{R}_{+},
\end{eqnarray*}%
where $$\xi (\frac{\sigma }{\eta })=\sum_{k=0}^{\infty }\frac{(-1)^{k}}{%
(2k+1)^{\sigma /\eta }}.$$ In view of (\ref{49}) and (\ref{57}), we have%
\[
||T_{1}||=||T_{2}||=\frac{2}{\eta \rho ^{\sigma /\eta }}\Gamma (\frac{\sigma
}{\eta })\xi (\frac{\sigma }{\eta }).
\]

(d) Set%
\begin{eqnarray*}
h(t) &=&k_{0}(1,t)=e^{-\rho t^{\eta }}\tan h(\rho t^{\eta })=e^{-\rho
t^{\eta }}\frac{e^{\rho t^{\eta }}-e^{-\rho t^{\eta }}}{e^{\rho t^{\eta
}}+e^{-\rho t^{\eta }}} \\
&=&\frac{1-e^{-2\rho t^{\eta }}}{e^{\rho t^{\eta }}+e^{-\rho t^{\eta }}}=%
\frac{e^{-\rho t^{\eta }}-e^{-3\rho t^{\eta }}}{1+e^{-2\rho t^{\eta }}}(\rho
,\eta ,\sigma >0).
\end{eqnarray*}%
Let $\tan h(\cdot )$ stand for the hyperbolic tangent function (cf. \cite{YM26}).
Then we have the kernels
\[
h(xy)=\frac{1-e^{-2\rho (xy)^{\eta }}}{e^{\rho (xy)^{\eta }}+e^{-\rho
(xy)^{\eta }}},\ \ k_{0}(x,y)=\frac{1-e^{-2\rho (\frac{y}{x})^{\eta }}}{e^{\rho (%
\frac{y}{x})^{\eta }}+e^{-\rho (\frac{y}{x})^{\eta }}}.
\]

By the Lebesgue term by term integration theorem, we obtain the constant factors%
\begin{eqnarray*}
k(\sigma )&=&k_{0}(\sigma )=\int_{0}^{\infty }\frac{(e^{-\rho t^{\eta
}}-e^{-3\rho t^{\eta }})t^{\sigma -1}}{1+e^{-2\rho t^{\eta }}}dt \\
&=&\int_{0}^{\infty }t^{\sigma -1}\sum_{k=0}^{\infty }(-1)^{k}e^{-(2k+1)\rho
t^{\eta }}dt-\int_{0}^{\infty }t^{\sigma -1}\sum_{k=0}^{\infty
}(-1)^{k}e^{-(2k+3)\rho t^{\eta }}dt \\
&=&\int_{0}^{\infty }t^{\sigma -1}\sum_{k=0}^{\infty }[e^{-(4k+1)\rho
t^{\eta }}-e^{-(4k+3)\rho t^{\eta }}]dt \\
&&-\int_{0}^{\infty }t^{\sigma -1}\sum_{k=0}^{\infty }[e^{-(4k+3)\rho
t^{\eta }}-e^{-(4k+5)\rho t^{\eta }}]dt \\
&=&\sum_{k=0}^{\infty }\left\{ \int_{0}^{\infty }t^{\sigma
-1}[e^{-(4k+1)\rho t^{\eta }}-e^{-(4k+3)\rho t^{\eta }}]dt\right. \\
&&\left. -\int_{0}^{\infty }t^{\sigma -1}[e^{-(4k+3)\rho t^{\eta
}}-e^{-(4k+5)\rho t^{\eta }}]dt\}\right\} \\
&=&\frac{1}{\eta \rho ^{\sigma /\eta }}\sum_{k=0}^{\infty }(-1)^{k}[\frac{1}{%
(2k+1)^{\sigma /\eta }}-\frac{1}{(2k+3)^{\sigma /\eta }}]\int_{0}^{\infty
}e^{-u}u^{\frac{\sigma }{\eta }-1}du \\
&=&\frac{1}{\eta \rho ^{\sigma /\eta }}\Gamma \left(\frac{\sigma }{\eta }\right)\left[
2\sum_{k=0}^{\infty }\frac{(-1)^{k}}{(2k+1)^{\sigma /\eta }}-1\right] \\
&=&\frac{1}{\eta \rho ^{\sigma /\eta }}\Gamma \left(\frac{\sigma }{\eta }\right)\left(2\xi \left(%
\frac{\sigma }{\eta }\right)-1\right)\in \mathbf{R}_{+}.
\end{eqnarray*}

In view of (\ref{49}) and (\ref{57}), we have
\[
||T_{1}||=||T_{2}||=\frac{1}{\eta \rho ^{\sigma /\eta }}\Gamma \left(\frac{\sigma
}{\eta }\right)\left(2\xi \left(\frac{\sigma }{\eta }\right)-1\right).
\]

Then for (c)-(d), we can obtain the equivalent inequalities with the kernels
and the best possible constant factors in Theorems 1-4. Setting $\delta _{\:0}=%
\frac{\sigma }{2}>0,$ we can still obtain the equivalent reverse
inequalities with the kernels and the best possible constant factors in
Theorems 1-4.
\begin{lemma}
 Let $\mathbf{C}$ stand for the set of complex numbers. If $\mathbf{%
C}_{\infty }=\mathbf{C}\cup \{\infty \},$ $$z_{k}\in \mathbf{C}\backslash \{z|%
\func{Re}z\geq 0,\ \func{Im}z=0\}\ \ (k=1,\:2,\ldots,\:n)$$ are different
points, the function $f(z)$ is analytic in $\mathbf{C}_{\infty }$ except for
$z_{i},\ i=1,\:2,\ldots,\:n$, and $z=\infty $ is a zero point of $f(z)$ whose
order is not less than 1, then for $\alpha \in \mathbf{R},$ we have%
\begin{equation}
\int_{0}^{\infty }f(x)x^{\alpha -1}dx=\frac{2\pi i}{1-e^{2\pi \alpha i}}%
\sum_{k=1}^{n}\func{Re}s[f(z)z^{\alpha -1},z_{k}],  \label{62}
\end{equation}%
where, $0<\func{Im}\ln z=\arg z<2\pi $. In particular, if $z_{k},\ k=1,\ldots
,\:n,$ are all poles of order 1, setting $$\varphi
_{k}(z)=(z-z_{k})f(z)\ \ (\varphi _{k}(z_{k})\neq 0),$$ then
\begin{equation}
\int_{0}^{\infty }f(x)x^{\alpha -1}dx=\frac{\pi }{\sin \pi \alpha }%
\sum_{k=1}^{n}(-z_{k})^{\alpha -1}\varphi _{k}(z_{k}).  \label{63}
\end{equation}
\end{lemma}

\begin{proof}
In view of the theorem (cf. \cite{YM27}, p. 118), we obtain (%
\ref{62}). We have%
\begin{eqnarray*}
1-e^{2\pi \alpha i} &=&1-\cos 2\pi \alpha -i\sin 2\pi \alpha \\
&=&-2i\sin \pi \alpha (\cos \pi \alpha +i\sin \pi \alpha )=-2ie^{i\pi \alpha
}\sin \pi \alpha .
\end{eqnarray*}%
In particular, since $$f(z)z^{\alpha -1}=\frac{1}{z-z_{k}}\varphi
_{k}(z)z^{\alpha -1},$$ it is obvious that%
$$
\func{Re}s[f(z)z^{\alpha -1},-a_{k}]=z_{k}{}^{\alpha -1}\varphi
_{k}(z_{k})=-e^{i\pi \alpha }(-z_{k})^{\alpha -1}\varphi _{k}(z_{k}).
$$
Then by (\ref{62}), we obtain (\ref{63}). $\Box $
\end{proof}

\begin{example}
 (a) Set $$h(t)=k_{\lambda }(1,t)=\frac{1}{%
\prod_{k=1}^{s}(a_{k}+t^{\lambda /s})}$$
where $s\in \mathbf{N},\:0<a_{1}<\cdots
<a_{s},\:\mu ,\:\sigma >0,\:\mu +\sigma =\lambda .$
Then we have the kernels%
$$
h(xy)=\frac{1}{\prod_{k=1}^{s}[a_{k}+(xy)^{\lambda /s}]},\ \ k_{\lambda }(x,y)=%
\frac{1}{\prod_{k=1}^{s}(a_{k}x^{\lambda /s}+y^{\lambda /s})}.
$$
For $$f(z)=\frac{1}{\prod_{k=1}^{s}(z+a_{k})},\ \ z_{k}=-a_{k},$$ by (\ref{63}),
we get%
$$
\varphi _{k}(z_{k})=(z+a_{k})\frac{1}{\prod_{i=1}^{s}(z+a_{i})}%
\Big|_{z=-a_{k}}=\prod_{j=1(j\neq k)}^{s}\frac{1}{a_{j}-a_{k}},
$$
and obtain the constant factors
\begin{eqnarray*}
k(\sigma ) &=&k_{\lambda }(\sigma )=\int_{0}^{\infty }\frac{t^{\sigma -1}dt}{%
\prod_{k=1}^{s}(a_{k}+t^{\lambda /s})}
=\frac{s}{\lambda }\int_{0}^{\infty }\frac{u^{(s\sigma /\lambda )-1}du}{%
\prod_{k=1}^{s}(u+a_{k})} \\
&=&\frac{\pi s}{\lambda \sin \pi (s\sigma /\lambda )}\sum_{k=1}^{s}a_{k}^{s%
\sigma /\lambda }\prod_{j=1(j\neq k)}^{s}\frac{1}{a_{j}-a_{k}}\in \mathbf{R}%
_{+}.
\end{eqnarray*}%
In view of (\ref{49}) and (\ref{57}), we have
$$
||T_{1}||=||T_{2}||=\frac{\pi s}{\lambda \sin \pi (s\sigma /\lambda )}%
\sum_{k=1}^{s}a_{k}^{s\sigma /\lambda }\prod_{j=1(j\neq k)}^{s}\frac{1}{%
a_{j}-a_{k}}.
$$
\end{example}

In particular,
(i) if $s=1,a_{1}=a,h(t)=k_{\lambda }(1,t)=\frac{1}{%
a+t^{\lambda }}(a,\mu,\sigma>0$ $\mu+\sigma =\lambda ),$ then we have the
kernels $h(xy)=\frac{1}{a+(xy)^{\lambda }},k_{\lambda }(x,y)=\frac{1}{%
ax^{\lambda }+y^{\lambda }}, $ and
\[
||T_{1}||=||T_{2}||=\frac{\pi }{\lambda \sin \pi (\sigma /\lambda )}a^{\frac{%
\sigma }{\lambda }-1};
\]
(ii) if $s=2,\:a_{1}=a,\:a_{2}=b,$ $$h(t)=k_{\lambda }(1,t)=\frac{1}{(a+t^{\lambda
/2})(b+t^{\lambda /2})}\ \ (0<a<b,\:\mu ,\:\sigma >0,\:\mu +\sigma =\lambda ),$$
then we have the kernels
\[
h(xy) =\frac{1}{[a+(xy)^{\lambda /2}][b+(xy)^{\lambda /2}]},\ \  k_{\lambda
}(x,y) =\frac{1}{(ax^{\lambda /2}+y^{\lambda /2})(ax^{\lambda /2}+y^{\lambda
/2})},
\]
and
\[
||T_{1}||=||T_{2}||=\frac{2\pi }{\lambda \sin \pi (2\sigma /\lambda )}\frac{1%
}{b-a}(a^{\frac{2\sigma }{\lambda }-1}-b^{\frac{2\sigma }{\lambda }-1}).
\]

(b) Set $$h(t)=k_{\lambda }(1,t)=\frac{1}{t^{\lambda }+2ct^{\lambda
/2}\cos \gamma +c^{2}}$$ $(c>0,\:|\gamma |<\frac{\pi }{2},\:\mu,\:\sigma >0,\:\mu
,\:\sigma =\lambda ).$
 Then we have the kernels
$$
h(xy) =\frac{1}{(xy)^{\lambda }+2c(xy)^{\lambda /2}\cos \gamma +c^{2}},$$ 
$$
k_{\lambda }(x,y) =\frac{1}{y^{\lambda }+2c(xy)^{\lambda /2}\cos \gamma
+c^{2}x^{\lambda }}.$$
By (\ref{63}), we can find%
\begin{eqnarray*}
k(\sigma ) &=&k_{\lambda }(\sigma )=\int_{0}^{\infty }\frac{t^{\sigma -1}}{%
t^{\lambda }+2ct^{\lambda /2}\cos \gamma +c^{2}}dt \\
&=&\frac{2}{\lambda }\int_{0}^{\infty }\frac{u^{(2\sigma /\lambda )-1}du}{%
u^{2}+2cu\cos \gamma +c^{2}}=\frac{2}{\lambda }\int_{0}^{\infty }\frac{%
u^{(2\sigma /\lambda )-1}du}{(u+ce^{i\gamma })(u+ce^{-i\gamma })} \\
&=&\frac{2\pi }{\lambda \sin \pi (2\sigma /\lambda )}[\frac{(ce^{i\gamma
})^{(2\sigma /\lambda )-1}}{c(e^{-i\gamma }-e^{i\gamma })}+\frac{%
(ce^{-i\gamma })^{(2\sigma /\lambda )-1}}{c(e^{i\gamma }-e^{-i\gamma })}] \\
&=&\frac{2\pi \sin \gamma (1-2\sigma /\lambda )}{\lambda \sin \pi (2\sigma
/\lambda )\sin \gamma }c^{\frac{2\sigma }{\lambda }-2}\in \mathbf{R}_{+}.
\end{eqnarray*}%
In view of (\ref{49}) and (\ref{57}), we have
\[
||T_{1}||=||T_{2}||=\frac{2\pi \sin \gamma (1-2\sigma /\lambda )}{\lambda
\sin \pi (2\sigma /\lambda )\sin \gamma }c^{\frac{2\sigma }{\lambda }-2}.
\]

Then for (a)-(b), we can obtain the equivalent inequalities with the kernels
and the best possible constant factors in Theorems 1-4. Setting $\delta _{\:0}=%
\frac{\sigma }{2}>0,$ we can obtain the equivalent reverse
inequalities with the kernels and the best possible constant factors in
Theorems1-4.
\begin{remark}
Setting $p=q=2,\:\mu =\sigma =\frac{\lambda }{2}$ in Theorem 1
and Corollary 3, in view of Remark 3 and the above results, if $%
f(x),\:g(y)\geq 0,$ with $$0<\int_{0}^{\infty }x^{1-\lambda
}f^{2}(x)dx<\infty$$ and $$ 0<\int_{0}^{\infty }y^{1-\lambda
}g^{2}(y)dy<\infty,$$
 then we obtain the following 8 couples of simpler
equivalent inequalities with an independent parameter $\lambda $ and the
best possible constant factors:
\end{remark}

(a) For $\lambda >0,$%
\begin{equation}
\int_{0}^{\infty }\int_{0}^{\infty }\frac{f(x)g(y)}{x^{\lambda }+y^{\lambda }%
}dxdy <\frac{\pi }{\lambda }\left( \int_{0}^{\infty }x^{1-\lambda
}f^{2}(x)dx\int_{0}^{\infty }y^{1-\lambda }g^{2}(y)dy\right) ^{\frac{1}{2}},
\label{64}
\end{equation}%
\begin{equation}
\int_{0}^{\infty }\int_{0}^{\infty }\frac{f(x)g(y)}{1+(xy)^{\lambda }}dxdy <%
\frac{\pi }{\lambda }\left( \int_{0}^{\infty }x^{1-\lambda
}f^{2}(x)dx\int_{0}^{\infty }y^{1-\lambda }g^{2}(y)dy\right) ^{\frac{1}{2}};
\label{65}
\end{equation}%
\begin{equation}
\int_{0}^{\infty }\int_{0}^{\infty }\frac{f(x)g(y)}{(x+y)^{\lambda }}dxdy <B\left(%
\frac{\lambda }{2},\frac{\lambda }{2}\right)\left( \int_{0}^{\infty }x^{1-\lambda
}f^{2}(x)dx\int_{0}^{\infty }y^{1-\lambda }g^{2}(y)dy\right) ^{\frac{1}{2}},
\label{66}
\end{equation}%
\begin{equation}
\int_{0}^{\infty }\int_{0}^{\infty }\frac{f(x)g(y)}{(1+xy)^{\lambda }}dxdy
<B\left(\frac{\lambda }{2},\frac{\lambda }{2}\right)\left( \int_{0}^{\infty
}x^{1-\lambda }f^{2}(x)dx\int_{0}^{\infty }y^{1-\lambda }g^{2}(y)dy\right) ^{%
\frac{1}{2}},  \label{67}
\end{equation}%
\begin{equation}
\int_{0}^{\infty }\int_{0}^{\infty }\frac{f(x)g(y)}{(\max \{x,y\})^{\lambda }%
}dxdy <\frac{4}{\lambda }\left( \int_{0}^{\infty }x^{1-\lambda
}f^{2}(x)dx\int_{0}^{\infty }y^{1-\lambda }g^{2}(y)dy\right) ^{\frac{1}{2}},
\label{68}
\end{equation}%
\begin{equation}
\int_{0}^{\infty }\int_{0}^{\infty }\frac{f(x)g(y)}{(\max \{xy,1\})^{\lambda
}}dxdy <\frac{4}{\lambda }\left( \int_{0}^{\infty }x^{1-\lambda
}f^{2}(x)dx\int_{0}^{\infty }y^{1-\lambda }g^{2}(y)dy\right) ^{\frac{1}{2}};
\label{69}
\end{equation}%
\begin{equation}
\int_{0}^{\infty }\int_{0}^{\infty }\frac{|\ln (\frac{x}{y})|f(x)g(y)}{(\max
\{x,y\})^{\lambda }}dxdy <\frac{8}{\lambda ^{2}}\left( \int_{0}^{\infty
}x^{1-\lambda }f^{2}(x)dx\int_{0}^{\infty }y^{1-\lambda }g^{2}(y)dy\right) ^{%
\frac{1}{2}},  \label{70}
\end{equation}%
\begin{equation}
\int_{0}^{\infty }\int_{0}^{\infty }\frac{|\ln (xy)|f(x)g(y)}{(\max
\{1,xy\})^{\lambda }}dxdy <\frac{8}{\lambda ^{2}}\left( \int_{0}^{\infty
}x^{1-\lambda }f^{2}(x)dx\int_{0}^{\infty }y^{1-\lambda }g^{2}(y)dy\right) ^{%
\frac{1}{2}};  \label{71}
\end{equation}%
\begin{equation}
\int_{0}^{\infty }\int_{0}^{\infty }\frac{\ln (\frac{x}{y})f(x)g(y)}{%
x^{\lambda }-y^{\lambda }}dxdy <\left(\frac{\pi }{\lambda }\right)^{2}\left(
\int_{0}^{\infty }x^{1-\lambda }f^{2}(x)dx\int_{0}^{\infty }y^{1-\lambda
}g^{2}(y)dy\right) ^{\frac{1}{2}},  \label{72}
\end{equation}%
\begin{equation}
\int_{0}^{\infty }\int_{0}^{\infty }\frac{\ln (xy)f(x)g(y)}{(xy)^{\lambda }-1%
}dxdy <\left(\frac{\pi }{\lambda }\right)^{2}\left( \int_{0}^{\infty }x^{1-\lambda
}f^{2}(x)dx\int_{0}^{\infty }y^{1-\lambda }g^{2}(y)dy\right) ^{\frac{1}{2}};
\label{73}
\end{equation}

(b) for $0<\lambda <1,$%
\begin{equation}
\int_{0}^{\infty }\int_{0}^{\infty }\frac{f(x)g(y)}{|x-y|^{\lambda }}dxdy
<2B\left(1-\lambda ,\frac{\lambda }{2}\right)\left( \int_{0}^{\infty }x^{1-\lambda
}f^{2}(x)dx\int_{0}^{\infty }y^{1-\lambda }g^{2}(y)dy\right) ^{\frac{1}{2}},
\label{74}
\end{equation}%
\begin{equation}
\int_{0}^{\infty }\int_{0}^{\infty }\frac{f(x)g(y)}{|1-xy|^{\lambda }}dxdy
<2B\left(1-\lambda ,\frac{\lambda }{2}\right)\left( \int_{0}^{\infty }x^{1-\lambda
}f^{2}(x)dx\int_{0}^{\infty }y^{1-\lambda }g^{2}(y)dy\right) ^{\frac{1}{2}};
\label{75}
\end{equation}

(c) for $\lambda <0$,%
\begin{equation}
\int_{0}^{\infty }\int_{0}^{\infty }\frac{f(x)g(y)}{(\min \{x,y\})^{\lambda }%
}dxdy <\frac{-4}{\lambda }\left( \int_{0}^{\infty }x^{1-\lambda
}f^{2}(x)dx\int_{0}^{\infty }y^{1-\lambda }g^{2}(y)dy\right) ^{\frac{1}{2}},
\label{76}
\end{equation}%
\begin{equation}
\int_{0}^{\infty }\int_{0}^{\infty }\frac{f(x)g(y)}{(\min \{1,xy\})^{\lambda
}}dxdy <\frac{-4}{\lambda }\left( \int_{0}^{\infty }x^{1-\lambda
}f^{2}(x)dx\int_{0}^{\infty }y^{1-\lambda }g^{2}(y)dy\right) ^{\frac{1}{2}};
\label{77}
\end{equation}%
\begin{equation}
\int_{0}^{\infty }\int_{0}^{\infty }\frac{|\ln (\frac{x}{y})|f(x)g(y)}{(\min
\{x,y\})^{\lambda }}dxdy <\frac{8}{\lambda ^{2}}\left( \int_{0}^{\infty
}x^{1-\lambda }f^{2}(x)dx\int_{0}^{\infty }y^{1-\lambda }g^{2}(y)dy\right) ^{%
\frac{1}{2}},  \label{78}
\end{equation}%
\begin{equation}
\int_{0}^{\infty }\int_{0}^{\infty }\frac{|\ln (xy)|f(x)g(y)}{(\min
\{1,xy\})^{\lambda }}dxdy <\frac{8}{\lambda ^{2}}\left( \int_{0}^{\infty
}x^{1-\lambda }f^{2}(x)dx\int_{0}^{\infty }y^{1-\lambda }g^{2}(y)dy\right) ^{%
\frac{1}{2}}.  \label{79}
\end{equation}

\section{Yang-Hilbert-Type Integral Inequalities in the Whole Plane}

\label{sec:3} \bigskip In this section, we study some Yang-Hilbert-type
integral inequalities in the whole plane with parameters and the best
constant factors. The equivalent forms, the reverses, the Hardy-type
inequalities, the operator expressions and some particular examples are also
discussed.

\subsection{Weight Functions and a Lemma}

\begin{definition}
 If $\delta \in \{-1,1\},\sigma \in \mathbf{R}%
,\:H(t)$ is a non-negative measurable function in $\mathbf{R}$, define the
following weight functions:%
\begin{eqnarray}
\omega _{\delta }(\sigma ,y) &:&=|y|^{\sigma }\int_{-\infty }^{\infty
}H(x^{\delta }y)|x|^{\delta \sigma -1}dx(y\in \mathbf{R}\backslash \{0\}),
\label{80} \\
\varpi _{\delta }(\sigma ,x) &:&=|x|^{\delta \sigma }\int_{-\infty }^{\infty
}H(x^{\delta }y)|y|^{\sigma -1}dy(x\in \mathbf{R}\backslash \{0\}).
\label{81}
\end{eqnarray}
\end{definition}

Setting $t=x^{\delta }y$ in (\ref{80}), we obtain $x=y^{-\frac{1}{\delta }}t^{%
\frac{1}{\delta }},\:dx=\frac{1}{\delta }y^{-\frac{1}{\delta }}t^{\frac{1}{%
\delta }-1}dt$ and
\begin{eqnarray}
\omega _{\delta }(\sigma ,y) &=&|y|^{\sigma }\int_{-\infty }^{\infty
}H(t)|y^{-\frac{1}{\delta }}t^{\frac{1}{\delta }}|^{\delta \sigma -1}|y|^{-%
\frac{1}{\delta }}|t|^{\frac{1}{\delta }-1}dt  \nonumber \\
&=&K(\sigma ):=\int_{-\infty }^{\infty }H(t)|t|^{\sigma -1}dt.  \label{82}
\end{eqnarray}%
Setting $t=x^{\delta }y$ in (\ref{81}), we find $y=x^{-\delta
}t,\: dy=x^{-\delta }dt$ and%
\begin{equation}
\varpi _{\delta }(\sigma ,x)=|x|^{\delta \sigma }\int_{-\infty }^{\infty
}H(t)|x^{-\delta }t|^{\sigma -1}|x|^{-\delta }dt=K(\sigma ).  \label{83}
\end{equation}

\begin{remark}
 We can still get%
\begin{eqnarray}
K(\sigma ) &=&\int_{-\infty }^{0}H(t)(-t)^{\sigma -1}dt+\int_{0}^{\infty
}H(t)t^{\sigma -1}dt   \nonumber \\
&=&\int_{0}^{\infty }H(-u)u^{\sigma -1}du+\int_{0}^{\infty }H(t)t^{\sigma
-1}dt  \nonumber \\
&=&\int_{0}^{\infty }(H(-t)+H(t))t^{\sigma -1}dt.  \label{84}
\end{eqnarray}%
If $H(t)=H(-t),$ then $$K(\sigma )=2\int_{0}^{\infty }H(t)t^{\sigma -1}dt,$$
and thus we obtain again cases of integrals in the first quadrant. In this
section, we assume that $H(t)\neq H(-t).$
\end{remark}

\begin{lemma}
 If $p>0\ (p\neq 1),\:\frac{1}{p}+\frac{1}{q}=1,\sigma \in
\mathbf{R},$  both $H(t)$ and $f(t)$
are non-negative measurable functions in $\mathbf{R},$ and $K(\sigma )$ is defined by (\ref{83}),
then, (i) for $p>1,$
we have the following inequality:%
\begin{eqnarray}
J &:&=\int_{-\infty }^{\infty }|y|^{p\sigma -1}\left( \int_{-\infty
}^{\infty }H(x^{\delta }y)f(x)dx\right) ^{p}dy   \nonumber \\
&\leq &K^{p}(\sigma )\int_{-\infty }^{\infty }|x|^{p(1-\delta \sigma
)-1}f^{p}(x)dx;  \label{85}
\end{eqnarray}%
(ii) for $0<p<1,$ we have the reverse of (\ref{85}).
\end{lemma}

\begin{proof}
 (i) By H\"{o}lder's weighted inequality (cf. \cite{YM23})
and (\ref{80}), it follows that%
\begin{eqnarray}
&&\int_{-\infty }^{\infty }H(x^{\delta }y)f(x)dx\nonumber\\
&=&\int_{-\infty }^{\infty }H(x^{\delta }y)[\frac{|x|^{(1-\delta \sigma )/q}%
}{|y|^{(1-\sigma )/p}}f(x)][\frac{|y|^{(1-\sigma )/p}}{|x|^{(1-\delta \sigma
)/q}}]dx\nonumber\\
&\leq &\left[ \int_{-\infty }^{\infty }H(x^{\delta }y)\frac{|x|^{(1-\delta
\sigma )p/q}}{|y|^{1-\sigma }}f^{p}(x)dx\right] ^{\frac{1}{p}}   \nonumber \\
&&\times \left[ \int_{-\infty }^{\infty }H(x^{\delta }y)\frac{|y|^{(1-\sigma
)q/p}}{|x|^{1-\delta \sigma }}dx\right] ^{\frac{1}{q}}   \nonumber \\
&=&\frac{(\omega _{\delta }(\sigma ,y))^{\frac{1}{q}}}{|y|^{\sigma -\frac{1}{%
p}}}\left[ \int_{-\infty }^{\infty }H(x^{\delta }y)\frac{|x|^{(1-\delta
\sigma )(p-1)}}{|y|^{1-\sigma }}f^{p}(x)dx\right] ^{\frac{1}{p}}.  \label{86}
\end{eqnarray}%
Then, by (\ref{82}) and Fubini's theorem (cf. \cite{YM24}), we get%
\begin{eqnarray}
J &\leq &K^{p-1}(\sigma )\int_{-\infty }^{\infty }\int_{-\infty }^{\infty
}H(x^{\delta }y)\frac{|x|^{(1-\delta \sigma )(p-1)}}{|y|^{1-\sigma }}%
f^{p}(x)dxdy   \nonumber \\
&=&K^{p-1}(\sigma )\int_{-\infty }^{\infty }\left[ \int_{-\infty }^{\infty
}H(x^{\delta }y)\frac{|x|^{(1-\delta \sigma )(p-1)}}{|y|^{1-\sigma }}dy%
\right] f^{p}(x)dx   \nonumber \\
&=&K^{p-1}(\sigma )\int_{-\infty }^{\infty }\varpi _{\delta }(\sigma
,x)|x|^{p(1-\delta \sigma )-1}f^{p}(x)dx.  \label{87}
\end{eqnarray}%
By (\ref{83}), we obtain (\ref{85}).

(ii) For $0<p<1,$ by the reverse of H\"{o}lder's weighted inequality (cf.
\cite{YM23}), we can similarly derive the reverses of (\ref{86}) and
(\ref{87}). Then we obtain the reverse of (\ref{85}).

 This completes the proof of the lemma.
\end{proof}

\subsection{Equivalent Inequalities with the Best Possible Constant Factors}

\begin{theorem}
Suppose that $p>1,\frac{1}{p}+\frac{1}{q}%
=1,\sigma\in \textbf{R}, H(t)\geq 0,$ and
$$
K(\sigma )=\int_{-\infty }^{\infty }H(t)|t|^{\sigma -1}dt\in \mathbf{R}_{+}.
$$
 If $f(x),\:g(y)\geq 0,$ such that $$
0<\int_{-\infty }^{\infty }|x|^{p(1-\delta \sigma )-1}f^{p}(x)dx<\infty$$ and $$ 0<\int_{0}^{\infty }y^{q(1-\sigma )-1}g^{q}(y)dy<\infty,$$
then we obtain the following equivalent inequalities:%
\begin{eqnarray}
I &:&=\int_{-\infty }^{\infty }\int_{-\infty }^{\infty }H(x^{\delta
}y)f(x)g(y)dxdy   \nonumber \\
&<&K(\sigma )\left[ \int_{-\infty }^{\infty }|x|^{p(1-\delta \sigma
)-1}f^{p}(x)dx\right] ^{\frac{1}{p}}\left[ \int_{0}^{\infty }y^{q(1-\sigma
)-1}g^{q}(y)dy\right] ^{\frac{1}{q}},  \label{88}
\end{eqnarray}%
\begin{eqnarray}
J &=&\int_{-\infty }^{\infty }|y|^{p\sigma -1}\left( \int_{-\infty }^{\infty
}H(x^{\delta }y)f(x)dx\right) ^{p}dy   \nonumber \\
&<&K^{p}(\sigma )\int_{-\infty }^{\infty }|x|^{p(1-\delta \sigma
)-1}f^{p}(x)dx,  \label{89}
\end{eqnarray}%
where the constant factors $K(\sigma )$ and $K^{p}(\sigma )$ are the best
possible.
\end{theorem}

In particular, for $\delta =1,$ we have%
\begin{eqnarray}
I &:&=\int_{-\infty }^{\infty }\int_{-\infty }^{\infty }H(xy)f(x)g(y)dxdy
\nonumber \\
&<&K(\sigma )\left[ \int_{-\infty }^{\infty }|x|^{p(1-\sigma )-1}f^{p}(x)dx%
\right] ^{\frac{1}{p}}\left[ \int_{0}^{\infty }y^{q(1-\sigma )-1}g^{q}(y)dy%
\right] ^{\frac{1}{q}},  \label{90}
\end{eqnarray}%
\begin{eqnarray}
J &=&\int_{-\infty }^{\infty }|y|^{p\sigma -1}\left( \int_{-\infty }^{\infty
}H(xy)f(x)dx\right) ^{p}dy  \nonumber \\
&<&K^{p}(\sigma )\int_{-\infty }^{\infty }|x|^{p(1-\sigma )-1}f^{p}(x)dx.
\label{91}
\end{eqnarray}

\begin{proof}
 We shall initially prove that (\ref{86}) preserves the form of
strict inequality for any $y\in \mathbf{R}\backslash \{0\}$.
There exist two constants $A$ and $B,$ such that they are not both
zero, and (cf. \cite{YM23})%
$$
A\frac{|x|^{(1-\delta \sigma )p/q}}{|y|^{1-\sigma }}f^{p}(x)=B\frac{%
|y|^{(1-\sigma )q/p}}{|x|^{1-\sigma }}\text{\,\, a. e.\,\, in\,\, }\mathbf{R}.
$$
If $A=0,$ then $B=0,$ which is impossible. Suppose that $A\neq 0.$ Then it
follows that%
$$
|x|^{p(1-\delta \sigma )-1}f^{p}(x)=|y|^{(1-\sigma )q}\frac{B}{A|x|}\text{\,\,
\,\,a. e.\,\, in\,\, }\mathbf{R},
$$
which contradicts the fact that $$0<\int_{-\infty }^{\infty }|x|^{p(1-\delta
\sigma )-1}f^{p}(x)dx<\infty,$$ in virtue of $$\int_{-\infty }^{\infty }\frac{1}{|x|}%
dx=\infty .$$ Hence, both (\ref{86}) and (\ref{87}) preserve the forms of strict
inequality, and thus we obtain (\ref{89}).

By H\"{o}lder's inequality (cf. \cite{YM23}), we find
\begin{eqnarray}
I &=&\int_{-\infty }^{\infty }\left( |y|^{\sigma -\frac{1}{p}}\int_{-\infty
}^{\infty }H(x^{\delta }y)f(x)dx\right) (|y|^{\frac{1}{p}-\sigma }g(y))dy
 \nonumber \\
&\leq &J^{\frac{1}{p}}\left[ \int_{-\infty }^{\infty }|y|^{q(1-\sigma
)-1}g^{q}(y)dy\right] ^{\frac{1}{q}}.  \label{92}
\end{eqnarray}%
Then by (\ref{89}), we have (\ref{88}). On the other hand, assuming that (%
\ref{88}) is valid, we set%
$$
g(y):=|y|^{p\sigma -1}\left( \int_{-\infty }^{\infty }H(x^{\delta
}y)f(x)dx\right) ^{p-1},\ y\in \mathbf{R}.
$$
Then we get $$J=\int_{-\infty }^{\infty }|y|^{q(1-\sigma )-1}g^{q}(y)dy.$$ By
(\ref{85}), in view of $$
0<\int_{-\infty }^{\infty }|x|^{p(1-\delta \sigma )-1}f^{p}(x)dx<\infty ,
$$
it follows that $J<\infty .$

If $J=0,$ then, (\ref{89}) is trivially valid;
if $J>0,$ then by (\ref{88}), we have
\begin{eqnarray*}
0 &<&\int_{-\infty }^{\infty }|y|^{q(1-\sigma )-1}g^{q}(y)dy=J=I \\
&<&K(\sigma )\left[ \int_{-\infty }^{\infty }|x|^{p(1-\delta \sigma
)-1}f^{p}(x)dx\right] ^{\frac{1}{p}}\left[ \int_{-\infty }^{\infty
}|y|^{q(1-\sigma )-1}g^{q}(y)dy\right] ^{\frac{1}{q}},
\end{eqnarray*}%
$$
J^{\frac{1}{p}}=\left[ \int_{-\infty }^{\infty }y^{q(1-\sigma
)-1}g^{q}(y)dy\right] ^{\frac{1}{p}}
<K(\sigma )\left[ \int_{-\infty }^{\infty }|x|^{p(1-\delta \sigma
)-1}f^{p}(x)dx\right] ^{\frac{1}{p}},
$$
and then (\ref{89}) follows, which is equivalent to (\ref{88}).

For any $n\in \mathbf{N}$, we define two sets $$E_{\delta }:=\{x\in \mathbf{R}%
;|x|^{\delta }\geq 1\},\ \ E_{\delta }^{+}:=\{x\in \mathbf{R}_{+};x^{\delta
}\geq 1\},$$ and the functions $f_{n}(x)$, $g_{n}(y)$ as follows:%
$$
f_{n}(x):=\left\{
\begin{array}{c}
0,\ x\in \mathbf{R}\backslash E_{\delta } \\
|x|^{\delta (\sigma -\frac{1}{np})-1},\ x\in E_{\delta }%
\end{array}%
\ \ \ \ \right. g_{n}(y):=\left\{
\begin{array}{c}
|y|^{\sigma +\frac{1}{nq}-1},\ y\in \lbrack -1,1] \\
0,\ \ y\in \mathbf{R}\backslash \lbrack -1,1]%
\end{array}%
\right.
$$
Then we find%
\begin{eqnarray*}
L_{n} &:&=\left[ \int_{-\infty }^{\infty }|x|^{p(1-\delta \sigma
)-1}f_{n}^{p}(x)dx\right] ^{\frac{1}{p}}\left[ \int_{-\infty }^{\infty
}|y|^{q(1-\sigma )-1}g_{n}^{q}(y)dy\right] ^{\frac{1}{q}} \\
&=&\left( \int_{E_{\delta }}|x|^{-\frac{\delta }{n}-1}dx\right) ^{\frac{1}{p}%
}\left( \int_{-1}^{1}|y|^{\frac{1}{n}-1}dy\right) ^{\frac{1}{q}} \\
&=&\left( 2\int_{E_{\delta }^{+}}x^{-\frac{\delta }{n}-1}dx\right) ^{\frac{1%
}{p}}\left( 2\int_{0}^{1}y^{\frac{1}{n}-1}dy\right) ^{\frac{1}{q}}=2n.
\end{eqnarray*}

Setting $Y=-y,$ we obtain%
\begin{eqnarray*}
I(x) &:&=\int_{-1}^{1}H(x^{\delta }y)|y|^{\sigma +\frac{1}{nq}-1}dy \\
&=&\int_{-1}^{1}H((-x)^{\delta }Y)|Y|^{\sigma +\frac{1}{nq}-1}dY=I(-x),
\end{eqnarray*}%
and then $I(x)$ is an even function. For $x>0,$ we find
\begin{eqnarray*}
I(x) &=&\int_{-1}^{0}H(x^{\delta }y)(-y)^{\sigma +\frac{1}{nq}%
-1}dy+\int_{0}^{1}H(x^{\delta }y)y^{\sigma +\frac{1}{nq}-1}dy \\
&=&x^{-\delta \sigma -\frac{\delta }{nq}}\int_{0}^{x^{\delta
}}(H(-t)+H(t))t^{\sigma +\frac{1}{nq}-1}dt.
\end{eqnarray*}%
By the above results and Fubini's theorem, it follows that%
\begin{eqnarray*}
I_{n} &:&=\int_{-\infty }^{\infty }\int_{-\infty }^{\infty }H(x^{\delta
}y)f_{n}(x)g_{n}(y)dxdy \\
&=&\int_{E_{\delta }}|x|^{\delta (\sigma -\frac{1}{np})-1}\left(
\int_{-1}^{1}H(x^{\delta }y)|y|^{\sigma +\frac{1}{nq}-1}dy\right) dx \\
&=&\int_{E_{\delta }}|x|^{\delta (\sigma -\frac{1}{np})-1}I(x)dx=2\int_{E_{%
\delta }^{+}}x^{\delta (\sigma -\frac{1}{np})-1}I(x)dx\\
&=&2\int_{E_{\delta }^{+}}x^{-\frac{\delta }{n}-1}\left(
\int_{0}^{1}(H(-t)+H(t))t^{\sigma +\frac{1}{nq}-1}dt\right) dx \\
&&+2\int_{E_{\delta }^{+}}x^{-\frac{\delta }{n}-1}\left( \int_{1}^{x^{\delta
}}(H(-t)+H(t))t^{\sigma +\frac{1}{nq}-1}dt\right) dx\\
&=&2n\left( \int_{0}^{1}(H(-t)+H(t))t^{\sigma +\frac{1}{nq}-1}dt\right) \\
&&+2\int_{1}^{\infty }\left( \int_{\{x>0;x^{\delta }\geq t\}}x^{-\frac{%
\delta }{n}-1}dx\right) (H(-t)+H(t))t^{\sigma +\frac{1}{nq}-1}dt\\
&=&2n\left[ \int_{0}^{1}(H(-t)+H(t))t^{\sigma +\frac{1}{nq}-1}dt\right.
\left. +\int_{1}^{\infty }(H(-t)+H(t))t^{\sigma -\frac{1}{np}-1}dt\right] .
\end{eqnarray*}%

If there exists a positive number $k\leq K(\sigma ),$ such that (\ref{88})
is still valid when replacing $K(\sigma )$ by $k,$ then in particular, it
follows that $$\frac{1}{2n}I_{n}<k\frac{1}{2n}L_{n},$$ and%
$$
\int_{0}^{1}(H(-t)+H(t))t^{\sigma +\frac{1}{nq}-1}dt+\int_{1}^{\infty
}(H(-t)+H(t))t^{\sigma -\frac{1}{np}-1}dt<k.
$$
Since both $$\{(H(-t)+H(t))t^{\sigma +\frac{1}{nq}-1}\}_{n=1}^{\infty }\ (t\in
(0,1])$$ and
$$
\{(H(-t)+H(t))t^{\sigma -\frac{1}{np}-1}\}_{n=1}^{\infty }\ (t\in (1,\infty ))
$$
are non-negative and increasing, then by Levi's theorem (cf. \cite{YM24}), it
follows that
\begin{eqnarray*}
&&K(\sigma )=\int_{0}^{1}(H(-t)+H(t))t^{\sigma -1}dt+\int_{1}^{\infty
}(H(-t)+H(t))t^{\sigma -1}dt\\
&=&\lim_{n\rightarrow \infty }\left( \int_{0}^{1}(H(-t)+H(t))t^{\sigma +%
\frac{1}{nq}-1}dt+\int_{1}^{\infty }(H(-t)+H(t))t^{\sigma -\frac{1}{np}%
-1}dt\right) \\
&\leq &k,
\end{eqnarray*}%
and thus $k=K(\sigma )$ is the best possible constant factor of (\ref{88}).

The constant factor in (\ref{89}) is still the best possible. Otherwise, we
would reach a contradiction by (\ref{92}).
This completes the proof of the theorem.
\end{proof}

\begin{theorem}
 Replacing $p>1$ by $0<p<1$ in Theorem 1, we obtain the
equivalent reverses of (\ref{88}) and (\ref{89}). If there exists a constant
$\delta _{\:0}>0,$ such that for any $\widetilde{\sigma }\in (\sigma -\delta
_{\:0},\sigma ],$  $$K(\widetilde{\sigma })=\int_{-\infty }^{\infty }H(t)|t|^{%
\widetilde{\sigma }-1}dt\in \mathbf{R}_{+},$$ then the constant factors in
the reverses of (\ref{88}) and (\ref{89}) are the best possible.
\end{theorem}

\begin{proof}
By Lemma 3 and the reverse of H\"{o}lder's inequality, we get
the reverses of (\ref{88}), (\ref{89}) and (\ref{92}). Similarly, we
can set $g(y)$ as in Theorem 5, and prove that the reverses of (\ref{88}) and (%
\ref{89}) are equivalent.

For $n>\frac{2}{\delta _{\:0}|q|}\ (n\in \mathbf{N}),$ we set $f_{n}(x)$ and $%
g_{n}(y)$ as in Theorem 1. If there exists a positive number $k\geq K(\sigma ),$
such that the reverse of (\ref{88}) is valid when replacing $K(\sigma )$ by $%
k,$ then it follows that $$\frac{1}{2n}I_{n}>k\frac{1}{2n}%
L_{n},$$ and%
\begin{equation}
\int_{0}^{1}(H(-t)+H(t))t^{\sigma +\frac{1}{nq}-1}dt+\int_{1}^{\infty
}(H(-t)+H(t))t^{\sigma -\frac{1}{np}-1}dt>k.  \label{93}
\end{equation}

Since $$\{(H(-t)+H(t))t^{\sigma -\frac{1}{np}-1}\}_{n=1}^{\infty }\ (t\in
(1,\infty ))$$ is still a non-negative and increasing sequence, then by Levi's theorem, it
follows that
$$
\lim_{n\rightarrow \infty }\int_{1}^{\infty }(H(-t)+H(t))t^{\sigma -\frac{1}{%
np}-1}dt=\int_{1}^{\infty }(H(-t)+H(t))t^{\sigma -1}dt.
$$
Due to the fact that $$
0\leq (H(-t)+H(t))t^{\sigma +\frac{1}{nq}-1}\leq (H(-t)+H(t))t^{(\sigma -%
\frac{\delta _{\:0}}{2})-1}
$$
$(t\in (0,1],n>\frac{2}{\delta _{\:0}|q|}),$ and $$
0\leq \int_{0}^{1}(H(-t)+H(t))t^{(\sigma -\frac{\delta _{\:0}}{2})-1}dt\leq
K(\sigma -\frac{\delta _{\:0}}{2})<\infty,
$$
 then, by Lebesgue's dominated convergence theorem (cf. \cite{YM24}), it follows
that%
$$
\lim_{n\rightarrow \infty }\int_{0}^{1}(H(-t)+H(t))t^{\sigma +\frac{1}{nq}%
-1}dt=\int_{0}^{1}(H(-t)+H(t))t^{\sigma -1}dt.
$$

In view of the above results and (\ref{93}), we have%
\begin{eqnarray*}
K(\sigma ) &=&\int_{0}^{1}(H(-t)+H(t))t^{\sigma -1}dt+\int_{1}^{\infty
}(H(-t)+H(t))t^{\sigma -1}dt \\
&=&\lim_{n\rightarrow \infty }\left( \int_{0}^{1}(H(-t)+H(t))t^{\sigma +%
\frac{1}{nq}-1}dt\right.\\
&&\left. +\int_{1}^{\infty }(H(-t)+H(t))t^{\sigma -\frac{1}{np}-1}dt\right)
\geq k,
\end{eqnarray*}%
and then $k=K(\sigma )$ is the best possible constant factor in the reverse
of (\ref{88}).

Simiraly, we can prove that the constant factor in the reverse of (%
\ref{89}) is the best possible by using the reverse of (\ref{92}). $\Box $
\end{proof}

\subsection{Yang-Hilbert-Type Integral Inequalities in the Whole Plane with
Multi-Variables}

\begin{theorem}
 Suppose that $p>1,\:\frac{1}{p}+\frac{1}{q}%
=1,\:H(t)\geq 0,\:K(\sigma )\in \mathbf{R}_{+},\:\delta \in \{-1,1\},$ $-\infty
\leq a_{i}<b_{i}\leq \infty ,v_{i}^{\prime }(s)>0\ (s\in
(a_{i},b_{i})),\:v_{i}(a_{i}^{+})=-\infty ,\:v_{i}(b_{i}^{-})=\infty\ (i=1,2).$
If $f(x),g(y)\geq 0,$ such that%
$$
0 <\int_{a_{1}}^{b_{1}}\frac{|v_{1}(x)|^{p(1-\delta \sigma )-1}}{%
(v_{1}^{\prime }(x))^{p-1}}f^{p}(x)dx<\infty$$ and $$
0<\int_{a_{2}}^{b_{2}}\frac{|v_{2}(y)|^{q(1-\sigma )-1}}{(v_{2}^{\prime
}(y))^{q-1}}g^{q}(y)dy<\infty ,
$$
then we have the following equivalent inequalities:%
\begin{eqnarray}
&&\int_{a_{2}}^{b_{2}}\int_{a_{1}}^{b_{1}}H(v_{1}^{\delta
}(x)v_{2}(y))f(x)g(y)dxdy   \nonumber \\
&<&K(\sigma )\left[ \int_{a_{1}}^{b_{1}}\frac{|v_{1}(x)|^{p(1-\delta \sigma
)-1}}{(v_{1}^{\prime }(x))^{p-1}}f^{p}(x)dx\right] ^{\frac{1}{p}}   \nonumber \\
&&\times \left[ \int_{a_{2}}^{b_{2}}\frac{|v_{2}(y)|^{q(1-\sigma )-1}}{%
(v_{2}^{\prime }(y))^{q-1}}g^{q}(y)dy\right] ^{\frac{1}{q}},  \label{94}
\end{eqnarray}%
\begin{eqnarray}
&&\int_{a_{2}}^{b_{2}}\frac{v_{2}^{\prime }(y)}{|v_{2}(y)|^{1-p\sigma }}%
\left( \int_{a_{1}}^{b_{1}}H(v_{1}^{\delta }(x)v_{2}(y))f(x)dx\right) ^{p}dy
 \nonumber \\
&<&K^{p}(\sigma )\int_{a_{1}}^{b_{1}}\frac{|v_{1}(x)|^{p(1-\delta \sigma )-1}%
}{(v_{1}^{\prime }(x))^{p-1}}f^{p}(x)dx,  \label{95}
\end{eqnarray}%
where the constant factors $K(\sigma )$ and $K^{p}(\sigma )$ are the best
possible.
\end{theorem}

\begin{proof}
Setting $x=v_{1}(s),\:y=v_{2}(t)$ in (\ref{88}), we get $%
dx=v_{1}^{\prime }(s)ds,$ $dy=v_{2}^{\prime }(t)dt,$ and%
\begin{eqnarray*}
I &=&\int_{a_{2}}^{b_{2}}\int_{a_{1}}^{b_{1}}H(v_{1}^{\delta
}(s)v_{2}(t))f(v_{1}(s))g(v_{2}(t))v_{1}^{\prime }(s)v_{2}^{\prime }(t)dsdt
\\
&=&\int_{a_{2}}^{b_{2}}\int_{a_{1}}^{b_{1}}H(v_{1}^{\delta
}(s)v_{2}(t))(f(v_{1}(s))v_{1}^{\prime }(s))(g(v_{2}(t))v_{2}^{\prime
}(t))dsdt,\\
I_{1} &:&=\int_{-\infty }^{\infty }|x|^{p(1-\delta \sigma
)-1}f^{p}(x)dx=\int_{a_{1}}^{b_{1}}|v_{1}(s)|^{p(1-\delta \sigma
)-1}f^{p}(v_{1}(s))v_{1}^{\prime }(s)ds, \\
I_{2} &:&=\int_{-\infty }^{\infty }y^{q(1-\sigma
)-1}g^{q}(y)dy=\int_{a_{2}}^{b_{2}}|v_{2}(t)|^{q(1-\sigma
)-1}g^{q}(v_{2}(t))v_{2}^{\prime }(t)dt.
\end{eqnarray*}%
If we set $$F(s)=f(v_{1}(s))v_{1}^{\prime }(s)\ \ \text{and}\ \ G(t)=g(v_{2}(t))v_{2}^{\prime
}(t),$$ we obtain%
$$
f^{p}(v_{1}(s))=(v_{1}^{\prime
}(s))^{-p}F^{p}(s),g^{q}(v_{2}(t))=(v_{2}^{\prime }(t))^{-q}G^{q}(t),
$$
and then it follows that
\begin{eqnarray*}
I &=&\int_{a_{2}}^{b_{2}}\int_{a_{1}}^{b_{1}}H(v_{1}^{\delta
}(s)v_{2}(t))F(s)G(t)dsdt, \\
I_{1} &=&\int_{a_{1}}^{b_{1}}\frac{|v_{1}(s)|^{p(1-\delta \sigma )-1}}{%
(v_{1}^{\prime }(s))^{p-1}}F^{p}(s)ds, \\
I_{2} &=&\int_{a_{2}}^{b_{2}}\frac{|v_{2}(t)|^{q(1-\sigma )-1}}{%
(v_{2}^{\prime }(t))^{q-1}}G^{q}(t)dt.
\end{eqnarray*}

 Substitution of the above results to (\ref{88}), with $$%
s=x,\:t=y,\:F(s)=f(x),\ \text{and}\ G(t)=g(y),$$ we obtain (\ref{94}). Similarly,
we derive (\ref{95}). On the other hand, if we set $$v_{1}(x)=x,\:v_{2}(y)=y,\;
a_{i}=-\infty ,\:b_{i}=\infty $$ in (\ref{94}), we get (\ref{88}). Hence, the
inequalities (\ref{94}) and (\ref{88}) are equivalent. It is evident that the
inequalities (\ref{95}) and (\ref{89}) are equivalent. Hence, the inequalities (%
\ref{94}) and (\ref{95}) are equivalent. Since the constant factors in (\ref%
{88}) and (\ref{89}) are the best possible, it follows that the constant factors in (%
\ref{94}) and (\ref{95}) are also the best possible.
This completes the proof of the theorem.
\end{proof}
\begin{theorem}
 Replacing $p>1$ by $0<p<1$ in Theorem 7, we
obtain the equivalent reverses of (\ref{94}) and (\ref{95}). If there exists a
constant $\delta _{\:0}>0,$ such that for any \mbox{$\widetilde{\sigma }\in (\sigma
-\delta _{\:0},\sigma ],$} $$K(\widetilde{\sigma })=\int_{-\infty }^{\infty
}H(t)|t|^{\widetilde{\sigma }-1}dt\in \mathbf{R}_{+},$$ then the constant
factors in the reverses of (\ref{94}) and (\ref{95}) are the best possible.
\end{theorem}

\begin{remark}
 We list the following $v_{i}(s)\ (i=1,2)$ that satisfy the
conditions of Theorem 7 and Theorem 8:
\end{remark}

(a) $v_{i}(s)=s^{\gamma },\ s\in (-\infty ,\infty )\ (\gamma \in \{a;a=\frac{1}{%
2k-1},2k+1\ (k\in \mathbf{N)}\}),$ satisfying \mbox{$v_{i}^{\prime }(s)=\gamma
s^{\gamma -1}>0;$}

(b) $v_{i}(s)=\tan ^{\gamma }s,\ s\in (-\frac{\pi }{2},\frac{\pi }{2})\ (\gamma
\in \{a;a=\frac{1}{2k-1},2k+1\ (k\in \mathbf{N)}\}),$ satisfying $%
v_{i}^{\prime }(s)=\gamma \tan ^{\gamma -1}s\sec ^{2}s>0;$

(c) $v_{i}(s)=\ln ^{\gamma }s,\ s\in (0,\infty )\ (\gamma \in \{a;a=\frac{1}{2k-1%
},2k+1\ (k\in \mathbf{N)}\}),$ satisfying \mbox{$v_{i}^{\prime }(s)=\frac{\gamma }{s}%
\ln ^{\gamma -1}s>0;$}

(d) $v_{i}(s)=(e^{|s|}-1)sgn(s),\ s\in (-\infty ,\infty ),$ satisfying $%
v_{i}^{\prime }(s)=e^{|s|}>0.$

\begin{definition}
 If $\lambda \in \mathbf{R},K_{\lambda }(x,y)$ is a
non-negative measurable function in\textbf{\ }$\mathbf{R}^{2}$, satisfying
$$K_{\lambda
}(tx,ty)=|t|^{-\lambda }K_{\lambda }(x,y),$$
for any $t\in \mathbf{R}\backslash \{0\},\:x,\:y\in \mathbf{R,}$ then $K_{\lambda }(x,y)$
is said to be the homogeneous function of degree $-\lambda $ in $\mathbf{R}^{2}$.
\end{definition}

In particular, by Theorem 7 and Theorem 8, with $\delta =1$, we obtain the
following integral inequalities with the non-homogeneous kernel:

\begin{corollary}
 Suppose that $p>1,\frac{1}{p}+\frac{1}{q}=1,\sigma\in \textbf{R}, H(t)\geq
0,$%
$$
K(\sigma )=\int_{-\infty }^{\infty }H(t)|t|^{\sigma -1}dt\in \mathbf{R}_{+},
$$
$-\infty \leq a_{i}<b_{i}\leq \infty ,\:v_{i}^{\prime }(s)>0\ (s\in
(a_{i},b_{i})),v_{i}(a_{i}^{+})=-\infty ,\:v_{i}(b_{i}^{-})=\infty\ (i=1,2).$
If $f(x),\:g(y)\geq 0,$ such that%
$$
0<\int_{a_{1}}^{b_{1}}\frac{|v_{1}(x)|^{p(1-\sigma )-1}}{(v_{1}^{\prime
}(x))^{p-1}}f^{p}(x)dx<\infty $$ and $$
0<\int_{a_{2}}^{b_{2}}\frac{|v_{2}(y)|^{q(1-\sigma )-1}}{(v_{2}^{\prime
}(y))^{q-1}}g^{q}(y)dy<\infty,
$$
then we have the following equivalent inequalities:%

\begin{eqnarray}
&&\int_{a_{2}}^{b_{2}}\int_{a_{1}}^{b_{1}}H(v_{1}(x)v_{2}(y))f(x)g(y)dxdy\nonumber\\
&<&K(\sigma )\left[ \int_{a_{1}}^{b_{1}}\frac{|v_{1}(x)|^{p(1-\sigma )-1}}{%
(v_{1}^{\prime }(x))^{p-1}}f^{p}(x)dx\right] ^{\frac{1}{p}}   \nonumber \\
&&\times \left[ \int_{a_{2}}^{b_{2}}\frac{|v_{2}(y)|^{q(1-\sigma )-1}}{%
(v_{2}^{\prime }(y))^{q-1}}g^{q}(y)dy\right] ^{\frac{1}{q}},  \label{96}
\end{eqnarray}%
\begin{eqnarray}
&&\int_{a_{2}}^{b_{2}}\frac{v_{2}^{\prime }(y)}{|v_{2}(y)|^{1-p\sigma }}%
\left( \int_{a_{1}}^{b_{1}}H(v_{1}(x)v_{2}(y))f(x)dx\right) ^{p}dy   \nonumber \\
&<&K^{p}(\sigma )\int_{a_{1}}^{b_{1}}\frac{|v_{1}(x)|^{p(1-\sigma )-1}}{%
(v_{1}^{\prime }(x))^{p-1}}f^{p}(x)dx,  \label{97}
\end{eqnarray}%
where the constant factors $K(\sigma )$ and $K^{p}(\sigma )$ are the best
possible.
\end{corollary}
Replacing $p>1$ by $0<p<1$ in Corollary 12, we derive the equivalent reverses
of (\ref{96}) and (\ref{97}). If there exists a constant $\delta _{\:0}>0,$
such that for any $\widetilde{\sigma }\in (\sigma -\delta _{\:0},\sigma ],$%
\[
K(\widetilde{\sigma })=\int_{-\infty }^{\infty }H(t)|t|^{\widetilde{\sigma }%
-1}dt\in \mathbf{R}_{+},
\]
then the constant factors in the reverses of (\ref{96}) and (\ref{97}) are
the best possible.

In particular, for $\delta =-1$ in Theorem 7 and Theorem 8, setting $%
H(t)=K_{\lambda }(1,t)$ (cf. Definition 10),   we get%
\[
H(\frac{v_{2}(y)}{v_{1}(x)})=K_{\lambda }(1,\frac{v_{2}(y)}{v_{1}(x)}%
)=|v_{1}(x)|^{\lambda }K_{\lambda }(v_{1}(x),v_{2}(y)).
\]
Replacing $f(x)$ by $|v_{1}(x)|^{-\lambda
}f(x),$ it follows that $%
|v_{1}(x)|^{p(1+\sigma )-1}f^{p}(x)$ is replaced by
\[| v_{1}(x)|^{p(1+\sigma )-1}[|v_{1}(x)|^{-\lambda
}f(x)]^{p}=|v_{1}(x)|^{p(1-\mu )-1}f^{p}(x),
\]
and we have the following integral inequalities with the
homogeneous kernel:

\begin{corollary}
 Let $p>1,\frac{1}{p}+\frac{1}{q}=1,\:\mu
,\:\sigma \in \mathbf{R},\:\mu +\sigma =\lambda ,\:K_{\lambda }(x,y)$ is a
homogeneous function in $\mathbf{R}^{2}$ of degree $-\lambda ,$%
\begin{eqnarray*}
K_{\lambda }(\sigma ):=\int_{-\infty }^{\infty }K_{\lambda }(1,t)|t|^{\sigma
-1}dt\in \mathbf{R}_{+},
\end{eqnarray*}
$0\leq a_{i}<b_{i}\leq \infty ,v_{i}^{\prime }(s)>0$ $(s\in
(a_{i},b_{i})),\:v_{i}(a_{i}^{+})=-\infty ,v_{i}(b_{i}^{-})=\infty\ (i=1,2).$
If $f(x),g(y)\geq 0,$ such that%
$$
0<\int_{a_{1}}^{b_{1}}\frac{|v_{1}(x)|^{p(1-\mu )-1}}{(v_{1}^{\prime
}(x))^{p-1}}f^{p}(x)dx<\infty$$ and $$
0<\int_{a_{2}}^{b_{2}}\frac{|v_{2}(y)|^{q(1-\sigma )-1}}{(v_{2}^{\prime
}(y))^{q-1}}g^{q}(y)dy<\infty,$$
then we have the following equivalent inequalities:%
\begin{eqnarray}
&&\int_{a_{2}}^{b_{2}}\int_{a_{1}}^{b_{1}}K_{\lambda
}(v_{1}(x),v_{2}(y))f(x)g(y)dxdy\nonumber\\
&<&K_{\lambda }(\sigma )\left[ \int_{a_{1}}^{b_{1}}\frac{|v_{1}(x)|^{p(1-\mu
)-1}}{(v_{1}^{\prime }(x))^{p-1}}f^{p}(x)dx\right] ^{\frac{1}{p}}   \nonumber \\
&&\times \left[ \int_{a_{2}}^{b_{2}}\frac{|v_{2}(y)|^{q(1-\sigma )-1}}{%
(v_{2}^{\prime }(y))^{q-1}}g^{q}(y)dy\right] ^{\frac{1}{q}},  \label{98}
\end{eqnarray}%
\begin{eqnarray}
&&\int_{a_{2}}^{b_{2}}\frac{v_{2}^{\prime }(y)}{|v_{2}(y)|^{1-p\sigma }}%
\left( \int_{a_{1}}^{b_{1}}K_{\lambda }(v_{1}(x),v_{2}(y))f(x)dx\right)
^{p}dy   \nonumber \\
&<&K_{\lambda }^{p}(\sigma )\int_{a_{1}}^{b_{1}}\frac{|v_{1}(x)|^{p(1-\mu
)-1}}{(v_{1}^{\prime }(x))^{p-1}}f^{p}(x)dx,  \label{99}
\end{eqnarray}%
where the constant factors $K_{\lambda }(\sigma )$ and $K_{\lambda
}^{p}(\sigma )$ are the best possible.
\end{corollary}

Replacing $p>1$ by $0<p<1$ in the above cases, we obtain the equivalent
reverses of (\ref{98}) and (\ref{99}). If there exists a constant $\delta
_{\:0}>0,$ such that for any $\widetilde{\sigma }\in (\sigma -\delta
_{\:0},\sigma ],$ $$K_{\lambda }(\widetilde{\sigma })=\int_{-\infty }^{\infty
}K_{\lambda }(1,t)|t|^{\widetilde{\sigma }-1}dt\in \mathbf{R}_{+},$$ then the
constant factors in the reverses of (\ref{98}) and (\ref{99}) are the best
possible.

Setting $a_{i}=-\infty ,\:b_{i}=\infty\ (i=1,2),\:v_{1}(x)=x,\:v_{2}(y)=y$ in
Corollary 13, we obtain the following corollary:

\begin{corollary}
Let $p>1,\:\frac{1}{p}+\frac{1}{q}%
=1,\:\mu,\:\sigma \in \mathbf{R},\:\mu +\sigma =\lambda ,$ $K_{\lambda }(x,y)$ is
a homogeneous function in $\mathbf{R}^{2}$ of degree $-\lambda,$
\begin{eqnarray*}
K_{\lambda }(\sigma )=\int_{-\infty }^{\infty }K_{\lambda }(1,t)t^{\sigma
-1}dt\in \mathbf{R}_{+}.
\end{eqnarray*}
If $f(x),g(y)\geq 0,$
$$
0<\int_{0}^{\infty }x^{p(1-\mu )-1}f^{p}(x)dx<\infty $$ and 
$$ 0<\int_{0}^{\infty
}g^{q(1-\sigma )-1}(y)dy<\infty ,
$$
then we have the following equivalent inequalities:%
\begin{eqnarray}
&&\int_{-\infty }^{\infty }\int_{-\infty }^{\infty }K_{\lambda
}(x,y)f(x)g(y)dxdy   \nonumber\\
&<&K_{\lambda }(\sigma )\left[ \int_{-\infty }^{\infty }x^{p(1-\mu
)-1}f^{p}(x)dx\right] ^{\frac{1}{p}}\left[ \int_{-\infty }^{\infty
}y^{q(1-\sigma )-1}g^{q}(y)dy\right] ^{\frac{1}{q}},  \label{100}
\end{eqnarray}%
\begin{equation}
\int_{-\infty }^{\infty }y^{p\sigma -1}\left( \int_{-\infty }^{\infty
}K_{\lambda }(x,y)f(x)dx\right) ^{p}dy
<K_{\lambda }^{p}(\sigma )\int_{0}^{\infty }x^{p(1-\mu )-1}f^{p}(x)dx,
\label{101}
\end{equation}%
where the constant factors $K_{\lambda }(\sigma )$ and $K_{\lambda
}^{p}(\sigma )$ are the best possible.
\end{corollary}

Replacing $p>1$ by $0<p<1$ in the above cases, we get the equivalent
reverses of (\ref{100}) and (\ref{101}). If there exists a constant $\delta
_{\:0}>0,$ such that for any $\widetilde{\sigma }\in (\sigma -\delta
_{\:0},\sigma ],$
\[
K_{\lambda }(\widetilde{\sigma })=\int_{-\infty }^{\infty }K_{\lambda
}(1,t)|t|^{\widetilde{\sigma }-1}dt\in \mathbf{R}_{+},
\]
then the constant factors in the reverses of (\ref{100}) and (\ref{101}) are
the best possible.

\begin{remark}
 It is evident that (\ref{90}) and (\ref{100}) are
equivalent for $H(t)=K_{\lambda }(1,t).$ The same holds for (\ref{91}) and (\ref{101}).
\end{remark}

\subsection{Hardy-Type Integral Inequalities in the Whole Plane}

In the following two sections, if the constant factors in the inequalities
(operator inequalities) are related to $K^{(1)}(\sigma)$ ( or $%
K_{\lambda}^{(1)}(\sigma)$), then we shall call them Hardy-type
inequalities (operator) of the first kind; if the constant factors in the inequalities
(operator inequalities) are related to $K^{(2)}(\sigma)$ ( or $%
K_{\lambda}^{(2)}(\sigma)$), then we shall call them Hardy-type
inequalities (operator) of the second kind.

\bigskip If $H(t)=0\ (|t|>1)$, then $H(xy)=0\ (|x|>\frac{1}{|y|}>0),$ and%
\[
K(\sigma )=\int_{-\infty }^{\infty }H(t)|t|^{\sigma
-1}dt=\int_{-1}^{1}H(t)|t|^{\sigma -1}dt.
\]
Set%
\begin{equation}
K^{(1)}(\sigma ):=\int_{-1}^{1}H(t)|t|^{\sigma
-1}dt=\int_{0}^{1}(H(-t)+H(t))t^{\sigma -1}dt.  \label{102}
\end{equation}%
\ Then, by Theorem 7 and Theorem 8 ($\delta =1$), we have the following
Hardy-type integral inequalities of the first kind, with non-homogeneous kernel in the whole
plane:

\begin{corollary}
 Suppose that $p>1,\:\frac{1}{p}+\frac{1}{q}%
=1,\:H(t)\geq 0,\sigma\in \textbf{R},$%
$$
K^{(1)}(\sigma )=\int_{-1}^{1}H(t)|t|^{\sigma -1}dt\in \mathbf{R}_{+}.
$$
If $f(x),g(y)\geq 0,$%
$$
0<\int_{-\infty }^{\infty }|x|^{p(1-\sigma )-1}f^{p}(x)dx<\infty
$$ and $$ 0<\int_{-\infty }^{\infty }|y|^{q(1-\sigma )-1}g^{q}(y)dy<\infty ,
$$
then we have the following equivalent inequalities:%
\begin{eqnarray}
&&\int_{-\infty }^{\infty }\left( \int_{-\frac{1}{|y|}}^{\frac{1}{|y|}%
}H(xy)f(x)dx\right) g(y)dy
=\int_{0}^{\infty }\left( \int_{-\frac{1}{|x|}}^{\frac{1}{|x|}%
}H(xy)g(y)dy\right) f(x)dx\nonumber\\
&<&K^{(1)}(\sigma )\left[ \int_{-\infty }^{\infty }|x|^{p(1-\sigma
)-1}f^{p}(x)dx\right] ^{\frac{1}{p}}\left[ \int_{-\infty }^{\infty
}|y|^{q(1-\sigma )-1}g^{q}(y)dy\right] ^{\frac{1}{q}},  \label{103}
\end{eqnarray}%
\begin{equation}
\int_{-\infty }^{\infty }|y|^{p\sigma -1}\left( \int_{-\frac{1}{|y|}}^{%
\frac{1}{|y|}}H(xy)f(x)dx\right) ^{p}dy
<(K^{(1)}(\sigma ))^{p}\int_{-\infty }^{\infty }|x|^{p(1-\sigma
)-1}f^{p}(x)dx;  \label{104}
\end{equation}%
where the constant factors $K^{(1)}(\sigma )$ and $(K^{(1)}(\sigma ))^{p}$
are the best possible.
\end{corollary}

Replacing $p>1$ by $0<p<1$ in the above inequalities, we obtain the equivalent
reverses of (\ref{103}) and (\ref{104}). If there exists a constant $\delta
_{\:0}>0,$ such that for any $\widetilde{\sigma }\in (\sigma -\delta
_{\:0},\sigma ],$
\[
K^{(1)}(\widetilde{\sigma })=\int_{-1}^{1}H(t)|t|^{\widetilde{\sigma }%
-1}dt\in \mathbf{R}_{+},
\]
then the constant factors in the reverses of (\ref{103}) and (\ref{104}) are
the best possible.

If we have $H(t)=0\ (|t|>1)$ in Corollary 12, then $$H(v_{1}(x)v_{2}(y))=0(|v_{1}(x)|>%
\frac{1}{|v_{2}(y)|}>0),$$ and thus we derive the following general results:

\begin{corollary}
 Suppose that $p>1,\frac{1}{p}+\frac{1}{q}=1,\:H(t)\geq
0,\sigma\in \textbf{R},$%
\begin{eqnarray*}
K^{(1)}(\sigma )=\int_{-1}^{1}H(t)|t|^{\sigma -1}dt\in \mathbf{R}_{+},
\end{eqnarray*}
$-\infty \leq a_{i}<b_{i}\leq \infty ,\:v_{i}^{\prime }(s)>0\ (s\in
(a_{i},b_{i})),\:v_{i}(a_{i}^{+})=-\infty ,\:v_{i}(b_{i}^{-})=\infty\ (i=1,2).$
If $f(x),g(y)\geq 0,$ such that%
$$
0<\int_{a_{1}}^{b_{1}}\frac{|v_{1}(x)|^{p(1-\sigma )-1}}{(v_{1}^{\prime
}(x))^{p-1}}f^{p}(x)dx<\infty$$ and $$
0<\int_{a_{2}}^{b_{2}}\frac{|v_{2}(y)|^{q(1-\sigma )-1}}{(v_{2}^{\prime
}(y))^{q-1}}g^{q}(y)dy<\infty ,
$$
then we obtain the following equivalent inequalities:%
\begin{eqnarray}
&&\int_{a_{2}}^{b_{2}}\left( \int_{v_{1}^{-1}(\frac{-1}{|v_{2}(y)|}%
)}^{v_{1}^{-1}(\frac{1}{|v_{2}(y)|})}H(v_{1}(x)v_{2}(y))f(x)dx\right) g(y)dy
 \nonumber \\
&<&K^{(1)}(\sigma )\left[ \int_{a_{1}}^{b_{1}}\frac{|v_{1}(x)|^{p(1-\sigma
)-1}}{(v_{1}^{\prime }(x))^{p-1}}f^{p}(x)dx\right] ^{\frac{1}{p}}   \nonumber \\
&&\times \left[ \int_{a_{2}}^{b_{2}}\frac{|v_{2}(y)|^{q(1-\sigma )-1}}{%
(v_{2}^{\prime }(y))^{q-1}}g^{q}(y)dy\right] ^{\frac{1}{q}},  \label{105}
\end{eqnarray}%
\begin{eqnarray}
&&\int_{a_{2}}^{b_{2}}\frac{v_{2}^{\prime }(y)}{|v_{2}(y)|^{1-p\sigma }}%
\left( \int_{v_{1}^{-1}(\frac{-1}{|v_{2}(y)|})}^{v_{1}^{-1}(\frac{1}{%
|v_{2}(y)|})}H(v_{1}(x)v_{2}(y))f(x)dx\right) ^{p}dy   \nonumber \\
&<&(K^{(1)}(\sigma ))^{p}\int_{a_{1}}^{b_{1}}\frac{|v_{1}(x)|^{p(1-\sigma
)-1}}{(v_{1}^{\prime }(x))^{p-1}}f^{p}(x)dx,  \label{106}
\end{eqnarray}%
where the constant factors $K^{(1)}(\sigma )$ and $(K^{(1)}(\sigma ))^{p}$
are the best possible.
\end{corollary}

Replacing $p>1$ by $0<p<1$ in the above inequalities, we get the equivalent
reverses of (\ref{105}) and (\ref{106}). If there exists a constant $\delta
_{\:0}>0,$ such that for any $\widetilde{\sigma }\in (\sigma -\delta
_{\:0},\sigma ],$
\[
K^{(1)}(\widetilde{\sigma })=\int_{-1}^{1}H(t)|t|^{\widetilde{\sigma }%
-1}dt\in \mathbf{R}_{+},
\]
then the constant factors in the reverses of (\ref{105}) and (\ref{106}) are
the best possible.

If $H(t)=0\ (0<|t|<1)$, then $H(xy)=0\ (0<|x|<\frac{1}{|y|}),$ and%
\begin{eqnarray}
K(\sigma ) &=&\int_{-\infty }^{\infty }H(t)|t|^{\sigma -1}dt=\int_{-\infty
}^{-1}H(t)(-t)^{\sigma -1}dt+\int_{1}^{\infty }H(t)t^{\sigma -1}dt  \nonumber
\\
&=&\int_{1}^{\infty }(H(-t)+H(t))t^{\sigma -1}dt.  \label{107}
\end{eqnarray}%
If we set $$E_{y}:=\left\{x\in \mathbf{R};x\geq \frac{1}{|y|}, \text{or}\ x\leq \frac{-1}{%
|y|}\right\},$$ and
\[
K^{(2)}(\sigma ):=\int_{E_{1}}H(t)|t|^{\sigma -1}dt=\int_{1}^{\infty
}(H(-t)+H(t))t^{\sigma -1}dt,
\]
then,
by Theorem 7 and Theorem 8 ($\delta =1$), we obtain the following Hardy-type integral inequalities
of the second kind with the non-homogeneous kernel in the
whole plane:

\begin{corollary}
Let $p>1,\:\frac{1}{p}+\frac{1}{q}=1,\:H(t)\geq
0,\sigma\in \textbf{R},$%
\begin{eqnarray*}
K^{(2)}(\sigma )=\int_{1}^{\infty }(H(-t)+H(t))t^{\sigma -1}dt\in \mathbf{R}%
_{+}.
\end{eqnarray*}
If $f(x),g(y)\geq 0,$%
$$
0<\int_{-\infty }^{\infty }|x|^{p(1-\sigma )-1}f^{p}(x)dx<\infty
$$ and $$ 0<\int_{-\infty }^{\infty }|y|^{q(1-\sigma )-1}g^{q}(y)dy<\infty ,
$$
then we have the following equivalent inequalities:%
\begin{eqnarray}
&&\int_{-\infty }^{\infty }\left( \int_{E_{y}}H(xy)f(x)dx\right) g(y)dy
\nonumber\\
&<&K^{(2)}(\sigma )\left[ \int_{-\infty }^{\infty }|x|^{p(1-\sigma
)-1}f^{p}(x)dx\right] ^{\frac{1}{p}}\left[ \int_{-\infty }^{\infty
}|y|^{q(1-\sigma )-1}g^{q}(y)dy\right] ^{\frac{1}{q}},  \label{108}
\end{eqnarray}%
\begin{equation}
\int_{-\infty }^{\infty }|y|^{p\sigma -1}\left(
\int_{E_{y}}H(xy)f(x)dx\right) ^{p}dy
<(K^{(2)}(\sigma ))^{p}\int_{-\infty }^{\infty }|x|^{p(1-\sigma
)-1}f^{p}(x)dx;  \label{109}
\end{equation}%
where the constant factors $K^{(2)}(\sigma )$ and $(K^{(2)}(\sigma ))^{p}$
are the best possible.
\end{corollary}

Replacing $p>1$ by $0<p<1$ in the above inequalities, we get the equivalent
reverses of (\ref{108}) and (\ref{109}). If there exists a constant $\delta
_{\:0}>0,$ such that for any $\widetilde{\sigma }\in (\sigma -\delta
_{\:0},\sigma ],$%
\[
K^{(2)}(\widetilde{\sigma })=\int_{1}^{\infty }(H(-t)+H(t))t^{\widetilde{%
\sigma }-1}dt\in \mathbf{R}_{+},
\]
then the constant factors in the reverses of (\ref{108}) and (\ref{109}) are
the best possible.

If we have $H(t)=0\ (0<|t|<1)$ in Corollary 12, then $$%
H(v_{1}(x)v_{2}(y))=0\ (0<|v_{1}(x)|<\frac{1}{|v_{2}(y)|}).$$ Setting
\[
\widetilde{E}_{y}:=\{x\in (a_{1},b_{1});\:x\geq v_{1}^{-1}(\frac{1}{|v_{2}(y)|}%
)\text{\,\,or \,\,}x\leq v_{1}^{-1}(\frac{-1}{|v_{2}(y)|})\},
\]
we obtain the following general results:

\begin{corollary}
 Let $p>1,\:\frac{1}{p}+\frac{1}{q}=1,\:H(t)\geq
0,\sigma\in \textbf{R},$%
\begin{eqnarray*}
K^{(2)}(\sigma )=\int_{1}^{\infty }(H(-t)+H(t))t^{\sigma -1}dt\in \mathbf{R}%
_{+},
\end{eqnarray*}
$-\infty \leq a_{i}<b_{i}\leq \infty ,\:v_{i}^{\prime }(s)>0\ (s\in
(a_{i},b_{i})),\:v_{i}(a_{i}^{+})=-\infty ,\:v_{i}(b_{i}^{-})=\infty\ (i=1,2).$
If $f(x),g(y)\geq 0,$ such that%
$$
0<\int_{a_{1}}^{b_{1}}\frac{|v_{1}(x)|^{p(1-\sigma )-1}}{(v_{1}^{\prime
}(x))^{p-1}}f^{p}(x)dx<\infty$$ and $$
0<\int_{a_{2}}^{b_{2}}\frac{|v_{2}(y)|^{q(1-\sigma )-1}}{(v_{2}^{\prime
}(y))^{q-1}}g^{q}(y)dy<\infty,$$
then we have the following equivalent inequalities:%
\begin{eqnarray}
&&\int_{a_{2}}^{b_{2}}\left( \int_{\widetilde{E}%
_{y}}H(v_{1}(x)v_{2}(y))f(x)dx\right) g(y)dy   \nonumber \\
&<&K^{(2)}(\sigma )\left[ \int_{a_{1}}^{b_{1}}\frac{|v_{1}(x)|^{p(1-\sigma
)-1}}{(v_{1}^{\prime }(x))^{p-1}}f^{p}(x)dx\right] ^{\frac{1}{p}}   \nonumber \\
&&\times \left[ \int_{a_{2}}^{b_{2}}\frac{|v_{2}(y)|^{q(1-\sigma )-1}}{%
(v_{2}^{\prime }(y))^{q-1}}g^{q}(y)dy\right] ^{\frac{1}{q}},  \label{110}
\end{eqnarray}%
\begin{eqnarray}
&&\int_{a_{2}}^{b_{2}}\frac{v_{2}^{\prime }(y)}{|v_{2}(y)|^{1-p\sigma }}%
\left( \int_{\widetilde{E}_{y}}H(v_{1}(x)v_{2}(y))f(x)dx\right) ^{p}dy
 \nonumber \\
&<&(K^{(2)}(\sigma ))^{p}\int_{a_{1}}^{b_{1}}\frac{|v_{1}(x)|^{p(1-\sigma
)-1}}{(v_{1}^{\prime }(x))^{p-1}}f^{p}(x)dx,  \label{111}
\end{eqnarray}%
where the constant factors $K^{(2)}(\sigma )$ and $(K^{(2)}(\sigma ))^{p}$
are the best possible.
\end{corollary}

Replacing $p>1$ by $0<p<1$ in the above inequalities, we get the equivalent
reverses of (\ref{110}) and (\ref{111}). If there exists a constant $\delta
_{\:0}>0,$ such that for any $\widetilde{\sigma }\in (\sigma -\delta
_{\:0},\sigma ],$
\[
K^{(2)}(\widetilde{\sigma })=\int_{1}^{\infty }(H(-t)+h(t))t^{\widetilde{%
\sigma }-1}dt\in \mathbf{R}_{+},
\]
then the constant factors in the reverses of (\ref{110}) and (\ref{111}) are
the best possible.

Similarly, if $K_{\lambda }(1,t)=0\ (|t|>1),$ then $$K_{\lambda
}(x,y)=|x|^{-\lambda }K_{\lambda }(1,\frac{y}{x})=0\ (|y|>|x|).$$ By Corollary
14, setting
$$F_{y}:=\left\{x\in \mathbf{R};x\geq |y|\text{\,\, or\,\, }x\leq
-|y|\right\}, $$
we obtain the following Hardy-type integral
inequalities of the first kind, with the homogeneous kernel in the whole plane:

\begin{corollary}
Let $p>1,\frac{1}{p}+\frac{1}{q}%
=1,\:\mu ,\:\sigma \in \mathbf{R},\:\mu +\sigma =\lambda ,$ $K_{\lambda }(x,y)$ is
a homogeneous function of degree $-\lambda $ in $\mathbf{R}^{2},$
\begin{eqnarray*}
K_{\lambda }^{(1)}(\sigma )=\int_{-1}^{1}K_{\lambda }(1,t)|t|^{\sigma
-1}dt\in \mathbf{R}_{+}.
\end{eqnarray*}
If $f(x),\:g(y)\geq 0,$ such that
$$
0<\int_{-\infty }^{\infty }|x|^{p(1-\mu )-1}f^{p}(x)dx<\infty
$$ and $$ 0<\int_{-\infty }^{\infty }|y|^{q(1-\sigma )-1}g^{q}(y)dy<\infty ,
$$
then we have the following equivalent inequalities:%
\begin{eqnarray}
&&\int_{-\infty }^{\infty }\left( \int_{F_{y}}K_{\lambda }(x,y)f(x)dx\right)
g(y)dy   \nonumber \\
&<&K_{\lambda }^{(1)}(\sigma )\left[ \int_{-\infty }^{\infty }|x|^{p(1-\mu
)-1}f^{p}(x)dx\right] ^{\frac{1}{p}}\left[ \int_{-\infty }^{\infty
}|y|^{q(1-\sigma )-1}g^{q}(y)dy\right] ^{\frac{1}{q}},  \label{112}
\end{eqnarray}%
\begin{equation}
\int_{-\infty }^{\infty }|y|^{p\sigma -1}\left( \int_{F_{y}}K_{\lambda
}(x,y)f(x)dx\right) ^{p}dy
<(K_{\lambda }^{(1)}(\sigma ))^{p}\int_{-\infty }^{\infty }|x|^{p(1-\mu
)-1}f^{p}(x)dx,  \label{113}
\end{equation}%
where the constant factors $K_{\lambda }^{(1)}(\sigma )$ and $(K_{\lambda
}^{(1)}(\sigma ))^{p}$ are the best possible.
\end{corollary}

Replacing $p>1$ by $0<p<1$ in the above inequalities, we get the equivalent
reverses of (\ref{112}) and (\ref{113}). If there exists a constant $\delta
_{\:0}>0,$ such that for any $\widetilde{\sigma }\in (\sigma -\delta
_{\:0},\sigma ],$ $$K_{\lambda }^{(1)}(\widetilde{\sigma })=\int_{-1}^{1}K_{%
\lambda }(1,t)|t|^{\widetilde{\sigma }-1}dt\in \mathbf{R}_{+},$$ then the
constant factors in the reverses of (\ref{112}) and (\ref{113}) are the best
possible.

If $K_{\lambda }(1,t)=0\ (|t|>1)$ in Corollary 12, then $$K_{\lambda
}(v_{1}(x),v_{2}(y))=0\ (0<|v_{1}(x)|<|v_{2}(y)|).$$ Setting
\[
\widetilde{F}_{y}:=\{x\in (a_{1},b_{1});\:x\geq v_{1}^{-1}(|v_{2}(y)|)\,\,%
\text{or }\,\,x\leq v_{1}^{-1}(-|v_{2}(y)|)\},
\]
we obtain the following general results:

\begin{corollary}
 Let $p>1,\:\frac{1}{p}+\frac{1}{q}=1,\:\sigma
,\:\mu \in \mathbf{R},\:\sigma +\mu =\lambda ,\:K_{\lambda }(x,y)$ is a
homogeneous function of degree $-\lambda $ in $\mathbf{R}^{2},$
\begin{eqnarray*}
K_{\lambda }^{(1)}(\sigma )=\int_{-1}^{1}K_{\lambda }(1,t)|t|^{\sigma
-1}dt\in \mathbf{R}_{+},
\end{eqnarray*}
$-\infty \leq a_{i}<b_{i}\leq \infty ,\:v_{i}^{\prime }(s)>0$ $(s\in
(a_{i},b_{i})),\:v_{i}(a_{i}^{+})=-\infty ,\:v_{i}(b_{i}^{-})=\infty\ (i=1,2).$
If $f(x),\:g(y)\geq 0,$ such that%
$$
0<\int_{a_{1}}^{b_{1}}\frac{|v_{1}(x)|^{p(1-\mu )-1}}{(v_{1}^{\prime
}(x))^{p-1}}f^{p}(x)dx<\infty$$ and $$
0<\int_{a_{2}}^{b_{2}}\frac{|v_{2}(y)|^{q(1-\sigma )-1}}{(v_{2}^{\prime
}(y))^{q-1}}g^{q}(y)dy<\infty,$$
then we have the following equivalent inequalities:%
\begin{eqnarray}
&&\int_{a_{2}}^{b_{2}}\left( \int_{\widetilde{F}_{y}}K_{\lambda
}(v_{1}(x),v_{2}(y))f(x)dx\right) g(y)dy   \nonumber \\
&<&K_{\lambda }^{(1)}(\sigma )\left[ \int_{a_{1}}^{b_{1}}\frac{%
|v_{1}(x)|^{p(1-\mu )-1}}{(v_{1}^{\prime }(x))^{p-1}}f^{p}(x)dx\right] ^{%
\frac{1}{p}}   \nonumber \\
&&\times \left[ \int_{a_{2}}^{b_{2}}\frac{|v_{2}(y)|^{q(1-\sigma )-1}}{%
(v_{2}^{\prime }(y))^{q-1}}g^{q}(y)dy\right] ^{\frac{1}{q}},  \label{114}
\end{eqnarray}%
\begin{eqnarray}
&&\int_{a_{2}}^{b_{2}}\frac{v_{2}^{\prime }(y)}{|v_{2}(y)|^{1-p\sigma }}%
\left( \int_{\widetilde{F}_{y}}K_{\lambda }(v_{1}(x),v_{2}(y))f(x)dx\right)
^{p}dy   \nonumber \\
&<&(K_{\lambda }^{(1)}(\sigma ))^{p}\int_{a_{1}}^{b_{1}}\frac{%
|v_{1}(x)|^{p(1-\mu )-1}}{(v_{1}^{\prime }(x))^{p-1}}f^{p}(x)dx,  \label{115}
\end{eqnarray}%
where the constant factors $K_{\lambda }^{(1)}(\sigma )$ and $(K_{\lambda
}^{(1)}(\sigma ))^{p}$ are the best possible.
\end{corollary}

Replacing $p>1$ by $0<p<1$ in the above inequalities, we obtain the equivalent
reverses of (\ref{114}) and (\ref{115}). If there exists a constant $\delta
_{\:0}>0,$ such that for any $\widetilde{\sigma }\in (\sigma -\delta
_{\:0},\sigma ],$%
\[
K_{\lambda }^{(1)}(\widetilde{\sigma })=\int_{-1}^{1}K_{\lambda }(1,t)|t|^{%
\widetilde{\sigma }-1}dt\in \mathbf{R}_{+},
\]
then the constant factors in the reverses of (\ref{114}) and (\ref{115}) are
the best possible.

Similarly, if $K_{\lambda }(1,t)=0\ (0<|t|<1)$ in Corollary 14, then $$%
K_{\lambda }(x,y)=0\ (|x|>|y|>0).$$
The following Hardy-type integral inequalities of the second kind with the homogeneous kernel, hold true:

\begin{corollary}
Let $p>1,\:\frac{1}{p}+\frac{1}{q}=1,\:\mu
,\:\sigma \in \mathbf{R},\:\mu +\sigma =\lambda ,$ $K_{\lambda }(x,y)$ is a
homogeneous function of degree $-\lambda $ in $\mathbf{R}^{2},$
$$
K_{\lambda }^{(2)}(\sigma )=\int_{1}^{\infty }(K_{\lambda }(1,-t)+K_{\lambda
}(1,t))t^{\sigma -1}dt\in \mathbf{R}_{+}.
$$
If $f(x),g(y)\geq 0,$%
$$
0<\int_{-\infty }^{\infty }|x|^{p(1-\mu )-1}f^{p}(x)dx<\infty
$$ and $$ 0<\int_{-\infty }^{\infty }|y|^{q(1-\sigma )-1}g^{q}(y)dy<\infty,$$
then we have the following equivalent inequalities:%
$$
\int_{-\infty }^{\infty }\left( \int_{-|y|}^{|y|}K_{\lambda
}(x,y)f(x)dx\right)g(y)dy
=\int_{-\infty }^{\infty }\left( \int_{-|x|}^{|x|}K_{\lambda
}(x,y)g(y)dy\right) f(x)dx
$$
\begin{equation}
<K_{\lambda }^{(2)}(\sigma )\left[ \int_{-\infty }^{\infty }|x|^{p(1-\mu
)-1}f^{p}(x)dx\right] ^{\frac{1}{p}}\left[ \int_{-\infty }^{\infty
}|y|^{q(1-\sigma )-1}g^{q}(y)dy\right] ^{\frac{1}{q}},  \label{116}
\end{equation}%
\begin{equation}
\int_{-\infty }^{\infty }|y|^{p\sigma -1}\left(
\int_{-|y|}^{|y|}K_{\lambda }(x,y)f(x)dx\right) ^{p}dy
<(K_{\lambda }^{(2)}(\sigma ))^{p}\int_{-\infty }^{\infty }|x|^{p(1-\mu
)-1}f^{p}(x)dx;  \label{117}
\end{equation}%
where the constant factors $K_{\lambda }^{(2)}(\sigma )$ and $(K_{\lambda
}^{(2)}(\sigma ))^{p}$ are the best possible.
\end{corollary}

Replacing $p>1$ by $0<p<1$ in the above inequalities, we obtain the equivalent
reverses of (\ref{116}) and (\ref{117}). If there exists a constant $\delta
_{\:0}>0,$ such that for any $\widetilde{\sigma }\in (\sigma -\delta
_{\:0},\sigma ],$%
\[
K_{\lambda }^{(2)}(\widetilde{\sigma })=\int_{1}^{\infty }(K_{\lambda
}(1,-t)+K_{\lambda }(1,t))t^{\widetilde{\sigma }-1}dt\in \mathbf{R}_{+},
\]
then the constant factors in the reverses of (\ref{116}) and (\ref{117}) are
the best possible.

If $K_{\lambda }(1,t)=0\ (0<|t|<1)$ in Corollary 12, then $$K_{\lambda
}(v_{1}(x),v_{2}(y))=0\ (|v_{1}(x)|>|v_{2}(y)|>0),$$ we have the following
general results:

\begin{corollary}
Let $p>1,\:\frac{1}{p}+\frac{1}{q}=1,\:\mu
,\:\sigma \in \mathbf{R},\:\mu +\sigma =\lambda ,\:K_{\lambda }(x,y)$ is a
homogeneous function of degree $-\lambda $ in $\mathbf{R}^{2},$
$$
K_{\lambda }^{(2)}(\sigma )=\int_{1}^{\infty }(K_{\lambda }(1,-t)+K_{\lambda
}(1,t))t^{\sigma -1}dt\in \mathbf{R}_{+},
$$
$-\infty \leq a_{i}<b_{i}\leq \infty ,\:v_{i}^{\prime }(s)>0$ $(s\in
(a_{i},b_{i})),\:v_{i}(a_{i}^{+})=-\infty ,\:v_{i}(b_{i}^{-})=\infty\ (i=1,2).$
If $f(x),\:g(y)\geq 0,$ such that%
$$
0<\int_{a_{1}}^{b_{1}}\frac{|v_{1}(x)|^{p(1-\mu )-1}}{(v_{1}^{\prime
}(x))^{p-1}}f^{p}(x)dx<\infty$$ and 
$$
0 <\int_{a_{2}}^{b_{2}}\frac{|v_{2}(y)|^{q(1-\sigma )-1}}{(v_{2}^{\prime
}(y))^{q-1}}g^{q}(y)dy<\infty,
$$
then we have the following equivalent inequalities:%
\begin{eqnarray}
&&\int_{a_{2}}^{b_{2}}\left(
\int_{v_{1}^{-1}(-|v_{2}(y)|)}^{v_{1}^{-1}(|v_{2}(y)|)}K_{\lambda
}(v_{1}(x),v_{2}(y))f(x)dx\right) g(y)dy   \nonumber \\
&<&K_{\lambda }^{(2)}(\sigma )\left[ \int_{a_{1}}^{b_{1}}\frac{%
|v_{1}(x)|^{p(1-\mu )-1}}{(v_{1}^{\prime }(x))^{p-1}}f^{p}(x)dx\right] ^{%
\frac{1}{p}}   \nonumber \\
&&\times \left[ \int_{a_{2}}^{b_{2}}\frac{|v_{2}(y)|^{q(1-\sigma )-1}}{%
(v_{2}^{\prime }(y))^{q-1}}g^{q}(y)dy\right] ^{\frac{1}{q}},  \label{118}
\end{eqnarray}%
\begin{eqnarray}
&&\int_{a_{2}}^{b_{2}}\frac{v_{2}^{\prime }(y)}{|v_{2}(y)|^{1-p\sigma }}%
\left( \int_{v_{1}^{-1}(-|v_{2}(y)|)}^{v_{1}^{-1}(|v_{2}(y)|)}K_{\lambda
}(v_{1}(x),v_{2}(y))f(x)dx\right) ^{p}dy   \nonumber \\
&<&(K_{\lambda }^{(2)}(\sigma ))^{p}\int_{a_{1}}^{b_{1}}\frac{%
|v_{1}(x)|^{p(1-\mu )-1}}{(v_{1}^{\prime }(x))^{p-1}}f^{p}(x)dx,  \label{119}
\end{eqnarray}%
where the constant factors $K_{\lambda }^{(2)}(\sigma )$ and $(K_{\lambda
}^{(2)}(\sigma ))^{p}$ are the best possible.
\end{corollary}

Replacing $p>1$ by $0<p<1$ in the above inequalities, we have the equivalent
reverses of (\ref{118}) and (\ref{119}). If there exists a constant $\delta
_{\:0}>0,$ such that for any $\widetilde{\sigma }\in (\sigma -\delta
_{\:0},\sigma ],$%
\[
K_{\lambda }^{(2)}(\widetilde{\sigma })=\int_{1}^{\infty }(K_{\lambda
}(1,-t)+K_{\lambda }(1,t))t^{\widetilde{\sigma }-1}dt\in \mathbf{R}_{+},
\]
then the constant factors in the reverses of (\ref{118}) and (\ref{119}) are
the best possible.
\subsection{Yang-Hilbert-Type Operators and Hardy-Type Operators in the
Whole Plane}

Let $p>1,\:\frac{1}{p}+\frac{1}{q}=1,\:\mu ,\:\sigma \in \mathbf{R},\:\mu
+\sigma =\lambda .$ We define the following functions:
\[
\varphi (x):=|x|^{p(1-\sigma)-1},\psi (y):=|y|^{q(1-\sigma )-1},\phi
(x):=|x|^{p(1-\mu )-1}(x,y\in \mathbf{R}),
\]
wherefrom, $\psi ^{1-p}(y)=|y|^{p\sigma-1}.$ 

We define also the following real
normed linear space:
\[
L_{p,\varphi }(\mathbf{R}):=\left\{f:||f||_{p,\varphi }:=\left\{\int_{-\infty
}^{\infty }\varphi (x)|f(x)|^{p}dx\right\}^{\frac{1}{p}}<\infty\right\},
\]
wherefrom,%
\begin{eqnarray*}
L_{p,\psi ^{1-p}}(\mathbf{R}) &=&\left\{h:||h||_{p,\psi ^{1-p}}:=\left\{\int_{-\infty
}^{\infty }\psi ^{1-p}(y)|h(y)|^{p}dy\right\}^{\frac{1}{p}}<\infty\right\}, \\
L_{p,\phi}(\mathbf{R}) &=&\{g:||g||_{p,\phi}:=\left\{\int_{-\infty }^{\infty
}\phi(x)|g(x)|^{p}dx\}^{\frac{1}{p}}<\infty\right\}.
\end{eqnarray*}

(a) In view of Theorem 5 ($\delta =1)$, for $f\in L_{p,\varphi }(\mathbf{R}%
),$ $$H_{1}(y):=\int_{-\infty }^{\infty }H(xy)|f(x)|dx^{{}}\ \ (y\in \mathbf{R}%
_{+}), $$ by (\ref{91}), we have
\begin{equation}
||H_{1}||_{p,\psi ^{1-p}}:=\left( \int_{-\infty }^{\infty }\psi
^{1-p}(y)H_{1}^{p}(y)dy\right) ^{\frac{1}{p}}<K(\sigma )||f||_{p,\varphi
}<\infty .  \label{120}
\end{equation}
\begin{definition}
 Define Yang-Hilbert-type integral operator with the
non-homogeneous kernel in the whole plane $$T_{1}\::\:L_{p,\varphi }(\mathbf{R})\rightarrow L_{p,\psi ^{1-p}}(\mathbf{R})$$ as follows:
 
 For any $f\in L_{p,\varphi }(\mathbf{R}),$ there exists a unique representation $$T_{1}f=H_{1}\in L_{p,\psi ^{1-p}}(\mathbf{R}),$$ satisfying
$$T_{1}f(y)=H_{1}(y),$$
for any $y\in
\mathbf{R}$.
\end{definition}

In view of (\ref{120}), it follows that
\[
||T_{1}f||_{p,\psi ^{1-p}}=||H_{1}||_{p,\psi ^{1-p}}\leq K(\sigma
)||f||_{p,\varphi }.
\]
Therefore, the operator $T_{1}$ is bounded and it satisfies the following relation%
\[
||T_{1}||=\sup_{f(\neq \theta )\in L_{p,\varphi }(\mathbf{R})}\frac{%
||T_{1}f||_{p,\psi ^{1-p}}}{||f||_{p,\varphi }}\leq K(\sigma ).
\]
Since the constant factor $K(\sigma )$ in (\ref{120}) is the best possible,
we have%
\begin{equation}
||T_{1}||=K(\sigma )=\int_{-\infty }^{\infty }H(t)|t|^{\sigma -1}dt.
\label{121}
\end{equation}

If we define the formal inner product of $T_{1}f$ and $g$ as
\begin{eqnarray*}
(T_{1}f,g) &:&=\int_{-\infty }^{\infty }\left( \int_{-\infty }^{\infty
}H(xy)f(x)dx\right) g(y)dy \\
&=&\int_{-\infty }^{\infty }\int_{-\infty }^{\infty }H(xy)f(x)g(y)dxdy,
\end{eqnarray*}
then we can rewrite (\ref{90}) and (\ref{91}) as follows:%
\[
(T_{1}f,g)<||T_{1}||\cdot ||f||_{p,\varphi }||g||_{q,\psi
},\ \ ||T_{1}f||_{p,\psi ^{1-p}}<||T_{1}||\cdot ||f||_{p,\varphi }.
\]

(b) In view of Corollary 15, for $f\in L_{p,\varphi }(\mathbf{R}),$ setting%
\[
H_{1}^{(1)}(y):=\int_{\frac{-1}{|y|}}^{\frac{1}{|y|}}H(xy)|f(x)|dx^{{}}(y\in
\mathbf{R}\backslash \{0\}),
\]
by (\ref{104}), we obtain
\begin{equation}
||H_{1}^{(1)}||_{p,\psi ^{1-p}} =\left( \int_{-\infty }^{\infty }\psi
^{1-p}(y)(H_{1}^{(1)}(y))^{p}dy\right) ^{\frac{1}{p}} <K^{(1)}(\sigma
)||f||_{p,\varphi }<\infty .  \label{122}
\end{equation}

\begin{definition}
 Let us define the Hardy-type integral operator of the first kind,
with the non-homogeneous kernel in the whole plane $$T_{1}^{(1)}\::\:%
L_{p,\varphi }(\mathbf{R})\rightarrow L_{p,\psi ^{1-p}}(\mathbf{R})$$ as
follows:

For any $f\in L_{p,\varphi }(\mathbf{R}),$ there exists a unique
representation $$T_{1}^{(1)}f=H_{1}^{(1)}\in L_{p,\psi ^{1-p}}(\mathbf{R}),$$
satisfying $$T_{1}^{(1)}f(y)=H_{1}^{(1)}(y),$$
for any $y\in \mathbf{R}.$
\end{definition}

In view of (\ref{122}), it follows that
\[
|T_{1}^{(1)}f||_{p,\psi ^{1-p}}=||H_{1}^{(1)}||_{p,\psi ^{1-p}}\leq
K^{(1)}(\sigma )||f||_{p,\varphi }.
\]
Then, the operator $T_{1}^{(1)}$ is bounded satisfying%
\[
||T_{1}^{(1)}||=\sup_{f(\neq \theta )\in L_{p,\varphi }(\mathbf{R})}\frac{%
||T_{1}^{(1)}f||_{p,\psi ^{1-p}}}{||f||_{p,\varphi }}\leq K^{(1)}(\sigma ).
\]
Since the constant factor $K^{(1)}(\sigma )$ in (\ref{122}) is the best
possible, we have%
\begin{equation}
||T_{1}^{(1)}||=K^{(1)}(\sigma )=\int_{-1}^{1}H(t)|t|^{\sigma -1}dt.
\label{123}
\end{equation}

Setting the formal inner product of $T_{1}^{(1)}f$ and $g$ as
\[
(T_{1}^{(1)}f,g)=\int_{-\infty }^{\infty }\left( \int_{\frac{-1}{|y|}}^{%
\frac{1}{|y|}}H(xy)f(x)dx\right) g(y)dy,
\]
we can rewrite (\ref{103}) and (\ref{104}) as follows:%
\begin{equation}
(T_{1}^{(1)}f,g)<||T_{1}^{(1)}||\cdot ||f||_{p,\varphi }||g||_{q,\psi
},\ \ ||T_{1}^{(1)}f||_{p,\psi ^{1-p}}<||T_{1}^{(1)}||\cdot ||f||_{p,\varphi }.
\label{124}
\end{equation}

(c) In view of Corollary 17, for $f\in L_{p,\varphi }(\mathbf{R}),$ setting $$%
H_{1}^{(2)}(y):=\int_{E_{y}}H(xy)|f(x)|dx\ \ (y\in \mathbf{R}), $$ by (\ref%
{109}), we have
\begin{eqnarray}
||H_{1}^{(2)}||_{p,\psi ^{1-p}} &=&\left( \int_{-\infty }^{\infty }\psi
^{1-p}(y)(H_{1}^{(2)}(y))^{p}dy\right) ^{\frac{1}{p}} \nonumber\\
&<&K^{(2)}(\sigma
)||f||_{p,\varphi }<\infty .  \label{125}
\end{eqnarray}

\begin{definition}
Let us define the Hardy-type integral operator of the second kind
with the non-homogeneous kernel in the whole plane $$T_{1}^{(2)}\::\:%
L_{p,\varphi }(\mathbf{R})\rightarrow L_{p,\psi ^{1-p}}(\mathbf{R})$$ as
follows:

 For any $f\in L_{p,\varphi }(\mathbf{R}),$ there exists a unique
representation $$T_{1}^{(2)}f=H_{1}^{(2)}\in L_{p,\psi ^{1-p}}(\mathbf{R}),$$
satisfying $$T_{1}^{(2)}f(y)=H_{1}^{(2)}(y),$$
for any $y\in \mathbf{R}.$
\end{definition}

In view of (\ref{125}), it follows that
\[
||T_{1}^{(2)}f||_{p,\psi ^{1-p}}=||H_{1}^{(2)}||_{p,\psi ^{1-p}}\leq
K^{(2)}(\sigma )||f||_{p,\varphi }.
\]
Thus, the operator $T_{1}^{(2)}$ is bounded satisfying%
\[
||T_{1}^{(2)}||=\sup_{f(\neq \theta )\in L_{p,\varphi }(\mathbf{R})}\frac{%
||T_{1}^{(2)}f||_{p,\psi ^{1-p}}}{||f||_{p,\varphi }}\leq K^{(2)}(\sigma ).
\]
Since the constant factor $K^{(2)}(\sigma )$ in (\ref{125}) is the best
possible, we have%
\begin{equation}
||T_{1}^{(2)}||=K^{(2)}(\sigma )=\int_{1}^{\infty }(H(-t)+H(t))t^{\sigma
-1}dt.  \label{126}
\end{equation}

Setting the formal inner product of $T_{1}^{(2)}f$ and $g$ as
\[
(T_{1}^{(2)}f,g)=\int_{-\infty }^{\infty }\left(
\int_{E_{y}}H(xy)f(x)dx\right) g(y)dy,
\]
we can rewrite (\ref{108}) and (\ref{109}) as follows:%
\begin{equation}
(T_{1}^{(2)}f,g)<||T_{1}^{(2)}||\cdot ||f||_{p,\varphi }||g||_{q,\psi
},\ \ ||T_{1}^{(2)}f||_{p,\psi ^{1-p}}<||T_{1}^{(2)}||\cdot ||f||_{p,\varphi }.
\label{127}
\end{equation}

(d) In view of Corollary 14, for $f\in L_{p,\phi }(\mathbf{R}),$  $$%
H_{2}(y):=\int_{-\infty }^{\infty }K_{\lambda }(x,y)|f(x)|dx\ \ (y\in \mathbf{%
R}), $$ by (\ref{101}), we have
\begin{eqnarray}
||H_{2}||_{p,\psi ^{1-p}}&=&\left( \int_{-\infty }^{\infty }\psi
^{1-p}(y)H_{2}^{p}(y)dy\right) ^{\frac{1}{p}}\nonumber\\
&<&K_{\lambda }(\sigma
)||f||_{p,\phi }<\infty .  \label{128}
\end{eqnarray}

\begin{definition}
 We define the Yang-Hilbert-type integral operator with the
homogeneous kernel in the whole plane $$T_{2}\::\:L_{p,\phi }(\mathbf{R}%
)\rightarrow L_{p,\psi ^{1-p}}(\mathbf{R})$$ as follows:

For any $f\in
L_{p,\phi }(\mathbf{R}),$ there exists a unique representation $$%
T_{2}f=H_{2}\in L_{p,\psi ^{1-p}}(\mathbf{R}),$$ satisfying $$T_{2}f(y)=H_{2}(y),$$
for any $y\in
\mathbf{R}.$
\end{definition}

By (\ref{128}), it follows that
\[
||T_{2}f||_{p,\psi^{1-p} }=||H_{2}||_{p,\psi ^{1-p}}\leq K_{\lambda }(\sigma
)||f||_{p,\phi}.
\]
Hence, the operator $T_{2}$ is bounded satisfying%
\[
||T_{2}||=\sup_{f(\neq \theta )\in L_{p,\phi }(\mathbf{R})}\frac{%
||T_{2}f||_{p,\psi^{1-p}}}{||f||_{p,\phi }}\leq K_{\lambda }(\sigma ).
\]
Since the constant factor $K_{\lambda }(\sigma )$ in (\ref{128}) is the best
possible, we have%
\begin{equation}
||T_{2}||=K_{\lambda }(\sigma )=\int_{-\infty }^{\infty }K_{\lambda
}(1,t)|t|^{\sigma -1}dt.  \label{129}
\end{equation}

Setting the formal inner product of $T_{2}f$ and $g$ as
\begin{eqnarray*}
(T_{2}f,g) &=&\int_{-\infty }^{\infty }\left( \int_{-\infty }^{\infty
}K_{\lambda }(x,y)f(x)dx\right) g(y)dy \\
&=&\int_{-\infty }^{\infty }\int_{-\infty }^{\infty }K_{\lambda
}(x,y)f(x)g(y)dxdy,
\end{eqnarray*}%
we can rewrite (\ref{100}) and (\ref{101}) as follows:%
\[
(T_{2}f,g)<||T_{2}||\cdot ||f||_{p,\phi }||g||_{q,\psi },\ \ ||T_{2}f||_{p,\psi
^{1-p}}<||T_{2}||\cdot ||f||_{p,\phi }.
\]

(e) By Corollary 19, for $f\in L_{p,\phi }(\mathbf{R}),$  $$%
H_{2}^{(1)}(y):=\int_{F_{y}}K_{\lambda }(x,y)|f(x)|dx\ \ (y\in \mathbf{R}), $$
combined with (\ref{113}), we obtain
\begin{eqnarray}
||H_{2}^{(1)}||_{p,\psi ^{1-p}} &=&\left( \int_{-\infty }^{\infty }\psi
^{1-p}(y)(H_{2}^{(1)}(y))^{p}dy\right) ^{\frac{1}{p}} \nonumber\\
&<&K_{\lambda
}^{(1)}(\sigma )||f||_{p,\phi }<\infty .  \label{130}
\end{eqnarray}

\begin{definition}
 We define the Hardy-type integral operator of the first kind,
with the homogeneous kernel in the whole plane $$T_{2}^{(1)}\::\:L_{p,\phi
}(\mathbf{R})\rightarrow L_{p,\psi ^{1-p}}(\mathbf{R})$$ as follows:

For any $%
f\in L_{p,\phi }(\mathbf{R}),$ there exists a unique representation $$%
T_{2}^{(1)}f=H_{2}^{(1)}\in L_{p,\psi^{1-p}}(\mathbf{R}),$$ satisfying $$T_{2}^{(1)}f(y)=H_{2}^{(1)}(y),$$
for
any $y\in \mathbf{R}.$
\end{definition}

In view of (\ref{130}), it follows that
\[
||T_{2}^{(1)}f||_{p,\psi ^{1-p}}=||H_{2}^{(1)}||_{p,\psi ^{1-p}}\leq
K_{\lambda }^{(1)}(\sigma )||f||_{p,\phi }.
\]
Therefore, the operator $T_{2}^{(1)}$ is bounded satisfying%
\[
||T_{2}^{(1)}||=\sup_{f(\neq \theta )\in L_{p,\phi }(\mathbf{R})}\frac{%
||T_{2}^{(1)}f||_{p,\psi ^{1-p}}}{||f||_{p,\phi}}\leq K_{\lambda
}^{(1)}(\sigma ).
\]
Since the constant factor $K_{\lambda }^{(1)}(\sigma )$ in (\ref{130}) is
the best possible, we have%
\begin{equation}
||T_{2}^{(1)}||=K_{\lambda }^{(1)}(\sigma )=\int_{-1}^{1}K_{\lambda
}(1,t)|t|^{\sigma -1}dt.  \label{131}
\end{equation}

Setting the formal inner product of $T_{2}^{(1)}f$ and $g$ as
\[
(T_{2}^{(1)}f,g)=\int_{-\infty }^{\infty }\left( \int_{F_{y}}K_{\lambda
}(x,y)f(x)dx\right) g(y)dy,
\]
we can rewrite (\ref{112}) and (\ref{113}) as follows:%
\[
(T_{2}^{(1)}f,g)<||T_{2}^{(1)}||\cdot ||f||_{p,\phi}||g||_{q,\psi
},\ \ ||T_{2}^{(1)}f||_{p,\psi ^{1-p}}<||T_{2}^{(1)}||\cdot ||f||_{p,\phi }.
\]

\bigskip (f) In view of Corollary 21, for $f\in L_{p,\phi }(\mathbf{R}_{+}),$
$$H_{2}^{(2)}(y):=\int_{-|y|}^{|y|}K_{\lambda }(x,y)|f(x)|dx\ (y\in \mathbf{%
R}), $$ by (\ref{117}), we have
\begin{eqnarray}
||H_{2}^{(2)}||_{p,\psi ^{1-p}} &:&=\left( \int_{-\infty }^{\infty }\psi
^{1-p}(y)(H_{2}^{(2)}(y))^{p}dy\right) ^{\frac{1}{p}} \nonumber\\
&<&K_{\lambda
}^{(2)}(\sigma )||f||_{p,\phi }<\infty .  \label{132}
\end{eqnarray}

\begin{definition}
We define the Hardy-type integral operator of the second kind,
with the homogeneous kernel in the whole plane $$T_{2}^{(2)}\::\:L_{p,\phi
}(\mathbf{R})\rightarrow L_{p,\psi^{1-p}}(\mathbf{R})$$ as follows:

 For any $%
f\in L_{p,\phi }(\mathbf{R}),$ there exists a unique representation $$%
T_{2}^{(2)}f=H_{2}^{(2)}\in L_{p,\psi^{1-p}}(\mathbf{R}),$$ satisfying $$T_{2}^{(2)}f(y)=H_{2}^{(2)}(y),$$
for
any $y\in \mathbf{R}.$
\end{definition}

By (\ref{132}), it follows that
\[
||T_{2}^{(2)}f||_{p,\psi ^{1-p}}=||H_{2}^{(2)}||_{p,\psi ^{1-p}}\leq
K_{\lambda }^{(2)}(\sigma )||f||_{p,\phi}
\]
and then the operator $T_{2}^{(2)}$ is bounded satisfying%
\[
||T_{2}^{(2)}||=\sup_{f(\neq \theta )\in L_{p,\phi }(\mathbf{R})}\frac{%
||T_{2}^{(2)}f||_{p,\psi ^{1-p}}}{||f||_{p,\phi}}\leq K_{\lambda
}^{(2)}(\sigma ).
\]
Since the constant factor $K_{\lambda }^{(2)}(\sigma )$ in (\ref{132}) is
the best possible, we have%
\begin{equation}
||T_{2}^{(2)}||=K_{\lambda }^{(2)}(\sigma )=\int_{1}^{\infty }(K_{\lambda
}(1,-t)+K_{\lambda }(1,t))t^{\sigma -1}dt.  \label{133}
\end{equation}

Setting the formal inner product of $T_{2}^{(2)}f$ and $g$ as
\[
(T_{2}^{(2)}f,g)=\int_{-\infty }^{\infty }\left( \int_{-|y|}^{|y|}K_{\lambda
}(x,y)f(x)dx\right) g(y)dy,
\]
we can rewrite (\ref{116}) and (\ref{117}) as follows:%
\[
(T_{2}^{(2)}f,g)<||T_{2}^{(2)}||\cdot ||f||_{p,\phi }||g||_{q,\psi
},\ \ ||T_{2}^{(2)}f||_{p,\psi ^{1-p}}<||T_{2}^{(2)}||\cdot ||f||_{p,\phi }.
\]

\begin{remark}$ $
 (a) If $K_{\lambda }(x,y)$ is a symmetric function
satisfying $K_{\lambda }(y,x)=K_{\lambda }(x,y),$ then by setting
$$H(t)=:K_{\lambda
}(1,t)\arctan |t|^{\beta }\ (\beta \in \mathbf{R}),$$
and $\mu =\sigma =\frac{%
\lambda }{2}$ in (\ref{121}), we obtain%
\begin{eqnarray}
||T_{1}|| &=&\int_{-\infty }^{\infty }H(t)|t|^{\sigma -1}=\int_{-\infty
}^{\infty }K_{\lambda }(1,t)\arctan |t|^{\beta }|t|^{\sigma -1}dt\nonumber\\
&=&\frac{\pi }{4}K_{\lambda }\left(\frac{\lambda }{2}\right),  \label{134}
\end{eqnarray}%
where $$K_{\lambda }\left(\frac{\lambda }{2}\right)=\int_{-\infty }^{\infty }K_{\lambda
}(1,t)|t|^{\sigma -1}dt.$$ In fact, we obtain%
\begin{eqnarray*}
&&\int_{0}^{\infty }K_{\lambda }(1,t)(\arctan t^{\beta })t^{\frac{\lambda }{2%
}-1}dt \\
&=&\int_{0}^{1}K_{\lambda }(1,t)(\arctan t^{\beta })t^{\frac{\lambda }{2}%
-1}dt
+\int_{1}^{\infty }K_{\lambda }(1,u)(\arctan u^{\beta })u^{\frac{\lambda }{%
2}-1}du\\
&=&\int_{0}^{1}K_{\lambda }(1,t)(\arctan t^{\beta })t^{\frac{\lambda }{2}%
-1}dt
+\int_{0}^{1}K_{\lambda }(t,1)(\arctan t^{-\beta })t^{\frac{\lambda }{2}%
-1}dt\\
&=&\int_{0}^{1}K_{\lambda }(1,t)(\arctan t^{\beta }+\arctan t^{-\beta })t^{%
\frac{\lambda }{2}-1}dt \\
&=&\frac{\pi }{2}\int_{0}^{1}K_{\lambda }(1,t)t^{\frac{\lambda }{2}-1}dt=%
\frac{\pi }{4}\int_{0}^{\infty }K_{\lambda }(1,t)t^{\frac{\lambda }{2}-1}dt.
\end{eqnarray*}%
Similarly, we get%
\begin{eqnarray*}
&&\int_{-\infty }^{0}K_{\lambda }(1,t)\arctan (-t)^{\beta }(-t)^{\sigma -1}dt
\\
&=&\int_{0}^{\infty }K_{\lambda }(1,-u)(\arctan u^{\beta })u^{\sigma -1}du
=\frac{\pi }{4}\int_{0}^{\infty }K_{\lambda }(1,-t)t^{\frac{\lambda }{2}%
-1}dt.
\end{eqnarray*}%
Then we have
$$
||T_{1}||=\frac{\pi }{4}\int_{0}^{\infty }(K_{\lambda }(1,-t)+K_{\lambda
}(1,t))t^{\frac{\lambda }{2}-1}dt=\frac{\pi }{4}K_{\lambda }\left(\frac{\lambda }{%
2}\right).
$$
\end{remark}

(b) If we replace $H(t)$ by $$h(|t|^{\gamma }+t^{\gamma }\cos \alpha )(\gamma \in
\{b;b=\frac{1}{2k-1},2k+1\ (k\in \mathbf{N})\},\alpha \in (0,\pi ))$$ in (\ref%
{121}), where, $h(t)$ is a non-negative measurable function in $\mathbf{R}%
_{+},$ satisfying $$k\left(\frac{\sigma }{\gamma }\right)=\int_{0}^{\infty }h(t)t^{\frac{%
\sigma }{\gamma }-1}dt\in \mathbf{R}_{+}, $$ it follows that%
\begin{eqnarray}
||T_{1}|| &=&\int_{-\infty }^{\infty }h(|t|^{\gamma }+t^{\gamma }\cos \alpha
)|t|^{\sigma -1}  \nonumber \\
&=&\frac{1}{\gamma 2^{\sigma /\gamma }}[(\sec \frac{\alpha }{2})^{\frac{%
2\alpha }{\gamma }}+(\csc \frac{\alpha }{2})^{\frac{2\alpha }{\gamma }}]k\left(%
\frac{\sigma }{\gamma }\right),  \label{135}
\end{eqnarray}%
In particular, setting $h(t)=k_{\lambda }(1,t)\ (t\in \mathbf{R}_{+}),$ it
follows that
\begin{eqnarray}
||T_{1}|| &=&\int_{-\infty }^{\infty }k_{\lambda }(1,|t|^{\gamma }+t^{\gamma
}\cos \alpha )|t|^{\sigma -1}  \nonumber \\
&=&\frac{1}{\gamma 2^{\sigma /\gamma }}[(\sec \frac{\alpha }{2})^{\frac{%
2\alpha }{\gamma }}+(\csc \frac{\alpha }{2})^{\frac{2\alpha }{\gamma }%
}]k_{\lambda }(\frac{\sigma }{\gamma }),  \label{136}
\end{eqnarray}%
where $$k_{\lambda }\left(\frac{\sigma }{\gamma }\right)=\int_{0}^{\infty }k_{\lambda
}(1,t)t^{\frac{\sigma }{\gamma }-1}dt\in \mathbf{R}_{+}.$$

In fact, setting $u=t^{\gamma }(1+\cos \alpha )$, we get
\begin{eqnarray*}
&&\int_{0}^{\infty }h(|t|^{\gamma }+t^{\gamma }\cos \alpha )|t|^{\sigma -1}dt
\\
&=&\int_{0}^{\infty }h(t^{\gamma }(1+\cos \alpha ))t^{\sigma -1}dt=\frac{1}{%
\gamma 2^{\sigma /\gamma }}(\sec \frac{\alpha }{2})^{\frac{2\alpha }{\gamma }%
}k(\frac{\sigma }{\gamma });
\end{eqnarray*}%
Moreover, setting $u=-t,$ we have%
\begin{eqnarray*}
&&\int_{-\infty }^{0}h(|t|^{\gamma }+t^{\gamma }\cos \alpha )|t|^{\sigma
-1}dt =\int_{-\infty }^{0}h(-t^{\gamma }(1-\cos \alpha ))(-t)^{\sigma -1}dt
\\
&=&\int_{0}^{\infty }h(u^{\gamma }(1-\cos \alpha ))u^{\sigma -1}du=\frac{1}{%
\gamma 2^{\sigma /\gamma }}(\csc \frac{\alpha }{2})^{\frac{2\alpha }{\gamma }%
}k(\frac{\sigma }{\gamma }),
\end{eqnarray*}%
and then (\ref{135}) follows.

(c) Replacing $K_{\lambda }(1,t)$ by $k_{\lambda }(1,|t|^{\gamma }+t^{\gamma
}\cos \alpha )$ in (\ref{134}), where, $k_{\lambda }(x,y)$ is a homogeneous
function in $\mathbf{R}_{+},$ satisfying $$k_{\lambda }\left(\frac{\sigma }{\gamma
}\right)=\int_{0}^{\infty }k_{\lambda }(1,t)t^{\frac{\sigma }{\gamma }-1}dt\in
\mathbf{R}_{+}, $$ in view of (\ref{136}), we obtain
\begin{eqnarray}
||T_{1}|| &=&\int_{-\infty }^{\infty }k_{\lambda }(1,|t|^{\gamma }+t^{\gamma
}\cos \alpha )\arctan |t|^{\beta }|t|^{\frac{\lambda }{2}-1}dt  \nonumber \\
&=&\frac{\pi }{\gamma 2^{2+\sigma /\gamma }}\left[\left(\sec \frac{\alpha }{2}\right)^{\frac{%
2\alpha }{\gamma }}+\left(\csc \frac{\alpha }{2}\right)^{\frac{2\alpha }{\gamma }%
}\right]k_{\lambda }\left(\frac{\lambda }{2\gamma }\right).  \label{137}
\end{eqnarray}

\subsection{ Some Examples}

\begin{example}$ $
 (a) Set $$H(t)=K_{\lambda }(1,t)=\frac{1}{%
|1+t|^{\lambda }}\ (\mu ,\sigma >0,\mu +\sigma =\lambda <1).$$ Then we have the
kernels $$
H(xy)=\frac{1}{|1+xy|^{\lambda }},\ \ K_{\lambda }(x,y)=\frac{1}{|x+y|^{\lambda }%
} $$
and obtain the constant factors
\begin{eqnarray*}
K(\sigma ) &=&K_{\lambda }(\sigma )=\int_{-\infty }^{\infty }\frac{%
|t|^{\sigma -1}}{|1+t|^{\lambda }}dt \\
&=&\int_{0}^{\infty }\frac{t^{\sigma -1}}{|1-t|^{\lambda }}%
dt+\int_{0}^{\infty }\frac{t^{\sigma -1}}{(1+t)^{\lambda }}dt \\
&=&B(1-\lambda ,\sigma )+B(1-\lambda ,\mu )+B(\mu ,\sigma )\in \mathbf{R}%
_{+}.
\end{eqnarray*}%
By (\ref{121}) and (\ref{129}), we have (cf. \cite{YM13})%
\begin{equation}
||T_{1}||=||T_{2}||=B(1-\lambda ,\sigma )+B(1-\lambda ,\mu )+B(\mu ,\sigma ).
\label{138}
\end{equation}
\end{example}

(b) Set $$H(t)=K_{\lambda }(1,t)=\frac{|\ln |t|^{\beta }|}{|1+t|^{\lambda }%
},$$
where $\beta >-1,\mu ,\sigma >0,\mu +\sigma =\lambda <1+\beta.$
Then we have
the kernels
\[
H(xy)=\frac{|\ln |xy|^{\beta }|}{|1+xy|^{\lambda }},\ \ K_{\lambda }(x,y)=\frac{%
|\ln |x/y|^{\beta }|}{|x+y|^{\lambda }}
\]
and obtain the constant factors%
\begin{eqnarray*}
K(\sigma ) &=&K_{\lambda }(\sigma )=\int_{-\infty }^{\infty }\frac{|\ln
|t|^{\beta }||t|^{\sigma -1}}{|1+t|^{\lambda }}dt \\
&=&\int_{0}^{\infty }\frac{|\ln t^{\beta }|t^{\sigma -1}}{|1-t|^{\lambda }}%
dt+\int_{0}^{\infty }\frac{|\ln t^{\beta }|t^{\sigma -1}}{(1+t)^{\lambda }}dt
\\
&=&\int_{0}^{1}(-\ln t)^{\beta }\left[\frac{1}{(1-t)^{\lambda }}+\frac{1}{%
(1+t)^{\lambda }}\right](t^{\mu -1}+t^{\sigma -1})dt.
\end{eqnarray*}%
There exists a constant $\delta _{\:0}>0,$ such that $\sigma >\delta _{\:0},\ \mu
>\delta _{\:0}.$ Since
\[
\lim_{t\rightarrow 0^{+}}t^{\delta _{\:0}}\frac{(-\ln t)^{\beta }}{%
(1-t)^{\beta }}=0,\ \ \lim_{t\rightarrow 1^{-}}t^{\delta _{\:0}}\frac{(-\ln
t)^{\beta }}{(1-t)^{\beta }}=1,
\]
there exists a constant $L>0,$ such that $$t^{\delta _{\:0}}\frac{(-\ln
t)^{\beta }}{(1-t)^{\beta }}\leq L\ (0<t<1)$$ and
\begin{eqnarray*}
0 &<&K(\sigma )=\int_{0}^{1}(-\ln t)^{\beta }[\frac{1}{(1-t)^{\lambda }}+%
\frac{1}{(1+t)^{\lambda }}](t^{\mu -1}+t^{\sigma -1})dt \\
&\leq &2\int_{0}^{1}\frac{(-\ln t)^{\beta }}{(1-t)^{\lambda }}(t^{\mu
-1}+t^{\sigma -1})dt\leq 2L\int_{0}^{1}\frac{t^{\mu -\delta
_{\:0}-1}+t^{\sigma -\delta _{\:0}-1}}{(1-t)^{\lambda -\beta }}dt \\
&=&2L[B(\beta +1-\lambda ,\mu -\delta _{\:0})+B(\beta +1-\lambda ,\sigma
-\delta _{\:0})]<\infty .
\end{eqnarray*}%
Therefore, $K(\sigma )=K_{\lambda }(\sigma )\in \mathbf{R}_{+}.$

 Since $$%
\binom{-\lambda}{2k}=\binom{\lambda+2k-1}{2k}>0,$$ then by Lebesgue's term
by term theorem, it follows that
\begin{eqnarray*}
K(\sigma ) &=&\int_{0}^{1}(-\ln t)^{\beta }\sum_{k=0}^{\infty
}(_{\,\,k}^{-\lambda})[(-1)^{k}+1](t^{k+\mu -1}+t^{k+\sigma -1})dt \\
&=&2\int_{0}^{1}(-\ln t)^{\beta }\sum_{k=0}^{\infty }(_{\,\,2k}^{-\lambda})
(t^{k+\mu -1}+t^{k+\sigma -1})dt \\
&=&2\sum_{k=0}^{\infty }(_{\,\,2k}^{-\lambda})\int_{0}^{1}(-\ln t)^{\beta
}(t^{2k+\mu -1}+t^{2k+\sigma -1})dt
\end{eqnarray*}%
\[
=2\Gamma (\beta +1)\sum_{k=0}^{\infty }(_{\,\,2k}^{-\lambda})\left[\frac{1}{%
(2k+\mu )^{\beta }}+\frac{1}{(2k+\sigma )^{\beta }}\right].
\]
In view of (\ref{121}) and (\ref{129}), we have%
\begin{equation}
||T_{1}||=||T_{2}||=2\Gamma (\beta +1)\sum_{k=0}^{\infty
}(_{\,\,2k}^{-\lambda})\left[\frac{1}{(2k+\mu )^{\beta }}+\frac{1}{(2k+\sigma
)^{\beta }}\right].  \label{139}
\end{equation}

(c) Set $$H(t)=K_{\lambda }(1,t)=\frac{(\max \{1,|t|\})^{\beta }}{%
|1+t|^{\lambda +\beta }}\ (\beta <1,\mu ,\sigma >0,\mu +\sigma =\lambda
<1-\beta ).$$ Then we have the kernels%
\[
H(xy)=\frac{(\max \{1,|xy|\})^{\beta }}{|1+xy|^{\lambda +\beta }},\ \ K_{\lambda
}(x,y)=\frac{(\max \{|x|,|y|\})^{\beta }}{|x+y|^{\lambda +\beta }}
\]
and obtain the constant factors%
\begin{eqnarray*}
K(\sigma ) &=&K_{\lambda }(\sigma )=\int_{-\infty }^{\infty }\frac{(\max
\{1,|t|\})^{\beta }}{|1+t|^{\lambda +\beta }}|t|^{\sigma -1}dt \\
&=&\int_{0}^{\infty }\frac{(\max \{1,t\})^{\beta }}{|1-t|^{\lambda +\beta }}%
t^{\sigma -1}dt+\int_{0}^{\infty }\frac{(\max \{1,t\})^{\beta }}{%
(1+t)^{\lambda +\beta }}t^{\sigma -1}dt\\
&=&\int_{0}^{1}\left[\frac{1}{(1-t)^{\lambda +\beta }}+\frac{1}{(1+t)^{\lambda
+\beta }}\right](t^{\mu -1}+t^{\sigma -1})dt \\
&=&B(1-\lambda -\beta ,\mu )+B(1-\lambda -\beta ,\sigma )+\int_{0}^{1}\frac{%
t^{\mu -1}+t^{\sigma -1}}{(1+t)^{\lambda +\beta }}dt\in \mathbf{R}_{+}.
\end{eqnarray*}

By Taylor's formula, we still can obtain
\begin{eqnarray*}
&&\int_{0}^{1}\frac{t^{\mu -1}+t^{\sigma -1}}{(1+t)^{\lambda +\beta }}%
dt=\int_{0}^{1}\sum_{k=0}^{\infty }\left( _{_{{}}\,\,\,k}^{-\lambda -\beta
}\right) (t^{k+\mu -1}+t^{k+\sigma -1})dt \\
&=&\int_{0}^{1}\sum_{k=0}^{\infty }(-1)^{k}\left( _{_{{}}\,\,\,k}^{\lambda
+\beta +k-1}\right) (t^{k+\mu -1}+t^{k+\sigma -1})dt \\
&=&\int_{0}^{1}\sum_{k=0}^{\infty }\left[\left( _{_{{}}\,\,\,2k}^{\lambda +\beta
+2k-1}\right) -\left( _{_{{}\,\,\,}2k+1}^{\lambda +\beta +2k}\right)
t\right](t^{2k+\mu -1}+t^{2k+\sigma -1})dt.
\end{eqnarray*}%
Since we find%
\begin{eqnarray*}
\left( _{_{{}}\,\,\,2k}^{\lambda +\beta +2k-1}\right) -\left(
_{_{{}}\,\,\,2k+1}^{\lambda +\beta +2k}\right) t &=&\left(
_{_{{}}\,\,\,2k}^{\lambda +\beta +2k-1}\right) -\frac{(\lambda +\beta +2k)t}{%
2k+1}\left( _{_{{}}\,\,\,2k}^{\lambda +\beta +2k-1}\right) \\
&=&\left[1-\frac{(\lambda +\beta +2k)t}{2k+1}\right]\left( _{_{{}}\,\,\,2k}^{\lambda
+\beta +2k-1}\right) ,
\end{eqnarray*}%
there exists a number $k_{0}\in \mathbf{N}_{0}=\mathbf{N}\cup \{0\},$ such
that $\lambda +\beta +2k_{0}>0$, and for any $s\in \mathbf{N},$%
\begin{eqnarray*}
&&\left( _{_{{}}\,\,\,2(k_{0}+s)}^{\lambda +\beta +2(k_{0}+s)-1}\right)
-\left( _{_{{}}\,\,\,2(k_{0}+s)+1}^{\lambda +\beta +2(k_{0}+s)+1}\right) t \\
&=&[1-\frac{(\lambda +\beta +2k_{0}+2s)t}{2(k_{0}+s)+1}]\left(
_{_{{}}\,\,\,2(k_{0}+s)}^{\lambda +\beta +2(k_{0}+s)-1}\right) =[1-\frac{%
(\lambda +\beta +2k_{0}+2s )t}{2(k_{0}+s)+1}] \\
&&\times \frac{\lambda +\beta +2k_{0}+2s-1}{2k_{0}+2s}\cdots \frac{\lambda
+\beta +2k_{0}}{2k_{0}+1}\left( _{_{{}}\,\,\,2k_{0}}^{\lambda +\beta
+2k_{0}-1}\right) .
\end{eqnarray*}%
For $t\in (0,1],$ we get
\begin{eqnarray*}
1-\frac{(\lambda +\beta +2k_{0}+2s )t}{2(k_{0}+s)+1}&\geq &1-\frac{\lambda
+\beta +2k_{0}+2s }{2(k_{0}+s)+1}\\
&=&\frac{1-\lambda -\beta }{2(k_{0}+s)+1}>0.
\end{eqnarray*}
Then it follows that for any $s\in \mathbf{N},$%
\[
sgn\left(\left( _{_{{}}\,\,\,2(k_{0}+s)}^{\lambda +\beta +2(k_{0}+s)-1}\right)
-\left( _{_{{}}\,\,\,2(k_{0}+s)+1}^{\lambda +\beta +2(k_{0}+s)+1}\right)
t\right)=sgn\left( _{_{{}}\,\,\,2k_{0}}^{\lambda +\beta +2k_{0}-1}\right) .
\]
Hence by Lebesgue term by term integration theorem, we have%
\begin{eqnarray*}
K(\sigma ) &=&B(1-\lambda -\beta ,\mu )+B(1-\lambda -\beta ,\sigma ) \\
&&+\sum_{k=0}^{\infty } \left( _{_{{}}\,\,\,k}^{-\lambda -\beta
}\right)\int_{0}^{1}(t^{k+\mu -1}+t^{k+\sigma -1})dt \\
&=&B(1-\lambda -\beta ,\mu )+B(1-\lambda -\beta ,\sigma ) \\
&&+\sum_{k=0}^{\infty } \left( _{_{{}}\,\,\,k}^{-\lambda -\beta }\right)(%
\frac{1}{k+\mu }+\frac{1}{k+\sigma }).
\end{eqnarray*}
In view of (\ref{121}) and (\ref{129}), we have%
\begin{eqnarray}
||T_{1}|| &=&||T_{2}||=B(1-\lambda -\beta ,\beta +\mu )+B(1-\lambda -\beta
,\beta +\sigma )  \nonumber \\
&&+\sum_{k=0}^{\infty } \left( _{_{{}}\,\,\,k}^{-\lambda -\beta }\right)\left(%
\frac{1}{k+\mu }+\frac{1}{k+\sigma }\right).  \label{140}
\end{eqnarray}

For (a)-(c), we can obtain the equivalent inequalities with the kernels
and the best possible constant factors in Theorem 5-8. Setting $\delta _{\:0}=%
\frac{\sigma }{2}>0,$ we also obtain the equivalent reverse
inequalities with the kernels and the best possible constant factors in
Theorem 5-8.

(d) Set $$H(t)=K_{\lambda }(1,t)=\frac{(\min \{1,|t|\})^{\beta }}{%
|1+t|^{\lambda +\beta }},$$
with $\beta >-1,\mu ,\sigma >-\beta ,\mu +\sigma =\lambda
<1-\beta .$ Then we have the kernels%
\[
H(xy)=\frac{(\min \{1,|xy|\})^{\beta }}{|1+xy|^{\lambda +\beta }},\ \ K_{\lambda
}(x,y)=\frac{(\min \{|x|,|y|\})^{\beta }}{|x+y|^{\lambda +\beta }}
\]
and obtain the constant factors%
\begin{eqnarray*}
K(\sigma ) &=&K_{\lambda }(\sigma )=\int_{-\infty }^{\infty }\frac{(\min
\{1,|t|\})^{\beta }}{|1+t|^{\lambda +\beta }}|t|^{\sigma -1}dt \\
&=&\int_{0}^{\infty }\frac{(\min \{1,t\})^{\beta }}{|1-t|^{\lambda +\beta }}%
t^{\sigma -1}dt+\int_{0}^{\infty }\frac{(\min \{1,t\})^{\beta }}{%
(1+t)^{\lambda +\beta }}t^{\sigma -1}dt \\
&=&\int_{0}^{1}[\frac{1}{(1-t)^{\lambda +\beta }}+\frac{1}{(1+t)^{\lambda
+\beta }}](t^{\beta +\mu -1}+t^{\beta +\sigma -1})dt \\
&=&B(1-\lambda -\beta ,\beta +\mu )+B(1-\lambda -\beta ,\beta +\sigma ) \\
&&+\int_{0}^{1}\frac{t^{\beta +\mu -1}+t^{\beta +\sigma -1}}{(1+t)^{\lambda
+\beta }}dt\in \mathbf{R}_{+}.
\end{eqnarray*}

Similarly to the method followed in (c), we find
\begin{eqnarray*}
\int_{0}^{1}\frac{t^{\beta +\mu -1}+t^{\beta +\sigma -1}}{(1+t)^{\lambda
+\beta }}dt &=&\sum_{k=0}^{\infty } \left( _{_{{}}\,\,\,k}^{-\lambda -\beta
}\right)\int_{0}^{1}(t^{k+\beta +\mu -1}+t^{k+\beta +\sigma -1})dt \\
&=&\sum_{k=0}^{\infty } \left( _{_{{}}\,\,\,k}^{-\lambda -\beta }\right)\left(%
\frac{1}{k+\beta +\mu }+\frac{1}{k+\beta +\sigma }\right).
\end{eqnarray*}%
By the above results, (\ref{121}) and (\ref{129}), we have%
\begin{eqnarray}
||T_{1}|| &=&||T_{2}||=B(1-\lambda -\beta ,\beta +\mu )+B(1-\lambda -\beta
,\beta +\sigma )  \nonumber \\
&&+\sum_{k=0}^{\infty } \left( _{_{{}}\,\,\,k}^{-\lambda -\beta }\right)\left(%
\frac{1}{k+\beta +\mu }+\frac{1}{k+\beta +\sigma }\right).  \label{141}
\end{eqnarray}

Then in (d), we can obtain the equivalent inequalities with the kernels and
the best possible constant factors in Theorem 5-8. Setting $\delta _{\:0}=%
\frac{\sigma +\beta }{2}>0,$ we can still obtain the equivalent reverse
inequalities with the kernels and the best possible constant factors in
Theorem 5-8.

\begin{example}
 Set $$H(t)=K_{2}(1,t)=\frac{1}{1+2bt+(ct)^{2}}\
(|b|<|c|,\mu =\sigma =1).$$ Then we have the kernels $$
H(xy)=\frac{1}{1+2bxy+(cxy)^{2}},\ \ K_{\lambda }(x,y)=\frac{1}{%
x^{2}+2bxy+(cy)^{2}} $$
and obtain the constant factors
\begin{eqnarray*}
K(\sigma ) &=&K_{2}(\sigma )=\int_{-\infty }^{\infty }\frac{1}{1+2bt+(ct)^{2}%
}dt \\
&=&\frac{2}{\sqrt{4c^{2}-4b^{2}}}\arctan \frac{2c^{2}t+2b}{\sqrt{%
4c^{2}-4b^{2}}}\Big|_{-\infty }^{\infty }
=\frac{\pi }{\sqrt{c^{2}-b^{2}}}\in \mathbf{R}_{+}.
\end{eqnarray*}%
In view of (\ref{121}) and (\ref{129}), we have (cf. \cite{YM10}) $$
||T_{1}||=||T_{2}||=\frac{\pi }{\sqrt{c^{2}-b^{2}}}.$$
\end{example}

In particular, for $c=1,\:b=\cos \alpha\ (0<\alpha <\pi ),$ we have $$%
||T_{1}||=||T_{2}||=\frac{\pi }{\sin \alpha }.$$

\begin{example}$ $
 (a) Set
$$
H(t)=K_{2}(1,t)=\min_{i\in \{1,2\}}\frac{1}{1+2t\cos \alpha _{i}+t^{2}},
$$
with
$0<\alpha _{1}\leq \alpha _{2}<\pi ,0<\sigma <2.$ Then we have the kernels%
$$
H(xy)=\min_{i\in \{1,2\}}\frac{1}{1+2xy\cos \alpha _{i}+(xy)^{2}},\ \
K_{2}(x,y)=\min_{i\in \{1,2\}}\frac{1}{x^{2}+2xy\cos \alpha _{i}+y^{2}}
$$
and obtain the constant factors%
\begin{eqnarray*}
K(\sigma ) &=&K_{2}(\sigma )=\int_{-\infty }^{\infty }\min_{i\in \{1,2\}}%
\frac{1}{1+2t\cos \alpha _{i}+t^{2}}|t|^{\sigma -1}dt \\
&=&\int_{0}^{\infty }\min_{i\in \{1,2\}}\frac{t^{\sigma -1}}{1+2t\cos \alpha
_{i}+t^{2}}dt+\int_{0}^{\infty }\min_{i\in \{1,2\}}\frac{t^{\sigma -1}}{1-2t\cos \alpha
_{i}+t^{2}}dt
\end{eqnarray*}%
$$
=\int_{0}^{\infty }\frac{t^{\sigma -1}}{1+2t\cos \alpha _{1}+t^{2}}dt
+\int_{0}^{\infty }\frac{t^{\sigma -1}}{1+2t\cos (\pi -\alpha _{2})+t^{2}}%
dt.
$$
\end{example}

\noindent Set $$f(z)=\frac{1}{1+2z\cos \alpha _{1}+z^{2}}.$$ Then $$z_{1}=-e^{i\alpha
_{1}},\ \ z_{2}=-e^{-i\alpha _{1}}$$ are the poles of order 1. Setting $$\varphi
_{1}(z)=(z-z_{1})f(z)=\frac{1}{z-z_{2}},\ \varphi _{2}(z)=(z-z_{2})f(z)=\frac{1%
}{z-z_{1}},$$ by (\ref{63}), we have%
\begin{eqnarray*}
&&\int_{0}^{\infty }\frac{t^{\sigma -1}}{1+2t\cos \alpha _{1}+t^{2}}%
dt=\int_{0}^{\infty }f(t)t^{\sigma -1}dt \\
&=&\frac{\pi }{\sin \pi \sigma }[(-z_{1})^{\sigma -1}\varphi
_{1}(z_{1})+(-z_{2})^{\sigma -1}\varphi _{2}(z_{2})] \\
&=&\frac{\pi }{\sin \pi \sigma }\left[\frac{e^{i\alpha _{1}(\sigma -1)}}{%
-e^{i\alpha _{1}}+e^{-i\alpha _{1}}}+\frac{e^{-i\alpha _{1}(\sigma -1)}}{%
-e^{-i\alpha _{1}}+e^{i\alpha _{1}}}\right] =\frac{\pi \sin \alpha _{1}(1-\sigma )%
}{\sin \pi \sigma \sin \alpha _{1}}.
\end{eqnarray*}%
Similarly, it follows that
\begin{eqnarray*}
&&\int_{0}^{\infty }\frac{t^{\sigma -1}dt}{1+2t\cos (\pi -\alpha _{2})+t^{2}}
\\
&=&\frac{\pi \sin (\pi -\alpha _{2})(1-\sigma )}{\sin \pi \sigma \sin (\pi
-\alpha _{2})}=\frac{\pi \sin (\pi -\alpha _{2})(1-\sigma )}{\sin \pi \sigma
\sin \alpha _{2}},
\end{eqnarray*}%
and then%
\[
K(\sigma )=\frac{\pi }{\sin \pi \sigma }\left[\frac{\sin \alpha _{1}(1-\sigma )}{%
\sin \alpha _{1}}+\frac{\sin (\pi -\alpha _{2})(1-\sigma )}{\sin \alpha _{2}}%
\right]\in \mathbf{R}_{+}.
\]
In view of (\ref{121}) and (\ref{129}), we have (cf. \cite{YM18})%
\begin{equation}
||T_{1}||=||T_{2}||=\frac{\pi }{\sin \pi \sigma }\left[\frac{\sin \alpha
_{1}(1-\sigma )}{\sin \alpha _{1}}+\frac{\sin (\pi -\alpha _{2})(1-\sigma )}{%
\sin \alpha _{2}}\right].  \label{142}
\end{equation}

Then in (a), we can obtain the equivalent inequalities with the kernels and
the best possible constant factors in Theorem 5-8. Setting $\delta _{\:0}=%
\frac{\sigma }{2}>0,$ we can still obtain the equivalent reverse
inequalities with the kernels and the best possible constant factors in
Theorem 5-8.

In particular, if $\alpha _{1}=\alpha _{2}=\alpha \in (0,\pi ),$ then $$%
H(t)=K_{2}(1,t)=\frac{1}{1+2t\cos \alpha +t^{2}}$$ and%
\[
||T_{1}||=||T_{2}||=\frac{\pi }{\sin \pi \sigma \sin \alpha }[\sin \alpha
(1-\sigma )+\sin (\pi -\alpha )(1-\sigma )].
\]

(b) Set
\[
H(t)=K_{0}(1,t)=\min_{i\in \{1,2\}}\frac{\min \{1,|t|\}}{\sqrt{1+2t\cos
\alpha _{i}+t^{2}}},
\]
with $0<\alpha _{1}\leq \alpha _{2}<\pi ,\sigma =\mu =0.$ Then we have the
kernels%
\begin{eqnarray*}
H(xy) &=&\min_{i\in \{1,2\}}\frac{\min \{1,|xy|\}}{\sqrt{1+2xy\cos \alpha
_{i}+(xy)^{2}}}, \\
K_{0}(x,y) &=&\min_{i\in \{1,2\}}\frac{\min \{|x|,|y|\}}{\sqrt{x^{2}+2xy\cos
\alpha _{i}+y^{2}}}
\end{eqnarray*}%
and obtain the constant factors%
\begin{eqnarray*}
K(0) &=&K_{0}(0)=\int_{-\infty }^{\infty }\min_{i\in \{1,2\}}\frac{\min
\{1,|t|\}}{\sqrt{1+2t\cos \alpha _{i}+t^{2}}}|t|^{-1}dt \\
&=&\int_{0}^{\infty }\min_{i\in \{1,2\}}\frac{\min \{1,t\}t^{-1}}{\sqrt{%
1+2t\cos \alpha _{i}+t^{2}}}dt \\
&&+\int_{0}^{\infty }\min_{i\in \{1,2\}}\frac{\min \{1,t\}t^{-1}}{\sqrt{%
1-2t\cos \alpha _{i}+t^{2}}}dt \\
&=&\int_{0}^{\infty }\frac{\min \{1,t\}t^{-1}dt}{\sqrt{1+2t\cos \alpha
_{1}+t^{2}}}+\int_{0}^{\infty }\frac{\min \{1,t\}t^{-1}dt}{\sqrt{1+2t\cos
(\pi -\alpha _{2})+t^{2}}} \\
&=&2\left[ \int_{0}^{1}\frac{dt}{\sqrt{1+2t\cos \alpha _{1}+t^{2}}}%
+\int_{0}^{1}\frac{dt}{\sqrt{1+2t\cos (\pi -\alpha _{2})+t^{2}}}\right] .
\end{eqnarray*}%
We get%
\begin{eqnarray*}
&&\int_{0}^{1}\frac{dt}{\sqrt{1+2t\cos \alpha _{1}+t^{2}}} \\
&=&\ln (2t+2\cos \alpha _{1}+2\sqrt{1+2t\cos \alpha _{1}+t^{2}}%
)|_{0}^{1}=\ln \left(1+\sec \frac{\alpha _{1}}{2}\right),
\end{eqnarray*}%
and by the same way,%
\[
\int_{0}^{1}\frac{dt}{\sqrt{1+2t\cos (\pi -\alpha _{2})+t^{2}}} =\ln (1+\sec
\frac{\pi -\alpha _{2}}{2})=\ln \left(1+\csc \frac{\alpha _{2}}{2}\right).
\]
Then it follows that $$K(0)=K_{0}(0)=2\ln \left(1+\sec \frac{\alpha _{1}}{2}%
\right)\left(1+\csc \frac{\alpha _{2}}{2}\right). $$ By (\ref{121}) and (\ref{129}),
we have (cf. \cite{YM17})%
\begin{equation}
||T_{1}||=||T_{2}||=2\ln \left(1+\sec \frac{\alpha _{1}}{2}\right)\left(1+\csc \frac{\alpha
_{2}}{2}\right).  \label{143}
\end{equation}

In particular, if $\alpha _{1}=\alpha _{2}=\alpha \in (0,\pi ),$ then $$%
H(t)=K_{0}(1,t)=\frac{\min \{1,|t|\}}{\sqrt{1+2t\cos \alpha +t^{2}}}$$ and
\[
||T_{1}||=||T_{2}||=2\ln \left(1+\sec \frac{\alpha }{2}\right)\left(1+\csc \frac{\alpha }{2}%
\right).
\]
\begin{example}$ $
 Set
\begin{eqnarray*}
H(t)=K_{0}(1,t)=\Big|\ln \frac{1+2t\cos \alpha _{1}+t^{2}}{1+2t\cos \alpha
_{2}+t^{2}}\Big|,
\end{eqnarray*}%
with $0<\alpha _{1}\leq \alpha _{2}<\pi ,\ \mu =-\sigma \in (0,1).$ Then we have
the kernels%
\begin{eqnarray*}
H(xy) &=&\Big|\ln \frac{1+2xy\cos \alpha _{1}+(xy)^{2}}{1+2xy\cos \alpha
_{2}+(xy)^{2}}\Big|, \\
K_{0}(x,y) &=&\Big|\ln \frac{x^{2}+2xy\cos \alpha _{1}+y^{2}}{x^{2}+2xy\cos
\alpha _{2}+y^{2}}\Big|
\end{eqnarray*}%
and obtain the constant factors%
\begin{eqnarray*}
K(\sigma ) &=&K_{0}(\sigma )=\int_{-\infty }^{\infty }\Big|\ln \frac{1+2t\cos
\alpha _{1}+t^{2}}{1+2t\cos \alpha _{2}+t^{2}}\Big|\cdot |t|^{\sigma -1}dt \\
&=&\int_{0}^{\infty }t^{\sigma -1}\ln \frac{1+2t\cos \alpha _{1}+t^{2}}{%
1+2t\cos \alpha _{2}+t^{2}}dt\\
&&+\int_{0}^{\infty }t^{\sigma -1}\ln \frac{1-2t\cos \alpha _{2}+t^{2}}{%
1-2t\cos \alpha _{1}+t^{2}}dt\  \\
&=&\int_{0}^{\infty }t^{\sigma -1}\ln \frac{1+2t\cos \alpha _{1}+t^{2}}{%
1+2t\cos \alpha _{2}+t^{2}}dt \\
&&+\int_{0}^{\infty }t^{\sigma -1}\ln \frac{1+2t\cos (\pi -\alpha _{2})+t^{2}%
}{1+2t\cos (\pi -\alpha _{1})+t^{2}}dt.
\end{eqnarray*}
\end{example}

We find%
\begin{eqnarray*}
I_{1} &:&=\int_{0}^{\infty }t^{\sigma -1}\ln \frac{1+2t\cos \alpha _{1}+t^{2}%
}{1+2t\cos \alpha _{2}+t^{2}}dt =\frac{1}{\sigma }\int_{0}^{\infty }\ln
\frac{1+2t\cos \alpha _{1}+t^{2}}{1+2t\cos \alpha _{2}+t^{2}}dt^{\sigma } \\
&=&\frac{1}{\sigma }\left[ t^{\sigma }\ln \frac{1+2t\cos \alpha _{1}+t^{2}}{%
1+2t\cos \alpha _{2}+t^{2}}|_{0}^{\infty }\right. \\
&&\left. -2\int_{0}^{\infty }t^{\sigma }(\frac{\cos \alpha _{1}+t}{1+2t\cos
\alpha _{1}+t^{2}}-\frac{\cos \alpha _{2}+t}{1+2t\cos \alpha _{2}+t^{2}})dt%
\right] \\
&=&\frac{-2}{\sigma }\left[ \int_{0}^{\infty }\frac{(\cos \alpha
_{1}+t)t^{\sigma }}{1+2t\cos \alpha _{1}+t^{2}}dt-\int_{0}^{\infty }\frac{%
(\cos \alpha _{2}+t)t^{\sigma }}{1+2t\cos \alpha _{2}+t^{2}}dt\right] \\
&=&\frac{-2}{\sigma }(I_{1}^{(1)}-I_{1}^{(2)}), \\
I_{1}^{(i)} &=&\int_{0}^{\infty }\frac{(\cos \alpha _{i}+t)t^{(\sigma +1)-1}%
}{1+2t\cos \alpha _{i}+t^{2}}dt\ (i=1,2).
\end{eqnarray*}

By (\ref{63}), we have%
\begin{eqnarray*}
I_{1}^{(i)} &=&\frac{\pi }{\sin \pi (\sigma +1)}\left(e^{i\alpha _{i}\sigma }%
\frac{\cos \alpha _{i}-e^{i\alpha _{i}}}{-e^{i\alpha _{i}}+e^{-i\alpha _{i}}}%
+e^{-i\alpha _{i}\sigma }\frac{\cos \alpha _{i}-e^{-i\alpha _{i}}}{%
-e^{-i\alpha _{i}}+e^{i\alpha _{i}}}\right) \\
&=&\frac{\pi \cos \alpha _{i}\sigma }{\sin \pi (\sigma +1)}\ (i=1,2),
\end{eqnarray*}%
and then%
\begin{eqnarray*}
I_{1} &=&\frac{-2}{\sigma }\left[\frac{\pi \cos \alpha _{1}\sigma }{\sin \pi
(\sigma +1)}-\frac{\pi \cos \alpha _{2}\sigma }{\sin \pi (\sigma +1)}\right] \\
&=&\frac{-2\pi (\cos \alpha _{1}\sigma -\cos \alpha _{2}\sigma )}{\sigma
\sin \pi (\sigma +1)} \\
&=&\frac{4\pi }{\sigma \sin \pi (\sigma +1)}\sin \frac{\alpha _{1}+\alpha
_{2}}{2}\sigma \sin \frac{\alpha _{1}-\alpha _{2}}{2}\sigma .
\end{eqnarray*}

Similarly, we have
\begin{eqnarray*}
I_{2} &:&=\int_{0}^{\infty }t^{\sigma -1}\ln \frac{1+2t\cos (\pi -\alpha
_{2})+t^{2}}{1+2t\cos (\pi -\alpha _{1})+t^{2}}dt \\
&=&\frac{4\pi }{\sigma \sin \pi (\sigma +1)}\sin \left(\pi -\frac{\alpha
_{1}+\alpha _{2}}{2}\right)\sigma \sin \frac{\alpha _{1}-\alpha _{2}}{2}\sigma ,
\end{eqnarray*}%
and then%
\begin{eqnarray*}
K(\sigma ) &=&K_{0}(\sigma )=\frac{4\pi }{\sigma \sin \pi (\sigma +1)}\sin
\frac{\alpha _{1}+\alpha _{2}}{2}\sigma \\
&&\times \left[ \sin \frac{\alpha _{1}-\alpha _{2}}{2}\sigma +\sin (\pi -%
\frac{\alpha _{1}+\alpha _{2}}{2})\sigma \right] \\
&=&\frac{-4\pi \sin \frac{\sigma }{2}(\alpha _{1}-\alpha _{2})}{\sigma \cos
\pi (\sigma /2)}\cos \frac{\sigma }{2}(\alpha _{1}+\alpha _{2}-\pi )\in
\mathbf{R}_{+}.
\end{eqnarray*}%
In view of (\ref{121}) and (\ref{129}), we have (cf. \cite{YM15})%
\begin{equation}
||T_{1}||=||T_{2}||=\frac{-4\pi \sin \frac{\sigma }{2}(\alpha _{1}-\alpha
_{2})}{\sigma \cos \pi (\sigma /2)}\cos \frac{\sigma }{2}(\alpha _{1}+\alpha
_{2}-\pi ).  \label{144}
\end{equation}

Then we can obtain the equivalent inequalities with the kernels and the best
possible constant factors in Theorem 5-8. Setting $\delta _{\:0}=\frac{-\sigma
}{2}>0,$ we still can obtain the equivalent reverse inequalities with the
kernels and the best possible constant factors in Theorem 5-8.

\begin{remark}
 Since $K(0^{-})=2\pi (\alpha _{2}-\alpha _{1})\in \mathbf{%
R}_{+},$ then (\ref{144}) is valid for $\sigma \in (-1,0].$
\end{remark}

\begin{example}$ $
 (a) For $$\gamma \in \{a;\:a=\frac{1}{2k-1},2k+1(k\in
\mathbf{N})\},$$ we set
$$
H(t)=K_{\lambda }(1,t)=\min_{i\in \{1,2\}}\frac{1}{(1+t^{\gamma }\cos \alpha
_{i}+|t|^{\gamma })^{\lambda /\gamma }},
$$
where $0<\alpha _{1}\leq \alpha _{2}<\pi ,\:\mu ,\:\sigma >0,\:\mu +\sigma =\lambda .$
Then we have the kernels%
\begin{eqnarray*}
H(xy) &=&\min_{i\in \{1,2\}}\frac{1}{(1+(xy)^{\gamma }\cos \alpha
_{i}+|xy|^{\gamma })^{\lambda /\gamma }}, \\
K_{\lambda }(x,y) &=&\min_{i\in \{1,2\}}\frac{1}{(|x|^{\gamma }+y^{\gamma
}sgn\text{(}x)\cos \alpha _{i}+|y|^{\gamma })^{\lambda /\gamma }}
\end{eqnarray*}%
and obtain the constant factors%
\begin{eqnarray*}
K(\sigma )&=&K(\sigma )=K_{\lambda }(\sigma )=\int_{-\infty }^{\infty }\min_{i\in \{1,2\}}%
\frac{1}{(1+t^{\gamma }\cos \alpha _{i}+|t|^{\gamma })^{\lambda /\gamma }}%
|t|^{\sigma -1}dt  \\
&=&\int_{0}^{\infty }\frac{t^{\sigma -1}dt}{[1+t^{\gamma }(1+\cos \alpha
_{1})]^{\lambda /\gamma }}+\int_{0}^{\infty }\frac{t^{\sigma -1}dt}{%
[1+t^{\gamma }(1-\cos \alpha _{2})]^{\lambda /\gamma }} \\
&=&\frac{1}{\gamma }\left[\frac{1}{(1+\cos \alpha _{1})^{\lambda /\gamma }}+\frac{%
1}{(1-\cos \alpha _{2})^{\lambda /\gamma }}\right]\int_{0}^{\infty }\frac{%
u^{(\sigma /\gamma )-1}du}{(1+u)^{\lambda /\gamma }} \\
&=&\frac{1}{\gamma 2^{\sigma /\gamma }}\left[(\sec \frac{\alpha _{1}}{2})^{\frac{%
2\sigma }{\gamma }}+(\csc \frac{\alpha _{2}}{2})^{\frac{2\sigma }{\gamma }%
}\right]B\left(\frac{\mu }{\gamma },\frac{\sigma }{\gamma }\right).
\end{eqnarray*}
\end{example}

\bigskip In view of (\ref{121}) and (\ref{129}), we have%
\begin{equation}
||T_{1}||=||T_{2}||=\frac{1}{\gamma 2^{\sigma /\gamma }}\left[(\sec \frac{\alpha
_{1}}{2})^{\frac{2\sigma }{\gamma }}+(\csc \frac{\alpha _{2}}{2})^{\frac{%
2\sigma }{\gamma }}\right]B\left(\frac{\mu }{\gamma },\frac{\sigma }{\gamma }\right).
\label{145}
\end{equation}

In particular, if $\alpha _{1}=\alpha _{2}=\alpha \in (0,\pi ),$ then it follows that
\[
H(t)=K_{\lambda }(1,t)=\frac{1}{(1+t^{\gamma }\cos \alpha +|t|^{\gamma
})^{\lambda /\gamma }},
\]
and
\[
||T_{1}||=||T_{2}||=\frac{1}{\gamma 2^{\sigma /\gamma }}\left[(\sec \frac{\alpha
}{2})^{\frac{2\sigma }{\gamma }}+(\csc \frac{\alpha }{2})^{\frac{2\sigma }{%
\gamma }}\right]B\left(\frac{\mu }{\gamma },\frac{\sigma }{\gamma }\right).
\]

(b) For $$\gamma \in \{a;a=\frac{1}{2k-1},2k+1(k\in \mathbf{N})\},$$ we set
\[
H(t)=K_{\lambda }(1,t)=\min_{i\in \{1,2\}}\frac{1}{|1-t^{\gamma }\cos \alpha
_{i}-|t|^{\gamma }|^{\lambda /\gamma }},
\]
where $0<\alpha _{1}\leq \alpha _{2}<\pi ,\:\mu ,\:\sigma >0,\:\mu +\sigma =\lambda
<\gamma .$ Then we have the kernels%
\begin{eqnarray*}
H(xy) &=&\min_{i\in \{1,2\}}\frac{1}{|1-(xy)^{\gamma }\cos \alpha
_{i}-|xy|^{\gamma }|^{\lambda /\gamma }}, \\
K_{\lambda }(x,y) &=&\min_{i\in \{1,2\}}\frac{1}{||x|^{\gamma }-y^{\gamma
}sgn\text{(}x)\cos \alpha _{i}-|y|^{\gamma }|^{\lambda /\gamma }}
\end{eqnarray*}%
and obtain the constant factors%
\begin{eqnarray*}
K(\sigma ) &=&K_{\lambda }(\sigma )=\int_{-\infty }^{\infty }\min_{i\in
\{1,2\}}\frac{1}{|1-t^{\gamma }\cos \alpha _{i}-|t|^{\gamma }|^{\lambda
/\gamma }}|t|^{\sigma -1}dt \\
&=&\int_{0}^{\infty }\min_{i\in \{1,2\}}\frac{1}{|1-t^{\gamma }(1+\cos
\alpha _{i})|^{\lambda /\gamma }}t^{\sigma -1}dt \\
&&+\int_{0}^{\infty }\min_{i\in \{1,2\}}\frac{1}{|1-t^{\gamma }(1-\cos
\alpha _{i})|^{\lambda /\gamma }}t^{\sigma -1}dt \\
&=&\frac{1}{\gamma }\left[ \int_{0}^{\infty }\frac{1}{(1+\cos \alpha
_{1})^{\lambda /\gamma }}\int_{0}^{\infty }\frac{u^{(\sigma /\gamma )-1}du}{%
|1-u|^{\lambda /\gamma }}\right] \\
&&+\frac{1}{\gamma }\left[ \int_{0}^{\infty }\frac{1}{(1-\cos \alpha
_{2})^{\lambda /\gamma }}\int_{0}^{\infty }\frac{u^{(\sigma /\gamma )-1}du}{%
|1-u|^{\lambda /\gamma }}\right] \\
&=&\frac{1}{\gamma 2^{\sigma /\gamma }}\left[\left(\sec \frac{\alpha _{1}}{2}\right)^{\frac{%
2\sigma }{\gamma }}+\left(\csc \frac{\alpha _{2}}{2}\right)^{\frac{2\sigma }{\gamma }}\right]
\\
&&\times \left[B(1-\frac{\lambda }{\gamma },\frac{\mu }{\gamma })+B\left(1-\frac{%
\lambda }{\gamma },\frac{\sigma }{\gamma }\right)\right].
\end{eqnarray*}%
In view of (\ref{121}) and (\ref{129}), we have%
\begin{eqnarray}
||T_{1}|| &=&||T_{2}||=\frac{1}{\gamma 2^{\sigma /\gamma }}\left[\left(\sec \frac{%
\alpha _{1}}{2}\right)^{\frac{2\sigma }{\gamma }}+\left(\csc \frac{\alpha _{2}}{2}\right)^{%
\frac{2\sigma }{\gamma }}\right]  \nonumber \\
&&\times \left[ B\left(1-\frac{\lambda }{\gamma },\frac{\mu }{\gamma }\right)+B\left(1-%
\frac{\lambda }{\gamma },\frac{\sigma }{\gamma }\right)\right].  \label{146}
\end{eqnarray}

In particular, if $\alpha _{1}=\alpha _{2}=\alpha \in (0,\pi ),$ then we
have
\[
H(t)=K_{\lambda }(1,t)=\frac{1}{|1-t^{\gamma }\cos \alpha -|t|^{\gamma
}|^{\lambda /\gamma }},
\]
and
\begin{eqnarray*}
||T_{1}|| &=&||T_{2}||=\frac{1}{\gamma 2^{\sigma /\gamma }}\left[\left(\sec \frac{%
\alpha }{2}\right)^{\frac{2\sigma }{\gamma }}+\left(\csc \frac{\alpha }{2}\right)^{\frac{%
2\sigma }{\gamma }}\right] \\
&&\times \left[ B\left(1-\frac{\lambda }{\gamma },\frac{\mu }{\gamma }\right)+B\left(1-%
\frac{\lambda }{\gamma },\frac{\sigma }{\gamma }\right)\right].
\end{eqnarray*}

(c) For $$\gamma \in \{a;a=\frac{1}{2k-1},2k+1\ (k\in \mathbf{N})\},$$ we set
\[
H(t)=K_{\lambda }(1,t)=\min_{i\in \{1,2\}}\frac{\ln (t^{\gamma }\cos \alpha
_{i}+|t|^{\gamma })}{(t^{\gamma }\cos \alpha _{i}+|t|^{\gamma })^{\lambda
/\gamma }-1},
\]
with $0<\alpha _{1}\leq \alpha _{2}<\pi ,\:\mu ,\:\sigma >0,\:\mu +\sigma =\lambda.$
Then we have the kernels%
\begin{eqnarray*}
H(xy) &=&\min_{i\in \{1,2\}}\frac{\ln ((xy)^{\gamma }\cos \alpha
_{i}+|xy|^{\gamma })}{((xy)^{\gamma }\cos \alpha _{i}+|xy|^{\gamma
})^{\lambda /\gamma }-1}, \\
K_{\lambda }(x,y) &=&\min_{i\in \{1,2\}}\frac{\ln ((\frac{y}{x})^{\gamma
}\cos \alpha _{i}+|\frac{y}{x}|^{\gamma })}{(y^{\gamma }sgn(x)\cos \alpha
_{i}+|y|^{\gamma })^{\lambda /\gamma }-|x|^{\lambda }}
\end{eqnarray*}%
and obtain the constant factors%
\begin{eqnarray*}
K(\sigma ) &=&K_{\lambda }(\sigma )=\int_{-\infty }^{\infty }\min_{i\in
\{1,2\}}\frac{\ln (t^{\gamma }\cos \alpha _{i}+|t|^{\gamma })}{(t^{\gamma
}\cos \alpha _{i}+|t|^{\gamma })^{\lambda /\gamma }-1}|t|^{\sigma -1}dt \\
&=&\int_{0}^{\infty }\min_{i\in \{1,2\}}\frac{\ln [t^{\gamma }(1+\cos \alpha
_{i})]}{[t^{\gamma }(1+\cos \alpha _{i})]^{\lambda /\gamma }-1}t^{\sigma
-1}dt \\
&&+\int_{0}^{\infty }\min_{i\in \{1,2\}}\frac{\ln [t^{\gamma }(1-\cos \alpha
_{i})]}{[t^{\gamma }(1-\cos \alpha _{i})]^{\lambda /\gamma }-1}t^{\sigma
-1}dt \\
&=&\frac{\gamma }{\lambda ^{2}}\left[ \int_{0}^{\infty }\frac{1}{(1+\cos
\alpha _{1})^{\sigma /\gamma }}\int_{0}^{\infty }\frac{\ln u}{u-1}u^{\frac{%
\sigma }{\lambda }-1}du\right] \\
&&+\frac{\gamma }{\lambda ^{2}}\left[ \int_{0}^{\infty }\frac{1}{(1-\cos
\alpha _{2})^{\sigma /\gamma }}\int_{0}^{\infty }\frac{\ln u}{u-1}u^{\frac{%
\sigma }{\lambda }-1}du\right] \\
&=&\frac{\gamma }{\lambda ^{2}2^{\sigma /\gamma }}\left[\left(\sec \frac{\alpha _{1}}{2%
}\right)^{\frac{2\sigma }{\gamma }}+\left(\csc \frac{\alpha _{2}}{2}\right)^{\frac{2\sigma }{%
\gamma }}\right]\left[\frac{\pi }{\sin \pi \left(\sigma /\lambda \right)}\right]^{2}.
\end{eqnarray*}%
By (\ref{121}) and (\ref{129}), we have%
\[
||T_{1}||=||T_{2}||=\frac{\gamma }{\lambda ^{2}2^{\sigma /\gamma }}\left[\left(\sec
\frac{\alpha _{1}}{2}\right)^{\frac{2\sigma }{\gamma }}+\left(\csc \frac{\alpha _{2}}{2}%
\right)^{\frac{2\sigma }{\gamma }}\right]\left[\frac{\pi }{\sin \pi (\sigma /\lambda )}\right]^{2}.
\]

In particular, if $\alpha _{1}=\alpha _{2}=\alpha \in (0,\pi ),$ then we
have
\[
H(t)=K_{\lambda }(1,t)=\frac{\ln (t^{\gamma }\cos \alpha +|t|^{\gamma })}{%
(t^{\gamma }\cos \alpha +|t|^{\gamma })^{\lambda /\gamma }-1},
\]
and
\[
||T_{1}||=||T_{2}||=\frac{\gamma }{\lambda ^{2}2^{\sigma /\gamma }}\left[\left(\sec
\frac{\alpha }{2}\right)^{\frac{2\sigma }{\gamma }}+\left(\csc \frac{\alpha }{2}\right)^{%
\frac{2\sigma }{\gamma }}\right]\left[\frac{\pi }{\sin \pi \left(\sigma /\lambda \right)}\right]^{2}.
\]

(d) For $$\gamma \in \{a;a=\frac{1}{2k-1},2k+1(k\in \mathbf{N})\},$$ we set
\[
H(t)=K_{\lambda }(1,t)=\min_{i\in \{1,2\}}\frac{1}{\max \{1,(t^{\gamma }\cos
\alpha _{i}+|t|^{\gamma })^{\lambda /\gamma }\}},
\]
where
$ 0<\alpha _{1}\leq \alpha _{2}<\pi ,\:\mu ,\:\sigma >0,\:\mu +\sigma =\lambda .$
Then we have the kernels%
\begin{eqnarray*}
H(xy) &=&\min_{i\in \{1,2\}}\frac{1}{\max \{1,[(xy)^{\gamma }\cos \alpha
_{i}+|xy|^{\gamma }]^{\lambda /\gamma }\}}, \\
K_{\lambda }(x,y) &=&\min_{i\in \{1,2\}}\frac{1}{\max \{|x|,(y^{\gamma
}sgn(x)\cos \alpha _{i}+|y|^{\gamma })^{\lambda /\gamma }\}}
\end{eqnarray*}%
and obtain the constant factors%
\begin{eqnarray*}
K(\sigma ) &=&K_{\lambda }(\sigma )=\int_{-\infty }^{\infty }\min_{i\in
\{1,2\}}\frac{1}{\max \{1,(t^{\gamma }\cos \alpha _{i}+|t|^{\gamma
})^{\lambda /\gamma }\}}|t|^{\sigma -1}dt \\
&=&\int_{0}^{\infty }\min_{i\in \{1,2\}}\frac{1}{\max \{1,[t^{\gamma
}(1+\cos \alpha _{i})]^{\lambda /\gamma }\}}t^{\sigma -1}dt \\
&&+\int_{0}^{\infty }\min_{i\in \{1,2\}}\frac{1}{\max \{1,[t^{\gamma
}(1-\cos \alpha _{i})]^{\lambda /\gamma }\}}t^{\sigma -1}dt \\
&=&\frac{1}{\lambda }\left[ \int_{0}^{\infty }\min_{i\in \{1,2\}}\frac{1}{%
(1+\cos \alpha _{i})^{\sigma /\gamma }}\frac{1}{\max \{1,u\}}u^{\frac{\sigma
}{\lambda }-1}du\right] \\
&&+\frac{1}{\lambda }\left[ \int_{0}^{\infty }\min_{i\in \{1,2\}}\frac{1}{%
(1-\cos \alpha _{i})^{\sigma /\gamma }}\frac{1}{\max \{1,u\}}u^{\frac{\sigma
}{\lambda }-1}du\right] \\
&=&\frac{1}{\lambda }\left[ \int_{0}^{\infty }\frac{1}{(1+\cos \alpha
_{1})^{\sigma /\gamma }}\int_{0}^{\infty }\frac{1}{\max \{1,u\}}u^{\frac{%
\sigma }{\lambda }-1}du\right] \\
&&+\frac{1}{\lambda }\left[ \int_{0}^{\infty }\frac{1}{(1-\cos \alpha
_{2})^{\sigma /\gamma }}\int_{0}^{\infty }\frac{1}{\max \{1,u\}}u^{\frac{%
\sigma }{\lambda }-1}du\right] \\
&=&\frac{1}{2^{\sigma /\gamma }}\left[\left(\sec \frac{\alpha _{1}}{2}\right)^{\frac{2\sigma
}{\gamma }}+\left(\csc \frac{\alpha _{2}}{2}\right)^{\frac{2\sigma }{\gamma }}\right]\frac{%
\lambda }{\mu \sigma }.
\end{eqnarray*}%
By (\ref{121}) and (\ref{129}), we have%
\begin{equation}
||T_{1}||=||T_{2}||=\frac{1}{2^{\sigma /\gamma }}\left[\left(\sec \frac{\alpha _{1}}{2}%
\right)^{\frac{2\sigma }{\gamma }}+\left(\csc \frac{\alpha _{2}}{2}\right)^{\frac{2\sigma }{%
\gamma }}\right]\frac{\lambda }{\mu \sigma }.  \label{147}
\end{equation}

In particular, if $\alpha _{1}=\alpha _{2}=\alpha \in (0,\pi ),$ then we
have
\[
H(t)=K_{\lambda }(1,t)=\frac{1}{\max \{1,(t^{\gamma }\cos \alpha
+|t|^{\gamma })^{\lambda /\gamma }\}},
\]
and
\[
||T_{1}||=||T_{2}||=\frac{1}{2^{\sigma /\gamma }}\left[\left(\sec \frac{\alpha }{2}\right)^{%
\frac{2\sigma }{\gamma }}+\left(\csc \frac{\alpha }{2}\right)^{\frac{2\sigma }{\gamma }%
}\right]\frac{\lambda }{\mu \sigma }.
\]

Then, for (a)-(d) we can obtain the equivalent inequalities with the kernels
and the best possible constant factors in Theorem 5-8. Setting $\delta _{\:0}=%
\frac{\sigma }{2}>0,$ we can still obtain the equivalent reverse
inequalities with the kernels and the best possible constant factors in
Theorem 5-8.

\section*{Acknowledgments}
This work is supported by the National Natural Science Foundation of China (No. 61370186), and 2013
Knowledge Construction Special Foundation Item of Guangdong Institution of
Higher Learning College and University (No. 2013KJCX0140).\\
The authors wish to express their thanks to Professors Tserendorj Batbold, Mario Krnic and Jichang Kuang for their careful reading of the manuscript and their valuable suggestions.

\end{document}